\documentclass[a4paper,12pt,twoside]{article}
\usepackage{amsmath,amssymb,ifthen}
\usepackage[amsmath,thmmarks]{ntheorem}
\usepackage{paper}
\usepackage[all]{xy}
\usepackage{mathrsfs}
\newcommand{\rig}{\mathrm{rig}}
\newcommand{\red}{\mathrm{red}}
\newcommand{\ad}{\mathrm{ad}}
\DeclareMathOperator{\Spf}{Spf}
\DeclareMathOperator{\Spa}{Spa}

\DeclareMathOperator{\Sh}{Sh}
\DeclareMathOperator{\GSp}{GSp}
\DeclareMathOperator{\PGL}{PGL}
\DeclareMathOperator{\Sp}{Sp}
\DeclareMathOperator{\SL}{SL}
\DeclareMathOperator{\Fix}{Fix}
\DeclareMathOperator{\Fil}{Fil}
\DeclareMathOperator{\vol}{vol}
\DeclareMathOperator{\Stab}{Stab}
\DeclareMathOperator{\Res}{Res}
\DeclareMathOperator{\Nilp}{\mathbf{Nilp}}
\DeclareMathOperator{\Ad}{Ad}
\DeclareMathOperator{\Irr}{Irr}
\DeclareMathOperator{\cInd}{c-Ind}
\DeclareMathOperator{\supp}{supp}
\DeclareMathOperator{\simil}{sim}
\DeclareMathOperator{\Art}{Art}
\DeclareMathOperator{\height}{height}

\newcommand{\M}{\breve{\mathscr{M}}}
\newcommand{\X}{\mathbb{X}}

\SelectTips{cm}{11}
\title{Lefschetz trace formula and $\ell$-adic cohomology of Rapoport-Zink tower for GSp(4)}
\author{Yoichi Mieda}
\begin{document}
\maketitle

\begin{firstfootnote}
 The Hakubi Center for Advanced Research / Department of Mathematics, Kyoto University, Kyoto, 606--8502, Japan

 E-mail address: \texttt{mieda@math.kyoto-u.ac.jp}

\end{firstfootnote}

\begin{abstract}
 We investigate the alternating sum of the $\ell$-adic cohomology of the Rapoport-Zink tower for GSp(4)
 by the Lefschetz trace formula. Under some assumptions on $L$-packets of GSp(4) and its inner form,
 we observe that the local Jacquet-Langlands correspondence appears in the cohomology.
\end{abstract}

\section{Introduction}
A Rapoport-Zink space is a certain moduli space of deformations by quasi-isogenies of a $p$-divisible group
with additional structures. 
By using level structures on the universal $p$-divisible group,
we can construct a projective system of \'etale coverings over the rigid generic fiber of the Rapoport-Zink space.
This projective system is called the Rapoport-Zink tower.
It can be regarded as a local analogue of a tower of Shimura varieties of PEL type.

By taking a compactly supported $\ell$-adic cohomology of the tower, we obtain a representation 
$H^i_{\mathrm{RZ}}$ of
$\mathbf{G}(\Q_p)\times\mathbf{J}(\Q_p)\times W_{\Q_p}$, where $\mathbf{G}$ is the reductive group
which is naturally attached to the local Shimura datum defining the Rapoport-Zink space, $\mathbf{J}$ is an
inner form of $\mathbf{G}$, and $W_{\Q_p}$ is the Weil group of $\Q_p$.
It is expected that the alternating sum $H_{\mathrm{RZ}}=\sum_i(-1)^iH^i_{\mathrm{RZ}}$ of $H^i_{\mathrm{RZ}}$
can be described by the (still conjectural) local Langlands correspondence of $\mathbf{G}$ and $\mathbf{J}$
(\cf \cite{MR1403942}).

The most classical examples of the Rapoport-Zink tower are the Lubin-Tate tower and the Drinfeld tower.
In these cases, the expectation above is called the non-abelian Lubin-Tate theory (\cf \cite{MR1044827})
and has been already proven (\cite{MR1464867}, \cite{MR1876802}).
There are more precise studies on the individual cohomology $H^i_{\mathrm{RZ}}$;
see \cite{MR2511742} and \cite{MR2308851}.

In this paper, we consider the case where $\mathbf{G}=\GSp_4$. In this case, the Rapoport-Zink space 
$\M$ is the moduli space of deformations by quasi-isogenies of a principally polarized 2-dimensional
$p$-divisible group with slope $1/2$. We ignore the action of the Weil group $W_{\Q_p}$ and
concentrate on the action of $\mathbf{G}(\Q_p)\times\mathbf{J}(\Q_p)$ on $H^i_{\mathrm{RZ}}$.
Our main result can be summarized as follows:

\begin{thm}[Theorem \ref{thm:main}, Corollary \ref{cor:RZ-char}]\label{thm:main-intro}
 For an irreducible smooth representation $\rho$ of $\mathbf{J}(\Q_p)$, we put
 $H_{\mathrm{RZ}}[\rho]=\sum_{i,j\ge 0}(-1)^{i+j}\Ext^j_{\mathbf{J}(\Q_p)}(H^i_{\mathrm{RZ}},\rho)^{\mathrm{sm}}$,
 where $\Ext^j_{\mathbf{J}(\Q_p)}$ is taken in the category of smooth $\mathbf{J}(\Q_p)$-representations
 and $(-)^{\mathrm{sm}}$ denotes the set of $\mathbf{G}(\Q_p)$-smooth vectors.
 Let $\phi\colon W_{\Q_p}\times \SL_2(\C)\longrightarrow \GSp_4(\C)$ be an $L$-parameter
 which is relevant for $\mathbf{J}(\Q_p)$. Assume that the $L$-packets $\Pi^{\mathbf{G}(\Q_p)}_\phi$
 and $\Pi^{\mathbf{J}(\Q_p)}_\phi$ corresponding to $\phi$ are stable and satisfy the character relation
 (see Section \ref{subsec:LLC-G-J} for notation on $L$-parameters and $L$-packets).
 
 Then, for an element of the Hecke algebra $f\in\mathcal{H}(\mathbf{G}(\Q_p))$
 supported on regular elliptic elements, we have
 \[
 \sum_{\rho\in\Pi_\phi^{\mathbf{J}(\Q_p)}}\Tr(f;H_{\mathrm{RZ}}[\rho])=-4\sum_{\pi\in\Pi_\phi^{\mathbf{G}(\Q_p)}}\Tr(f;\pi).
 \]
 Moreover, if the $\mathbf{G}(\Q_p)$-representation 
 $\Ext^j_{\mathbf{J}(\Q_p)}(H^i_{\mathrm{RZ}},\rho)^{\mathrm{sm}}$ has finite length for every $i,j\ge 0$
 and $\rho\in\Pi_\phi^{\mathbf{J}(\Q_p)}$, we have
 \[
 \sum_{\rho\in\Pi_\phi^{\mathbf{J}(\Q_p)}}\theta_{H_{\mathrm{RZ}}[\rho]}(g)=-4\sum_{\pi\in\Pi_\phi^{\mathbf{G}(\Q_p)}}\theta_{\pi}(g)
 \]
 for every regular elliptic element $g$ of $\mathbf{G}(\Q_p)$. Here $\theta_{H_{\mathrm{RZ}}[\rho]}$
 and $\theta_\pi$ denote the distribution characters of $H_{\mathrm{RZ}}[\rho]$ and $\pi$ respectively,
 which are locally constant functions over regular elements of $\mathbf{G}(\Q_p)$.
\end{thm}

Very roughly speaking, this theorem says that the local Jacquet-Langlands correspondence 
$\Pi^{\mathbf{G}(\Q_p)}_\phi\leftrightarrow \Pi^{\mathbf{J}(\Q_p)}_\phi$ appears in $H_{\mathrm{RZ}}$.

To prove the theorem above, we will apply the Lefschetz trace formula for adic spaces
developed in \cite{adicLTF}; we count fixed points on the Rapoport-Zink space under the action of
elements in $\mathbf{G}(\Q_p)\times\mathbf{J}(\Q_p)$ to compute the trace on the cohomology $H_{\mathrm{RZ}}$.
Such a method goes back to a pioneering work of Faltings \cite{MR1302321},
in which he treated the Drinfeld tower. A similar study for the Lubin-Tate tower has been carried out
by Strauch \cite{MR2383890}. Needless to say, our work is strongly inspired by these two works.
However, our case is more difficult than the classical cases in the following two points.
First, any connected component of our Rapoport-Zink space $\M$ is neither quasi-compact nor $p$-adic,
therefore harder to deal with.
This is related to the fact that neither $\mathbf{G}(\Q_p)$ nor $\mathbf{J}(\Q_p)$ is
compact modulo center (in the classical cases, $\mathbf{G}(\Q_p)$ or $\mathbf{J}(\Q_p)$ is the multiplicative
group of a division algebra).
The other point is representation-theoretic one; the local Langlands correspondences for $\mathbf{G}(\Q_p)$
and $\mathbf{J}(\Q_p)$ are not bijective. Under the ``dictionary'' between irreducible representations
and conjugacy classes, this corresponds to the fact that conjugacy and stable conjugacy are different in
$\mathbf{G}(\Q_p)$ and $\mathbf{J}(\Q_p)$.
This difference makes our argument on counting points and harmonic analysis more subtle.

We sketch the outline of this paper. In Section 2, we introduce some notation on algebraic groups and
stable orbital integrals, which will be used throughout this paper. 
In section 3, after recalling basic definitions on the Rapoport-Zink tower for $\GSp_{2d}$,
we count fixed points on the Rapoport-Zink space under the action of an element $(g,h)\in \mathbf{G}(\Q_p)\times\mathbf{J}(\Q_p)$. Our method of counting is similar to \cite[\S 2.6]{MR2383890};
we use the period map introduced in \cite[Chapter 5]{MR1393439} and the $p$-adic Hodge theory
for $p$-divisible groups. In Section 4, we construct formal models of the Rapoport-Zink spaces
with some higher levels (more precisely, levels which are open normal subgroups of parahoric subgroups
of $\mathbf{G}(\Q_p)$). Moreover, we introduce ``boundary strata'' of these formal models
and investigate group actions on them. 
These constructions are extremely important for applying the Lefschetz trace formula such as
\cite[Theorem 4.5]{adicLTF}. Basically, the content of this section
(especially Proposition \ref{prop:formal-model-disjoint}) forces us to assume that $d=2$. 
In Section 5, we construct a nice open covering of the Rapoport-Zink space with
parahoric level. The construction is similar to the case of the Drinfeld upper half space,
which has an open covering indexed by vertices of the Bruhat-Tits building for $\PGL_n$.
In Section 6, we apply the Lefschetz trace formula to a finite union of open subsets belonging to
the open covering constructed in Section 5. Finally in Section 7, we briefly review
the local Langlands correspondence for $\mathbf{G}(\Q_p)$ and $\mathbf{J}(\Q_p)$
due to Gan-Takeda \cite{MR2800725} and Gan-Tantono \cite{Gan-Tantono} respectively,
and give a proof of the main theorem. We use the harmonic-analytic method introduced in \cite{LT-LTF}.

\bigbreak

\noindent{\bfseries Acknowledgment}\quad
The author would like to thank Tetsushi Ito and Matthias Strauch for valuable discussions.
He is also grateful to Takuya Konno for helpful comments.
This work was supported by JSPS KAKENHI Grant Numbers 21740022, 24740019.

\bigbreak

\noindent{\bfseries Notation}

Let $d\ge 1$ be an integer.
For a ring $A$, let $\langle\ ,\ \rangle\colon A^{2d}\times A^{2d}\longrightarrow A$ be the symplectic pairing
defined as follows:
for $x=(x_i), y=(y_i)\in A^{2d}$,
\[
 \langle x,y\rangle=x_1y_{2d}+\cdots+x_dy_{d+1}-x_{d+1}y_d-\cdots-x_{2d}y_1.
\]
We denote by $\GSp_{2d}(A)$ the symplectic similitude group with respect to the symplectic pairing
$\langle\ ,\ \rangle$. 

For a field $k$, we denote its algebraic closure by $\overline{k}$.
Fix a prime number $p$. 
For an integer $m\ge 1$, we denote by $\Q_{p^m}$ the unique degree $m$ unramified extension of $\Q_p$,
and by $\Z_{p^m}$ the ring of integers of $\Q_{p^m}$. We denote by $\Q_{p^\infty}$ the completion of
the maximal unramified extension of $\Q_p$, and by $\Z_{p^\infty}$ the ring of integers of $\Q_{p^\infty}$.
Let $\ell$ be a prime number distinct from $p$.
We fix an isomorphism $\overline{\Q}_\ell\cong \C$ and identify them.
Every representation is considered over $\C$, and every function is $\C$-valued.

For a totally disconnected locally compact group $G$ with a fixed Haar measure,
we denote by $\mathcal{H}(G)$ the Hecke algebra of $G$, namely,
the abelian group of locally constant compactly supported functions on $G$ with convolution product.
Put $\overline{\mathcal{H}}(G)=\mathcal{H}(G)/[\mathcal{H}(G),\mathcal{H}(G)]=\mathcal{H}(G)_G$
(the $G$-coinvariant quotient).
For smooth representations $\pi_1$, $\pi_2$ of $G$, we denote by
$\Ext^i_G(\pi_1,\pi_2)$ the higher Ext group in the category of smooth $G$-representations.

\section{Notation on stable conjugacy classes}\label{sec:notation}
In this section, we introduce some basic notation on algebraic groups and harmonic analysis.
Here we will work on slightly general situation;
let $F$ be a $p$-adic field and $\mathbf{G}$ a connected reductive group over $F$.
Put $G=\mathbf{G}(F)$.
Assume for simplicity that the derived group of $\mathbf{G}$ is simply connected.
We denote by $\mathbf{Z}_{\mathbf{G}}$ the center of $\mathbf{G}$ and put $Z_G=\mathbf{Z}_{\mathbf{G}}(F)$.
For $g\in G$, let $\mathbf{Z}(g)$ denote the centralizer of $g$ and put $Z(g)=\mathbf{Z}(g)(F)$.
Since we assume that the derived group of $\mathbf{G}$ is simply connected,
$\mathbf{Z}(g)$ is connected for a semisimple $g$.
We say that $g$ is regular if $\mathbf{Z}(g)$ is a maximal torus of $\mathbf{G}$
(note that a regular element is assumed to be semisimple).
We write $G^{\mathrm{reg}}$ for the set of regular elements of $G$.
For a maximal torus $\mathbf{T}$ of $\mathbf{G}$, we put $T^{\mathrm{reg}}=\mathbf{T}(F)\cap G^{\mathrm{reg}}$.
We say that $g$ is elliptic if it is contained in an elliptic maximal torus.
If $g$ is regular, this is equivalent to saying that $\mathbf{Z}(g)$ is an elliptic maximal torus.
We write $G^{\mathrm{ell}}$ for the set of regular elliptic elements of $G$.

Two elements $g_1,g_2\in G=\mathbf{G}(F)$ is said to be stably conjugate if they are conjugate 
in $\mathbf{G}(\overline{F})$. For $g\in G$, we write $\{g\}$ (resp.\ $\{g\}_{\mathrm{st}}$) for
the conjugacy class (resp.\ stable conjugacy class) of $g$.
It is well-known that $\{g\}_{\mathrm{st}}/{\sim}$, the set of conjugacy classes in $\{g\}_{\mathrm{st}}$,
is a finite set if $g$ is regular.

Two maximal tori $\mathbf{T}_1$, $\mathbf{T}_2$ of $\mathbf{G}$ are said to be stably conjugate
if $T_1=\mathbf{T}_1(F)$ and $T_2=\mathbf{T}_2(F)$ are conjugate in $\mathbf{G}(\overline{F})$. 
For such tori $\mathbf{T}_1$, $\mathbf{T}_2$ and elements
$g_1\in T_1^{\mathrm{reg}}$, $g_2\in T_2^{\mathrm{reg}}$ which are stably conjugate,
we can construct an isomorphism
$\iota_{g_1,g_2}\colon \mathbf{T}_1\yrightarrow{\cong}\mathbf{T}_2$ as follows.
Take $h\in \mathbf{G}(\overline{F})$ such that $g_2=h^{-1}g_1 h$.
Since $g_1$ and $g_2$ are $F$-valued points, such $h$ satisfies
$h\sigma(h)^{-1}\in \mathbf{T}_1(\overline{F})$ for every $\sigma\in \Gal(\overline{F}/F)$.
By this, it is easy to see that
$\mathbf{T}_1\otimes_F\overline{F}\longrightarrow \mathbf{T}_2\otimes_F\overline{F}$;
$g\longmapsto h^{-1}gh$ descends to an isomorphism $\iota_{g_1,g_2}\colon \mathbf{T}_1\yrightarrow{\cong}\mathbf{T}_2$ over $F$. It does not depend on the choice of $h$.
In particular, stably conjugate maximal tori are isomorphic, 
and thus a maximal tori which is stably conjugate to an elliptic torus is elliptic.

For a maximal torus $\mathbf{T}$, we write $\{\mathbf{T}\}$ (resp.\ $\{\mathbf{T}\}_{\mathrm{st}}$) for
its conjugacy class (resp.\ stable conjugacy class).
We denote the set of conjugacy classes of maximal tori
(resp.\ elliptic maximal tori) of $\mathbf{G}$ by $\mathcal{T}_{\mathbf{G}}$ 
(resp.\ $\mathcal{T}^{\mathrm{ell}}_{\mathbf{G}}$), and the set of stable conjugacy classes of maximal tori
(resp.\ elliptic maximal tori) of $\mathbf{G}$ by $\mathcal{T}_{\mathbf{G},\mathrm{st}}$
(resp.\ $\mathcal{T}^{\mathrm{ell}}_{\mathbf{G},\mathrm{st}}$). 

Fix a Haar measure on $G$. 
For an element $g\in G^{\mathrm{reg}}$, we also choose a Haar measure on $Z(g)$.
Then, for each $g'\in\{g\}_{\mathrm{st}}$,
$Z(g')$ is naturally equipped with a Haar measure induced by the isomorphism
$\iota_{g,g'}\colon Z(g)\longrightarrow Z(g')$. 
For a locally constant function $f$ on $G$ whose support is compact modulo $Z_G$, we set
\[
 O_g(f)=\int_{Z(g)\backslash G}f(h^{-1}gh)dh,\quad \mathit{SO}_g(f)=\sum_{g'\in \{g\}_{\mathrm{st}}/{\sim}}O_{g'}(f),
\]
and call them the orbital integral and the stable orbital integral of $f$, respectively.
It is well-known that $O_g(f)$ always converges (\cite{MR0320232}).

\bigbreak

Next we compare stable conjugacy classes between inner forms.
Let $\mathbf{G}'$ be an inner form of $\mathbf{G}$, and fix an inner twist
$\xi\colon \mathbf{G}'\otimes_{F}\overline{F}\yrightarrow{\cong}\mathbf{G}\otimes_{F}\overline{F}$.
For $g\in G$ and $g'\in G'=\mathbf{G}'(F)$, $g$ is said to be a transfer of $g'$ with respect to $\xi$
if $g$ and $\xi(g')$ are conjugate in $\mathbf{G}(\overline{F})$.
We also say that $g$ and $g'$ match, and write $g\leftrightarrow g'$.
For $g\in G^{\mathrm{reg}}$ and $g'\in G'^{\mathrm{reg}}$ with $g\leftrightarrow g'$,
we can construct an isomorphism $\iota_{g',g}\colon \mathbf{Z}(g')\yrightarrow{\cong}\mathbf{Z}(g)$
in the same way as above.
In particular, $g$ is elliptic if and only if $g'$ is elliptic.

By \cite[Lemma 10.2]{MR858284}, there is a natural bijection
$\mathcal{T}^{\mathrm{ell}}_{\mathbf{G}',\mathrm{st}}\yrightarrow{\cong} \mathcal{T}^{\mathrm{ell}}_{\mathbf{G},\mathrm{st}}$ such that $\{\mathbf{T}'\}_{\mathrm{st}}$ corresponds to $\{\mathbf{T}\}_{\mathrm{st}}$
if and only if $\xi(\mathbf{T}'(F))$ and $\mathbf{T}(F)$ are conjugate in $\mathbf{G}(\overline{F})$.
Therefore, for every $g'\in G'^{\mathrm{ell}}$ (resp.\ $g\in G^{\mathrm{ell}}$),
we can always find $g\in G^{\mathrm{ell}}$ (resp.\ $g'\in G'^{\mathrm{ell}}$) with $g\leftrightarrow g'$.
In particular, stable conjugacy classes of regular elliptic elements of $G$ are in bijection with
those of $G'$.

The following lemma, which is used in Section \ref{sec:computation-of-char}, would be well-known. 

\begin{lem}\label{lem:Weyl-group}
 Let $\mathbf{T}$ be an elliptic maximal torus of $\mathbf{G}$, and $\mathbf{T}'$ that of $\mathbf{G}'$
 such that $\{\mathbf{T}'\}_{\mathrm{st}}$ corresponds to $\{\mathbf{T}\}_{\mathrm{st}}$ under the bijection
 above.
 Let $W_{\mathbf{T}}$ denote the Weyl group 
 $N_{\mathbf{G}}(\mathbf{T})/\mathbf{T}$ of $\mathbf{T}$, which is an algebraic group over $F$.
 Similarly we define $W_{\mathbf{T}'}$.
 Then, we have an isomorphism $W_{\mathbf{T}}\cong W_{\mathbf{T}'}$.
 In particular, $\#W_{\mathbf{T}}(F)=\#W_{\mathbf{T}'}(F)$.
\end{lem}

\begin{prf}
 Take $t\in T^{\mathrm{reg}}$ and $t'\in T'^{\mathrm{reg}}$ such that $t\leftrightarrow t'$,
 and $h\in \mathbf{G}(\overline{F})$ such that $t=h^{-1}\xi(t')h$.
 Let $\xi_h\colon \mathbf{G}'\otimes_F\overline{F}\yrightarrow{\cong} \mathbf{G}\otimes_F\overline{F}$ be
 the composite $\Ad(h^{-1})\circ \xi$. It satisfies $\xi_h(\mathbf{T}')=\mathbf{T}$, thus induces
 $W_{\mathbf{T}'}\otimes_F\overline{F}\yrightarrow{\cong} W_{\mathbf{T}}\otimes_F\overline{F}$.
 It suffices to check that this isomorphism descends to an isomorphism over $F$.
 Take $\sigma\in \Gal(\overline{F}/F)$. Since $\xi$ is an inner twist, there exists 
 $c_\sigma\in \mathbf{G}(\overline{F})$ such that $\sigma\circ\xi=\Ad(c_\sigma)\circ \xi\circ\sigma$.
 Then $t=h^{-1}\xi(t')h$ implies that $t=\Ad(\sigma(h)^{-1})\Ad(c_\sigma)\xi(t')$ (note that
 $t$ and $t'$ are rational), and thus $\sigma(h)^{-1}c_\sigma h\in \mathbf{T}(\overline{F})$.
 Therefore, for $g'\in N_{\mathbf{G}'}(\mathbf{T}')(\overline{F})$ we have
 \[
 \sigma\bigl(\xi_h(g')\bigr)=\Ad\bigl(\sigma(h)^{-1}c_\sigma h)\xi_h\bigl(\sigma(g')\bigr)\in \xi_h\bigl(\sigma(g')\bigr)\mathbf{T}(\overline{F}).
 \]
 This means that $\xi_h\colon W_{\mathbf{T}'}\otimes_F\overline{F}\yrightarrow{\cong} W_{\mathbf{T}}\otimes_F\overline{F}$ commutes with the action of $\sigma$, as desired.
\end{prf}

\section{Rapoport-Zink tower for $\boldsymbol{\GSp(2d)}$}
\subsection{Definition of the Rapoport-Zink tower}\label{subsec:def-RZ}
In this subsection, we recall basic notions on the Rapoport-Zink tower.
General definitions are given in \cite{MR1393439}, but here we restrict ourselves to the Siegel case,
namely, the case for $\GSp(2d)$.

Fix a $d$-dimensional isoclinic $p$-divisible group $\X$ over $\overline{\F}_p$ with slope $1/2$,
and a (principal) polarization $\lambda_0\colon \X\yrightarrow{\cong} \X^\vee$ of $\X$, namely,
an isomorphism satisfying $\lambda_0^\vee=-\lambda_0$.
Let $\Nilp$ be the category of $\Z_{p^\infty}$-schemes on which $p$ is locally nilpotent.
For an object $S$ of $\Nilp$, we put $\overline{S}=S\otimes_{\Z_{p^\infty}}\overline{\F}_p$.
Consider the contravariant functor $\M\colon \Nilp\longrightarrow \mathbf{Set}$ that associates $S$
with the set of isomorphism classes of pairs $(X,\rho)$ consisting of
\begin{itemize}
 \item a $d$-dimensional $p$-divisible group $X$ over $S$,
 \item and a quasi-isogeny (\cf \cite[Definition 2.8]{MR1393439})
       $\rho\colon \X\otimes_{\overline{\F}_p}\overline{S}\longrightarrow X\otimes_S\overline{S}$,
\end{itemize}
such that there exists an isomorphism $\lambda\colon X\longrightarrow X^\vee$ which makes the following
diagram commutative up to multiplication by $\Q_p^\times$:
	\[
	 \xymatrix{%
	 \mathbb{X}\otimes_{\overline{\F}_p}\overline{S}\ar[r]^-{\rho}\ar[d]^-{\lambda_0\otimes\id}& X\otimes_S\overline{S}\ar[d]^-{\lambda\otimes \id}\\
	 \mathbb{X}^\vee\otimes_{\overline{\F}_p}\overline{S}& X^\vee\otimes_S\overline{S}\lefteqn{.}\ar[l]_-{\rho^\vee}
	}
	\]
Note that such $\lambda$ is uniquely determined by $(X,\rho)$ up to multiplication by $\Z_p^\times$
and gives a polarization of $X$.
It is proved by Rapoport-Zink that $\M$ is represented by a special formal scheme (\cf \cite{MR1395723})
over $\Spf \Z_{p^\infty}$.
Moreover, $\M$ is separated over $\Spf \Z_{p^\infty}$
(\cite[Lemme 2.3.23]{MR2074714}). However, each connected component of $\M$ is neither quasi-compact nor $p$-adic.
It is known that $\dim \M^{\red}=\lfloor d^2/4\rfloor$, where $\lfloor x\rfloor$ denotes the
greatest integer less than or equal to $x$ (for example, see \cite{MR2466182}), and
every irreducible component of $\M^{\red}$ is projective over $\overline{\F}_p$
(\cite[Proposition 2.32]{MR1393439}).

Let $J$ be the group consisting of self-quasi-isogenies $h$ on $\mathbb{X}$ which makes
the following diagram commutative up to multiplication by $\Q_p^\times$:
\[
 \xymatrix{%
 \mathbb{X}\ar[r]^-{h}\ar[d]^-{\lambda_0}&\mathbb{X}\ar[d]^-{\lambda_0}\\
 \mathbb{X}^\vee&\mathbb{X}^\vee\lefteqn{.}\ar[l]_-{h^{\vee}}
 }
\]
Then, we can define a right action of $J$ on $\M$ by $h\colon \M(S)\longrightarrow \M(S)$;
$(X,\rho)\longmapsto (X,\rho\circ h)$.
It is known that $J$ is the group of $\Q_p$-valued points of an inner form $\mathbf{J}$ of
$\GSp_{2d}$ (see the next subsection). In particular, $J$ is naturally endowed with a topology.

We denote the rigid generic fiber $\M^\rig$ of $\M$ by $M$.
It is defined as $t(\M)\setminus V(p)$, where $t(\M)$ is the adic space associated with $\M$
(\cf \cite[Proposition 4.1]{MR1306024}). 
It is locally of finite type, partially proper and smooth over $\Spa (\Q_{p^\infty},\Z_{p^\infty})$
(\cite[Lemme 2.3.24]{MR2074714}).
Moreover, we know that $\dim M=d(d+1)/2$; it can be proved by using \'etaleness of
the period map (\cite[Proposition 5.17]{MR1393439}) or the $p$-adic uniformization theorem
(\cite[Theorem 6.30]{MR1393439}).

Let $\widetilde{X}$ be the universal $p$-divisible group over $\M$ and $\widetilde{X}^\rig$ the induced
$p$-divisible group over $M$. For each geometric point $x$ of $M$,
the rational Tate module $V_p(\widetilde{X}^\rig_x)$ is endowed with a non-degenerate alternating pairing
$V_p(\widetilde{X}^\rig_x)\times V_p(\widetilde{X}^\rig_x)\longrightarrow \Q_p(1)$ 
induced by a polarization on $\widetilde{X}$.
It is well-defined up to $\Z_p^\times$-multiplication.
Therefore, by taking a trivialization of the Tate twist, we get a non-degenerate symplectic form
$V_p(\widetilde{X}^\rig_x)\times V_p(\widetilde{X}^\rig_x)\longrightarrow \Q_p$
which is well-defined up to $\Q_p^\times$-multiplication.
By considering $K$-level structures on $\widetilde{X}^\rig$ for each compact open subgroup 
$K\subset K_0=\GSp_{2d}(\Z_p)$, we can construct a projective system $\{M_K\}_{K\subset K_0}$ of
finite \'etale coverings of $M$, which is called the Rapoport-Zink tower.
If $K$ is a normal subgroup of $K_0$, $M_K$ is a finite \'etale Galois covering of $M$ with
Galois group $K_0/K$. In particular, $M_{K_0}$ is nothing but $M$.
For more precise description, see \cite[5.34]{MR1393439} or \cite[\S 3.1]{RZ-GSp4}.

The group $J$ naturally acts on the projective system $\{M_K\}_{K\subset K_0}$. On the other hand,
for $g\in G=\GSp_{2d}(\Q_p)$ and a compact open subgroup $K\subset K_0$ satisfying $g^{-1}Kg\subset K_0$,
we can define a natural morphism $M_K\longrightarrow M_{g^{-1}Kg}$ over $\Q_{p^\infty}$.
Therefore, we have a right action of $G$ on the pro-object
$\sideset{\text{``}}{\text{''}}\varprojlim M_K$.

\begin{defn}
 For an integer $i$, we put
 \[
 H^i_c(M_K)=H^i_c(M_K\otimes_{\Q_{p^\infty}}\overline{\Q}_{p^\infty},\Q_\ell)\otimes_{\Q_\ell}\overline{\Q}_\ell,
 \qquad 
 H^i_{\mathrm{RZ}}=\varinjlim_{K\subset K_0}H^i_c(M_K).
 \]
 Here $H^i_c(M_K\otimes_{\Q_{p^\infty}}\overline{\Q}_{p^\infty},\Q_\ell)$ denotes the compactly supported $\ell$-adic cohomology introduced in \cite{MR1626021}.
 The group $G\times J$ naturally acts on $H^i_{\mathrm{RZ}}$. It is known that this action is smooth
 (\cite[Corollary 7.7]{MR1262943}, \cite[Corollaire 4.4.7]{MR2074714}).
\end{defn}

\begin{rem}
 We can also define a natural action of the Weil group $W_{\Q_p}$ of $\Q_p$ on $H^i_{\mathrm{RZ}}$.
 This action is expected to be very interesting, but in this article we do not consider it.
\end{rem}

\begin{defn}\label{defn:H-RZ-Ext}
 For an irreducible smooth representation $\rho$ of $J$ and integers $i,j\ge 0$, we put
 \[
 H^{i,j}_{\mathrm{RZ}}[\rho]=\Ext^j_J(H^i_{\mathrm{RZ}},\rho)^{\mathrm{sm}},
 \]
 where $(-)^{\mathrm{sm}}$ denotes the set of $G$-smooth vectors. It is a smooth representation of $G$.
\end{defn}

\begin{rem}\label{rem:RZ-rho-K}
 By the same argument as in \cite[Lemma 3.1]{LJLC}, we can show that
 $H^{i,j}_{\mathrm{RZ}}[\rho]^K=\Ext^j_J((H^i_{\mathrm{RZ}})^K,\rho)=\Ext^j_J(H^i_c(M_K),\rho)$
 for a compact open subgroup $K$ of $K_0$.
\end{rem}

\begin{lem}\label{lem:RZ-twist}
 Let $\chi\colon \Q_p^\times\longrightarrow \overline{\Q}_\ell^\times$ be an unramified character, that is,
 a character which is trivial on $\Z_p^\times$.
 Denote the composite $G\yrightarrow{\simil}\Q_p^\times\yrightarrow{\chi}\overline{\Q}_\ell^\times$ (resp.\ $J\yrightarrow{\simil}\Q_p^\times\yrightarrow{\chi}\overline{\Q}_\ell^\times$) by $\chi_G$ (resp.\ $\chi_J$) ,
 where $\simil$ denotes the similitude character.
 Then, we have $H^{i,j}_{\mathrm{RZ}}[\rho\otimes\chi_J]\cong H^{i,j}_{\mathrm{RZ}}[\rho]\otimes\chi_G$.
\end{lem}

\begin{prf}
 First let us recall the natural partition of $\M$ into open and closed formal subschemes
 introduced in \cite[3.52]{MR1393439}. For an integer $\delta\in\Z$, let $\M^{(\delta)}$ be
 the open and closed subscheme consisting of $(X,\rho)$ such that $d^{-1}\cdot\height(\rho)=\delta$.
 Note that the left hand side is always an integer. Indeed, by the definition of $\M$,
 there exist a polarization $\lambda\colon X\longrightarrow X^\vee$ and an element $a\in\Q_p^\times$
 such that $a\lambda_0=\rho^\vee\circ(\lambda\bmod p)\circ \rho$. Taking the heights of both sides,
 we obtain $d^{-1}\cdot \height(\rho)=v_p(a)\in\Z$, where $v_p$ is the $p$-adic valuation.
 Denote by $M^{(\delta)}$ the rigid generic fiber of
 $\M^{(\delta)}$. For a compact open subgroup $K$ of $K_0$, let $M_K^{(\delta)}$ be the inverse image of
 $M^{(\delta)}$ under the map $M_K\longrightarrow M$.
 We have $M_K=\coprod_{\delta\in\Z}M^{(\delta)}_K$.
 Put $H^i_{\mathrm{RZ},\delta}=\varinjlim_{K\subset K_0}H^i_c(M_K^{(\delta)})$.
 Then $H^i_{\mathrm{RZ}}=\bigoplus_{\delta\in\Z} H^i_{\mathrm{RZ},\delta}$.

 For $(g,h)\in G\times J$ with $g^{-1}Kg\subset K_0$, it is known that $(g,h)\colon M_K\longrightarrow M_{g^{-1}Kg}$
 maps $M_K^{(\delta)}$ to $M_K^{(\delta-v_p(\simil g)+v_p(\simil h))}$.
 In particular, if we denote by $(G\times J)^0$ the kernel of the homomorphism
 $G\times J\longrightarrow \Z$; $(g,h)\longmapsto v_p(\simil g)-v_p(\simil h)$,
 $H^i_{\mathrm{RZ},0}$ is a smooth representation of $(G\times J)^0$
 and $H^i_{\mathrm{RZ}}$ is isomorphic to $\cInd_{(G\times J)^0}^{G\times J}H^i_{\mathrm{RZ},0}$
 (\cf \cite[\S 4.4.2]{MR2074714}).
 Since the character $\chi_G\otimes \chi_J^{-1}\colon G\times J\longrightarrow \overline{\Q}^\times_\ell$;
 $(g,h)\longmapsto \chi_G(g)\chi_J(h)^{-1}$ is trivial on $(G\times J)^0$, we have a natural isomorphism
 $H^i_{\mathrm{RZ}}\otimes\chi_G\otimes\chi_J^{-1}\cong H^i_{\mathrm{RZ}}$ of $G\times J$-representations.
 Hence we have 
 \[
  \Ext^j_J(H^i_{\mathrm{RZ}},\rho\otimes\chi_J)\cong\Ext^j_J(H^i_{\mathrm{RZ}}\otimes\chi_J^{-1},\rho)
 \cong \Ext^j_J(H^i_{\mathrm{RZ}}\otimes\chi_G^{-1},\rho)\cong \Ext^j_J(H^i_{\mathrm{RZ}},\rho)\otimes \chi_G,
 \]
 and thus $H^{i,j}_{\mathrm{RZ}}[\rho\otimes\chi_J]\cong H^{i,j}_{\mathrm{RZ}}[\rho]\otimes\chi_G$.
\end{prf}

Sometimes it is convenient to work on the quotient $M_K/p^\Z$ of $M_K$ by the discrete subgroup $p^\Z$ of $J$.
The cohomology of $M_K/p^\Z$ and $H^{i,j}_{\mathrm{RZ}}[\rho]$ are related by the following lemma.

\begin{lem}\label{lem:mod-center}
 Assume that an irreducible smooth representation $\rho$ of $J$ is trivial on the subgroup $p^\Z\subset J$.
 Then we have $\Ext^j_J(H^i_c(M_K),\rho)\cong \Ext^j_{J/p^\Z}(H^i_c(M_K/p^\Z),\rho)$
 for each compact open subgroup $K\subset K_0$.
\end{lem}

\begin{prf}
 We will use the notation in the proof of Lemma \ref{lem:RZ-twist}.
 Since $p\in J$ maps $M_K^{(\delta)}$ isomorphically onto $M_K^{(\delta+2)}$ for every integer $\delta$,
 we have $M_K/p^\Z\cong M_K^{(0)}\amalg M_K^{(1)}$.
 Under this isomorphism, the natural morphism from $M_K$ to $M_K/p^\Z$ is described as follows:
 \begin{itemize}
  \item If $\delta=2\delta'$ is even, the restriction to $M_K^{(\delta)}$ is given by $p^{-\delta'}\colon M_K^{(\delta)}\longrightarrow M_K^{(0)}$.
  \item If $\delta=2\delta'+1$ is odd, the restriction to $M_K^{(\delta)}$ is given by $p^{-\delta'}\colon M_K^{(\delta)}\longrightarrow M_K^{(1)}$.
 \end{itemize}
 From this description we deduce that the natural push-forward map
 $H_c^i(M_K)\longrightarrow H_c^i(M_K/p^\Z)$ induces a $J$-equivariant isomorphism
 $H_c^i(M_K)_{p^\Z}\cong H_c^i(M_K/p^\Z)$.
 On the other hand, it is immediate to see that $H_c^i(M_K)=\bigoplus_{\delta\in\Z}H_c^i(M^{(\delta)}_K)$
 is a free $\overline{\Q}_p[p^\Z]$-module, 
 where $\overline{\Q}_\ell[p^\Z]$ denotes the group algebra of $p^\Z$.
 In particular, $H_c^i(M_K)$ is acyclic for $(-)_{p^\Z}$ (namely, the higher left derived functor
 of $(-)_{p^\Z}$ vanishes).
 Therefore we have
 \[
  \Ext^j_J\bigl(H_c^i(M_K),\rho\bigr)=\Ext^j_{J/p^\Z}\bigl(H_c^i(M_K)_{p^\Z},\rho\bigr)
 \cong \Ext^j_{J/p^\Z}\bigl(H_c^i(M_K/p^\Z),\rho\bigr),
 \]
 as desired.
 \end{prf}

\begin{cor}\label{cor:RZ-adm-vanish}
 Let $\rho$ be an irreducible smooth representation of $J$.
 \begin{enumerate}
  \item For integers $i,j\ge 0$, $H^{i,j}_{\mathrm{RZ}}[\rho]$ is an admissible representation of
	$G$.
  \item If $j>d-1$, we have $H^{i,j}_{\mathrm{RZ}}[\rho]=0$.
 \end{enumerate}
\end{cor}

\begin{prf}
 First assume that $\rho$ is trivial on $p^\Z\subset J$. Then, for a compact open subgroup $K\subset K_0$, we have
 \[
  H^{i,j}_{\mathrm{RZ}}[\rho]^K=\Ext^j_{J/p^\Z}\bigl(H^i_c(M_K/p^\Z),\rho\bigr)
 \]
 by Remark \ref{rem:RZ-rho-K} and Lemma \ref{lem:mod-center}. As in \cite[Proposition 4.4.13]{MR2074714},
 $H^i_c(M_K/p^\Z)$ is a finitely generated
 $J/p^\Z$-module, and thus \cite[Corollary II.3.2]{MR1471867} tells us that
 $\Ext^j_{J/p^\Z}(H^i_c(M_K/p^\Z),\rho)$ is finite-dimensional and vanishes for $j>d-1$
 (here $d-1$ is the split semisimple rank of $\mathbf{J}$).
 Since $H^{i,j}_{\mathrm{RZ}}[\rho]=\varinjlim_{K\subset K_0}H^{i,j}_{\mathrm{RZ}}[\rho]^K$,
 we obtain i) and ii) for this case.

 Next we consider a general $\rho$. Let $\omega\colon \Q_p^\times\longrightarrow \overline{\Q}^\times_\ell$
 be the central character of $\rho$.
  Take $c\in\overline{\Q}^\times_\ell$ such that $c^2=\omega(p)$, and $\chi\colon \Q_p^\times\longrightarrow \overline{\Q}_\ell^\times$ the character given by $\chi(a)=c^{-v_p(a)}$.
 Lemma \ref{lem:RZ-twist} tells us that
 $H^{i,j}_{\mathrm{RZ}}[\rho]=H^{i,j}_{\mathrm{RZ}}[\rho\otimes\chi_J]\otimes \chi_G^{-1}$.
 Since $\rho\otimes \chi_J$ is trivial on $p^\Z$, the right hand side is admissible and vanishes
 for $j>d-1$. This concludes the proof.
\end{prf}

By the corollary above, we can take the alternating sum of $H^{i,j}_{\mathrm{RZ}}[\rho]$.

\begin{defn}\label{defn:H-RZ-alternating}
 For an irreducible smooth representation $\rho$ of $J$, we put 
 $H_{\mathrm{RZ}}[\rho]=\sum_{i,j\ge 0}(-1)^{i+j}H^{i,j}_{\mathrm{RZ}}[\rho]$, where the sum is taken
 in the Grothendieck group of admissible representations of $G$.
\end{defn}
The goal of this paper is to investigate $H_{\mathrm{RZ}}[\rho]$ by means of the Lefschetz trace formula.

In the sequel, we fix Haar measures on $G$ and $J$. For each $g\in G^{\mathrm{reg}}$,
we also fix a Haar measure on $Z(g)$. If $g$ is elliptic,
then we normalize the measure so that $\vol(Z(g)/p^\Z)=1$,
where $p^\Z\subset G$ is endowed with the counting measure.
Note that if $g_1,g_2\in G^{\mathrm{ell}}$ are stably conjugate, then the isomorphism
$\iota_{g_1,g_2}\colon Z(g_1)\yrightarrow{\cong} Z(g_2)$ preserves the measures.
For $g\in G^{\mathrm{ell}}$ and a locally constant function $f$ on $G$ whose support is compact modulo $Z_G$,
we have 
\[
 O_g(f)=\int_{G/p^\Z}f(h^{-1}gh)dh,\quad
 \mathit{SO}_g(f)=\sum_{g'\in \{g\}_{\mathrm{st}}/{\sim}}\int_{G/p^\Z}f(h^{-1}g'h)dh.
\]
Similarly we fix a Haar measure of the centralizer of
each regular element of $J$. For $g\in G^{\mathrm{ell}}$ and $h\in J^{\mathrm{ell}}$ with $g\leftrightarrow h$,
the isomorphism
$\iota_{h,g}\colon Z(h)\yrightarrow{\cong} Z(g)$ preserves the measures.

\subsection{Period space and period map}
The goal of the rest of this section is to count fixed points under the action of
$(g,h)\in G\times J$ on $M_K/p^\Z$.
As in \cite[\S 2.6]{MR2383890}, we use the period map introduced in \cite[Chapter 5]{MR1393439}.

Put $L_0=\Frac W(\overline{\F}_p)$ and denote the Frobenius automorphism on $L_0$ by $\sigma$.
Although $L_0$ is isomorphic to $\Q_{p^\infty}$, we distinguish them as in \cite{MR1393439}.
An isocrystal over $\overline{\F}_p$ is a finite-dimensional $L_0$-vector space
equipped with a bijective $\sigma$-linear endomorphism
(\cf \cite[\S 1.1]{MR1393439}).

Let $D(\X)_\Q=(N,\Phi)$ be the rational Dieudonn\'e module of $\X$, which is a $d$-dimensional isocrystal
over $\overline{\F}_p$. 
The fixed polarization $\lambda_0$ on $\X$ gives the alternating pairing
$\psi\colon N\times N\longrightarrow L_0$ 
satisfying $\psi(\Phi(x),\Phi(y))=p\sigma(\psi(x,y))$ for every $x,y\in N$.
We define the algebraic group $\mathbf{J}$ over $\Q_p$ as follows: for a $\Q_p$-algebra $R$, 
the group $\mathbf{J}(R)$ consists of elements $g\in \Aut_{R\otimes_{\Q_p}L_0}(R\otimes_{\Q_p}N)$ such that
\begin{itemize}
 \item $g$ commutes with $\Phi$, i.e., $g\circ (\id_R\otimes \Phi)=(\id_R\otimes \Phi)\circ g$,
 \item and $g$ preserves the pairing $\psi$ up to scalar multiplication,
       i.e., there exists $c(g)\in (R\otimes_{\Q_p}L_0)^\times$ such that 
       $\psi(gx,gy)=c(g)\psi(x,y)$ for every $x,y\in R\otimes_{\Q_p}N$.
\end{itemize}
Representability of $\mathbf{J}$ is shown in \cite[Proposition 1.12]{MR1393439}.
By the Dieudonn\'e theory, we have $\mathbf{J}(\Q_p)=J$.

Since the isocrystal $(N,\Phi)$ is basic, $\mathbf{J}$ is known to be an inner form of $\GSp_{2d}$
(\cite[Corollary 1.14, Remark 1.15]{MR1393439}). For later use, we will observe it directly.
Put $N^\circ=N^{p^{-1}\Phi^2}$. It is a $\Phi$-stable $\Q_{p^2}$-subspace of $N$ satisfying
$L_0\otimes_{\Q_{p^2}}N^\circ=N$.
For $x,y\in N^\circ$, we have $\sigma^2(\psi(x,y))=p^{-2}\psi(\Phi^2(x),\Phi^2(y))=p^{-2}\psi(px,py)=\psi(x,y)$,
and thus $\psi(x,y)\in \Q_{p^2}$. 
Therefore $\psi$ gives a perfect alternating bilinear pairing 
$\psi\colon N^\circ\times N^\circ\longrightarrow \Q_{p^2}$.
Its base change from $\Q_{p^2}$ to $L_0$ coincides with the original $\psi$.
By using $N^\circ$ and the restrictions of $\Phi$ and $\psi$ on it, we can describe $\mathbf{J}$ as follows:
for each $\Q_p$-algebra $R$, 
\[
 \mathbf{J}(R)=\bigl\{g\in \Aut_{R\otimes_{\Q_p}\Q_{p^2}}(R\otimes_{\Q_p}N^\circ)\bigm| \text{$g$ satisfies the similar conditions as above}\bigr\}.
\]
In the sequel, we always use this description of $\mathbf{J}$.
Now we can prove the following:

\begin{lem}\label{lem:J-inner-form}
 We have a natural isomorphism $\xi\colon \mathbf{J}\otimes_{\Q_p}\Q_{p^2}\yrightarrow{\cong}\GSp(N^{\circ},\psi)$
 of algebraic groups over $\Q_{p^2}$.
\end{lem}

\begin{prf}
 Take a $\Q_{p^2}$-algebra $R$. Then we have 
 $R\otimes_{\Q_p}N^\circ\cong (R\otimes_{\Q_{p^2}}N^\circ)\oplus (R\otimes_{\Q_{p^2}}{}^\sigma\!N^\circ)$,
 where ${}^\sigma\!N^\circ$ is the scalar extension of $N^\circ$ by
 $\sigma\colon \Q_{p^2}\longrightarrow \Q_{p^2}$. Under this isomorphism,
 $\id_R\otimes\Phi\colon R\otimes_{\Q_p}N^\circ\longrightarrow R\otimes_{\Q_p}N^\circ$ is expressed by
 the matrix
 \[
  \begin{pmatrix}0&\id_R\otimes \Phi_1\\\id_R\otimes \Phi_2&0\end{pmatrix},
 \]
 where $\Phi_1$ (resp.\ $\Phi_2$) denotes the $\Q_{p^2}$-homomorphism
 ${}^\sigma\!N^\circ\longrightarrow N^\circ$ (resp.\ $N^\circ\longrightarrow {}^\sigma\!N^\circ$)
 induced by $\Phi$. Therefore, every element $g\in \Aut_{R\otimes_{\Q_p}\Q_{p^2}}(R\otimes_{\Q_p}N^\circ)$
 can be written as $g'\oplus g''$ with $g'\in \Aut_R(R\otimes_{\Q_{p^2}}N^\circ)$ and 
 $g''\in \Aut_R(R\otimes_{\Q_{p^2}}{}^\sigma\!N^\circ)$,
 and the condition $g\circ (\id_R\otimes \Phi)=(\id_R\otimes \Phi)\circ g$ is equivalent to
 $g'\circ (\id_R\otimes \Phi_1)=(\id_R\otimes \Phi_1)\circ g''$ and 
 $g''\circ (\id_R\otimes \Phi_2)=(\id_R\otimes \Phi_2)\circ g'$.
 For $g'\in \Aut_{R\otimes_{\Q_p}\Q_{p^2}}(R\otimes_{\Q_p}N^\circ)$, put 
 $g''=(\id_R\otimes \Phi_1)^{-1}\circ g'\circ (\id_R\otimes \Phi_1)$. Then the pair $(g',g'')$ satisfies
 the conditions above (note that $\Phi_1\circ \Phi_2=p$).
 In other words, the group
 \[
 \bigl\{g\in \Aut_{R\otimes_{\Q_p}\Q_{p^2}}(R\otimes_{\Q_p}N^\circ)\bigm| g\circ (\id_R\otimes \Phi)=(\id_R\otimes \Phi)\circ g\bigr\}
 \]
 can be identified with the group $\Aut_R(R\otimes_{\Q_{p^2}}N^\circ)$.
 Now it is straightforward to see that $(g',g'')$ preserves the pairing $\psi$ on $R\otimes_{\Q_p}N^\circ$
 up to scalar if and only if $g\in \GSp(N^\circ,\psi)(R)$. This concludes the proof.
\end{prf}

By the construction of the isomorphism $\xi$, we have the following:

\begin{cor}\label{cor:J-inner-form}
 A natural homomorphism $\mathbf{J}\longrightarrow \Res_{\Q_{p^2}/\Q_p}\GSp(N^\circ,\psi)$
 corresponds to $\xi$ by the adjointness between base change and the Weil restriction.
 In particular, the composite 
 $\mathbf{J}(\Q_p)\hooklongrightarrow \mathbf{J}(\Q_{p^2})\yrightarrow[\cong]{\xi}\GSp(N^{\circ},\psi)$
 is nothing but the natural inclusion.
\end{cor}

\begin{rem}\label{rem:GU(d,D)}
 Actually, we can describe $\mathbf{J}$ more explicitly as follows.

 Let $\Sigma_2$ be a (unique) one-dimensional
 $p$-divisible group with slope $1/2$ over $\overline{\F}_p$. It is well-known that there exists
 a polarization $\lambda_{\Sigma_2}$ on $\Sigma_2$;
 for example, a principal polarization on a supersingular elliptic curve over $\overline{\F}_p$ induces
 such a polarization. Put $D=\End(\Sigma_2)\otimes_{\Z_p}\Q_p$. Then $D$ is a quaternion division algebra
 over $\Q_p$ and $\lambda_{\Sigma_2}$ induces an involution on it.
 By \cite[Lemma 4.1]{RZ-GSp4}, we know that $(\X,\lambda_0)$ and $(\Sigma_2^{\oplus d},\lambda_{\Sigma_2}^{\oplus d})$ are isogenous. Therefore, we can prove without difficulty that the algebraic group $\mathbf{J}$ is 
 isomorphic to the quaternionic unitary similitude group $\mathrm{GU}(d,D)$.
\end{rem}

Next we introduce the period space for $\GSp_{2d}$. 

\begin{defn}
 \begin{enumerate}
  \item Let $\boldsymbol{\mathcal{F}}$ be the Grassmannian over $L_0$ parameterizing $d$-dimensional
	subspaces $\Fil\subset N$ such that $\Fil^\perp=\Fil$.
  \item Let $L$ be a finite extension of $L_0$.
	An element $\Fil\subset L\otimes_{L_0}N$ of $\boldsymbol{\mathcal{F}}(L)$ is said to be weakly admissible
	if, for every subspace $N'$ of $N$ which is stable under $\Phi$, the following inequality holds:
	\[
	\dim_L\bigl((L\otimes_{L_0}N')\cap \Fil\bigr)\le \frac{1}{2}\dim N'.
	\]
	It is known that there exists a canonical open rigid subspace 
	$\Omega\subset \boldsymbol{\mathcal{F}}^{\ad}$
	such that $\Omega(L)=\{\Fil\in \boldsymbol{\mathcal{F}}(L)\mid \text{$\Fil$ is weakly admissible}\}$
	for every finite extension $L$ of $L_0$ (\cite[Proposition 1.36]{MR1393439}).
	We call this $\Omega$ a period space for $\GSp_{2d}$.
	The group $\GSp(N,\psi)$ naturally acts on $\boldsymbol{\mathcal{F}}$ and
	the induced action of $J=\mathbf{J}(\Q_p)\subset \GSp(N,\psi)$ preserves $\Omega\subset \boldsymbol{\mathcal{F}}$.
 \end{enumerate}
\end{defn}

The following theorem is due to Rapoport-Zink:

\begin{thm}\label{thm:RZ-period}
 \begin{enumerate}
  \item There exists a $J$-equivariant \'etale morphism
	$\wp\colon M\longrightarrow \Omega$ over $L_0$ called the period morphism.
	For a finite extension $L$ of $L_0$ and an $L$-valued point $x=(X,\rho)$ of $M$,
	$\wp(x)$ is given by the subspace $\rho_*^{-1}(\Fil_X)$ of $L\otimes_{L_0}N$,
	where $\rho_*\colon D(\X)_{\Q}\yrightarrow{\cong}D(X)_{\Q}$ is the isomorphism between rational
	Dieudonn\'e modules induced by $\rho$, and $\Fil_X\subset L\otimes_{L_0}D(X)_{\Q}$ is 
	the Hodge filtration of $X$.
  \item The period map $\wp$ induces a surjection on classical points. Namely, for every finite extension $L$
	of $L_0$ and every $L$-valued point $x$ of $\Omega$, there exist a finite extension $L'$ of $L$ and
	an $L'$-valued point $\widetilde{x}$ of $M$ such that $\wp(\widetilde{x})=x$.
 \end{enumerate}
\end{thm}

\begin{prf}
 The period map $\wp$ is constructed in \cite[5.16]{MR1393439}. Precisely speaking, our $\wp$ is
 the first factor $\breve{\pi}_1$ of the period map $\breve{\pi}$ defined by Rapoport-Zink.

 ii) follows from \cite[Proposition 5.28]{MR1393439}; note that Fontaine's conjecture assumed in
 the proposition has been solved by Kisin (\cite[Corollary 2.2.6]{MR2263197}).
\end{prf}

The following proposition is the first step of our point counting:

\begin{prop}\label{prop:point-counting-period}
 Let $h$ be a regular element of $J$.
 Then all fixed points of $\boldsymbol{\mathcal{F}}$ under $h$ are discrete with multiplicity one.
 If moreover $h$ is elliptic, then every fixed point lies in $\Omega$.
\end{prop}

The former part is well-known. In order to see the latter part, we will use the theory of 
Harder-Narasimhan filtrations. Let us fix a finite extension $L$ of $L_0$ and
an element $\Fil\in \boldsymbol{\mathcal{F}}(L)$.
For a non-zero subspace $N'$ of $N$ which is stable under $\Phi$, we put
\[
 \mu(N')=\frac{\dim_L\bigl((L\otimes_{L_0}N')\cap \Fil\bigr)-1/2\dim_{L_0} N'}{\dim_{L_0}N'}.
\]
We say that $N'\neq 0$ is semi-stable if every non-zero $\Phi$-stable subspace $N''\subset N'$
satisfies $\mu(N'')\le \mu(N')$.

The following proposition is a part of \cite[Proposition 1.4]{MR1393439}:

\begin{prop}\label{prop:HN}
 There exists a unique $\Phi$-stable subspace $N_0\subset N$ satisfying the following conditions:
 \begin{itemize}
  \item $N_0$ is semi-stable (in particular non-zero).
  \item For every $\Phi$-stable $N'$ with $N_0\subsetneq N'\subset N$, we have $\mu(N_0)>\mu(N')$.
 \end{itemize}
\end{prop}

\begin{prf}[of Proposition \ref{prop:point-counting-period}]
 Let $h\in J$ be a regular element, $L$ a finite extension of $L_0$
 and $\Fil$ an element of $\boldsymbol{\mathcal{F}}(L)$ which is fixed by $h$.
 We will assume $\Fil\notin \Omega$ and prove that $h$ is not elliptic.
 Since $\Fil\notin \Omega$, $\mu(N')>0=\mu(N)$ for some $\Phi$-stable subspace $N'\subset N$.
 Therefore $N$ is not semi-stable and thus $N_0\subsetneq N$, where $N_0$ is given in Proposition \ref{prop:HN}.
 Since $\Fil$ is fixed by $h$, we have $hN_0=N_0$ by the uniqueness.
 
 Let us prove that $N_0\subset N_0^\perp$. The following argument is inspired by
 \cite[Proposition 1.43]{MR1393439}.
 Put $W=N_0\cap N_0^\perp$ and assume that $W\subsetneq N_0$. Then $\psi$ induces a perfect alternating
 pairing $N_0/W\times N_0/W\longrightarrow L_0$. Denote the image of 
 $(L\otimes_{L_0}N_0)\cap \Fil$ under $L\otimes_{L_0}N_0\longrightarrow (L\otimes_{L_0}N_0)/(L\otimes_{L_0}W)$
 by $\Fil'$. Since $\Fil^\perp=\Fil$, we have $\Fil'\subset \Fil'^\perp$.
 Therefore we have $\dim_L\Fil'\le 1/2\dim_{L_0}(N_0/W)$. 
 On the other hand, by the definition of $N_0$, $N_0\subsetneq N$ implies $\mu(N_0)>\mu(N)=0$.
 Hence, if $W\neq 0$,
 \begin{align*}
 \mu(W)&=\frac{\dim_L\bigl((L\otimes_{L_0}W)\cap \Fil\bigr)}{\dim_{L_0}W}-\frac{1}{2}
 =\frac{\dim_L\bigl((L\otimes_{L_0}N_0)\cap \Fil\bigr)-\dim_L\Fil'}{\dim_{L_0}W}-\frac{1}{2}\\
  &\ge \frac{\dim_L\bigl((L\otimes_{L_0}N_0)\cap \Fil\bigr)-1/2\dim_{L_0}(N_0/W)}{\dim_{L_0}W}-\frac{1}{2}
  =\frac{\dim_{L_0}N_0}{\dim_{L_0}W}\mu(N_0)\\
  &>\mu(N_0),
 \end{align*}
 which contradicts to semi-stability of $N_0$. If $W=0$,
 \[
 \dim_L\bigl((L\otimes_{L_0}N_0)\cap \Fil\bigr)-\frac{1}{2}\dim_{L_0}N_0
 =\dim_L\Fil'-\frac{1}{2}\dim_{L_0}(N_0/W)\le 0,
 \]
 which contradicts to $\mu(N_0)>0$. Thus we get $N_0\subset N_0^\perp$.

 Put $N_0^\circ=(N_0)^{p^{-1}\Phi^2}\subset N^\circ$.
 Since $(N_0,\Phi)$ is isoclinic of slope $1/2$, we have $L_0\otimes_{\Q_{p^2}}N_0^\circ=N_0$.
 In particular $0\subsetneq N_0^\circ\subsetneq N^\circ$.
 Moreover we have $N^\circ_0\subset (N^\circ_0)^\perp$, since $N_0\subset N_0^\perp$.
 Consider the subgroup $\mathbf{P}$ of $\mathbf{J}$ given as follows ($R$ denotes a $\Q_p$-algebra):
 \[
  \mathbf{P}(R)=\bigl\{h'\in \mathbf{J}(R)\bigm| h'(R\otimes_{\Q_p}N_0^\circ)=R\otimes_{\Q_p}N_0^\circ\bigr\}.
 \]
 Here $\mathbf{J}(R)$ is regarded as a subgroup of $\Aut_{R\otimes_{\Q_p}\Q_{p^2}}(R\otimes_{\Q_p}N^\circ)$.
 In the similar way as in the proof of Lemma \ref{lem:J-inner-form},
 we can see that $\mathbf{P}\otimes_{\Q_p}\Q_{p^2}$ is isomorphic to the stabilizer subgroup 
 $\Stab_{\GSp(N^\circ,\psi)}(N_0^\circ)$ of $N_0^\circ$ in $\GSp(N^\circ,\psi)$.
 Therefore $\mathbf{P}\otimes_{\Q_p}\Q_{p^2}$ is a proper parabolic subgroup
 of $\mathbf{J}\otimes_{\Q_p}\Q_{p^2}$, 
 and thus $\mathbf{P}$ is a proper parabolic subgroup of $\mathbf{J}$.
 Since $h\in \mathbf{P}(\Q_p)$, the following lemma says that $h$ is not elliptic.
 This completes the proof.
\end{prf}

\begin{lem}\label{lem:proper-parab-non-elliptic}
 Let $F$ be a $p$-adic field, $\mathbf{G}$ a connected reductive group over $F$ and
 $g$ a regular elliptic element of $\mathbf{G}(F)$.
 Then, for every proper parabolic subgroup $\mathbf{P}$ of $\mathbf{G}$ defined over $F$,
 $g$ does not lie in $\mathbf{P}(F)$.
\end{lem}

\begin{prf}
 Assume that there exists a proper parabolic subgroup of $\mathbf{G}$ defined over $F$ such that
 $g\in \mathbf{P}(F)$.
 Then, since $g\in \mathbf{G}(F)$ is semisimple, there exists a Levi subgroup $\mathbf{L}$ of $\mathbf{P}$
 defined over $F$ such that
 $g\in \mathbf{L}(F)$ (\cf \cite[13.3.8 (i), 8.4.4, 16.1.4]{MR1642713}).
 By the restricted root decomposition, it is easy to see that the split center of $\mathbf{L}$ is
 strictly bigger than that of $\mathbf{G}$.
 In other words, the center $\mathbf{Z}_{\mathbf{L}}$ of $\mathbf{L}$ is not anisotropic
 modulo $\mathbf{Z}_{\mathbf{G}}$.
 Therefore the centralizer of $g$, that contains $\mathbf{Z}_{\mathbf{L}}$,
 is not anisotropic modulo $\mathbf{Z}_{\mathbf{G}}$.
 This contradicts to the assumption that $g$ is regular elliptic.
\end{prf}

\begin{rem}\label{rem:fixed-point-number-period}
 By the Bruhat decomposition, we can easily calculate the number of fixed points in 
 Proposition \ref{prop:point-counting-period}; the number is $d!2^{d-1}$.
\end{rem}

\subsection{Counting fixed points under the group action}
Let $(g,h)$ be an element of $G\times J$ and $K$ a compact open subgroup of $K_0$
normalized by $g$. Let $\wp_K\colon M_K\longrightarrow \Omega$ and
$\wp_{K,p}\colon M_K/p^\Z\longrightarrow \Omega$ be the \'etale morphisms induced
from the period map $\wp$.
In Proposition \ref{prop:point-counting-period}, we considered fixed points on the period space $\Omega$.
Thus, to count fixed points on $M_K/p^\Z$, it suffices to investigate the action of $(g,h)$
on the fiber $\wp^{-1}_{K,p}(x)$ of each point $x$ in $\Omega$ fixed by $h$.

\begin{defn}\label{defn:transfer}
 Let $L$ be a finite extension of $L_0$ and $x\in \Omega(L)$ a point fixed by $h$.
 We denote the subspace of $L\otimes_{L_0}N$ corresponding to $x$ by $\Fil_x$.
 Then, since $(N,\Phi,\Fil_x)$ is a weakly admissible filtered isocrystal, there exists
 a $2d$-dimensional $p$-adic Galois representation $V_x$ of $\Gal(\overline{L}/L)$ such that 
 $D_{\mathrm{crys}}(V_x)\cong (N,\Phi,\Fil_x)$ (\cf \cite{MR1779803}).
 As the functor $D_{\mathrm{crys}}$ is fully faithful and compatible with tensor products,
 the alternating bilinear pairing $\psi\colon N\times N\longrightarrow L_0$ induces
 an alternating bilinear pairing $\psi_x\colon V_x\times V_x\longrightarrow \Q_p(1)$.
 Since $h\colon N\yrightarrow{\cong} N$ commutes with $\Phi$, preserves $\Fil_x$ and preserves $\psi$ up to
 $\Q_p^\times$-multiplication, it induces a $\Gal(\overline{L}/L)$-automorphism $g_{h,x}$
 on $V_x$ preserving $\psi_x$
 up to $\Q_p^\times$-multiplication. 
 By choosing isomorphisms $\Q_p(1)\cong \Q_p$ and $(V_x,\psi_x)\cong (\Q_p^{2d},\langle\ ,\ \rangle)$,
 $g_{h,x}$ can be regarded as an element of $G$. Obviously, the conjugacy class of $g_{h,x}$
 is independent of the choice of the isomorphisms above.
\end{defn}

\begin{prop}\label{prop:stable-conj}
The element $g_{h,x}\in G$ is a transfer of $h\in J$ with respect to $\xi$.
 Namely, if we fix an isomorphism
 $(N^\circ,\psi)\cong (\Q_{p^2}^{2d},\langle\ ,\ \rangle)$, the image of $h$ under the composite
 \[
 \mathbf{J}(\Q_p)\hooklongrightarrow \mathbf{J}(\Q_{p^2})\yrightarrow[\cong]{\xi} \GSp(N^\circ,\psi)
 \cong \GSp_{2d}(\Q_{p^2})
 \]	
 and $g_{h,x}\in G$ are conjugate in $\GSp_{2d}(\overline{\Q}_p)$.

 In particular, $g_{h,x}\in G$ is regular (resp.\ regular elliptic) if and only if $h$ is regular
 (resp.\ regular elliptic).
\end{prop}

\begin{prf}
 Since $D_{\mathrm{crys}}(V_x)\cong (N,\Phi,\Fil_x)$, we have an isomorphism
 \[
 V_x\otimes_{\Q_p}B_{\mathrm{dR}}\cong N\otimes_{L_0}B_{\mathrm{dR}}=N^\circ\otimes_{\Q_{p^2}}B_{\mathrm{dR}}.
 \]
 Since we have fixed isomorphisms $(V_x,\psi_x)\cong (\Q_p^{2d},\langle\ ,\ \rangle)$ and
 $(N^\circ,\psi)\cong (\Q_{p^2}^{2d},\langle\ ,\ \rangle)$, the isomorphism above induces
 \[
 B_{\mathrm{dR}}^{2d}\yrightarrow{\cong}V_x\otimes_{\Q_p}B_{\mathrm{dR}}\yrightarrow{\cong}N^\circ\otimes_{\Q_{p^2}}B_{\mathrm{dR}}\yrightarrow{\cong}B_{\mathrm{dR}}^{2d},
 \]
 which gives an element $\alpha\in \GSp_{2d}(B_{\mathrm{dR}})$.
 By the definition of $g_{h,x}$ and Corollary \ref{cor:J-inner-form}, we have
 $g_{h,x}=\alpha^{-1}\xi(h)\alpha$ in $\GSp_{2d}(B_{\mathrm{dR}})$. Therefore
 $g_{h,x}$ and $\xi(h)$ are conjugate in $\GSp_{2d}(B_{\mathrm{dR}})$, and thus
 conjugate in $\GSp_{2d}(\overline{\Q}_p)$.
\end{prf}

We will consider how the conjugacy class of $g_{h,x}$ changes
when we vary $h$ inside its stable conjugacy class. Assume that $h$ is regular.
Let $\mathbf{T}_h=\mathbf{Z}(h)$ be the centralizer of $h$ in $\mathbf{J}$,
which is a maximal torus of $\mathbf{J}$.

\begin{lem}\label{lem:torus-stabilizer}
 For every $\Q_p$-algebra $R$ and an element $t\in \mathbf{T}_h(R)$, $t$ stabilizes the subspace
 $R\otimes_{\Q_p}\Fil_x\subset R\otimes_{\Q_p}(L\otimes_{L_0}N)=L\otimes_{L_0}(R\otimes_{\Q_p}N)$.
\end{lem}

\begin{prf}
 Let $\boldsymbol{\mathcal{F}}^\circ$ be the Grassmannian over $\Q_{p^2}$ parameterizing $d$-dimensional
 subspaces $\Fil\subset N^\circ$ such that $\Fil^\perp=\Fil$. 
 Then, clearly $\boldsymbol{\mathcal{F}}^\circ\otimes_{\Q_{p^2}}L_0=\boldsymbol{\mathcal{F}}$.
 Moreover, since $h$ is regular, the fixed point $x$ in $\boldsymbol{\mathcal{F}}$ comes from a closed point
 of $\boldsymbol{\mathcal{F}}^\circ$; namely, there exists a finite extension $F$ of $\Q_{p^2}$ contained in $L$
 and $x'\in \boldsymbol{\mathcal{F}}^\circ(F)$ such that $\Fil_x=L\otimes_{F}\Fil_{x'}$.
 Therefore, we may consider the subspace $\Fil_{x'}\subset F\otimes_{\Q_{p^2}}N^\circ$ instead of
 $\Fil_x\subset L\otimes_{L_0}N$.
 
 Put $\mathbf{G}=\GSp(N^\circ,\psi)\otimes_{\Q_{p^2}}F$. This algebraic group acts on
 $\boldsymbol{\mathcal{F}}_{\!\!F}^\circ=\boldsymbol{\mathcal{F}}^\circ\otimes_{\Q_{p^2}}F$.
 Let $\mathbf{P}$ be the stabilizer of $x'\in \boldsymbol{\mathcal{F}}_{\!\!F}^\circ(F)$ in $\mathbf{G}$.
 It is a parabolic subgroup of $\mathbf{G}$.
 We have homomorphisms
 $\mathbf{J}\longrightarrow \Res_{\Q_{p^2}/\Q_p}\GSp(N^\circ,\psi)\longrightarrow \Res_{F/\Q_p}\mathbf{G}$,
 which induce an action of $\mathbf{J}$ on $\Res_{F/\Q_p}\boldsymbol{\mathcal{F}}_{\!\!F}^\circ$.
 What we would like to show is that the subgroup $\mathbf{T}_h$ of $\mathbf{J}$
 stabilizes $x'\in (\Res_{F/\Q_p}\boldsymbol{\mathcal{F}}_{\!\!F}^\circ)(\Q_p)$;
 in other words, $\mathbf{T}_h\subset \Res_{F/\Q_p}\mathbf{P}$.
 Let us denote by $h'$ the image of $h\in\mathbf{J}(\Q_p)$ in $(\Res_{F/\Q_p}\mathbf{G})(\Q_p)=\mathbf{G}(F)$.
 Note that $h'\in\mathbf{P}(F)$, for $h$ stabilizes $x'$.
 Moreover $h'$ is regular, since it is the image of $\xi(h)\in \GSp(N^\circ,\psi)$
 under $\GSp(N^\circ,\psi)\longrightarrow \GSp(N^\circ\otimes_{\Q_{p^2}}F,\psi)=\mathbf{G}(F)$.
 It suffices to show that the centralizer $\mathbf{S}_{h'}$ of $h'$ in $\mathbf{G}$ is contained in $\mathbf{P}$.
 We can pass to an algebraic closure $\overline{F}$ of $F$; we simply write $\mathbf{G}$ and $\mathbf{P}$
 for $\mathbf{G}\otimes_F\overline{F}$ and $\mathbf{P}\otimes_F\overline{F}$, respectively.
 Take a maximal torus $\mathbf{T}'$ of $\mathbf{P}$ containing $h'$.
 As $\mathbf{P}$ is a parabolic subgroup of $\mathbf{G}$,
 it contains a maximal torus $\mathbf{T}''$ of $\mathbf{G}$.
 Since $\mathbf{T}'$ and $\mathbf{T}''$ are conjugate in $\mathbf{P}$,
 $\mathbf{T}'$ is also a maximal torus of $\mathbf{G}$.
 This implies that $\mathbf{T}'=\mathbf{S}_{h'}$.
 In particular $\mathbf{S}_{h'}$ is contained in $\mathbf{P}$.
\end{prf}

\begin{defn}
 Fix isomorphisms $\Q_p(1)\cong \Q_p$ and $(V_x,\psi_x)\cong (\Q_p^{2d},\langle\ ,\ \rangle)$ as in Definition \ref{defn:transfer}.
 We define a homomorphism $\iota_{h,x}\colon \mathbf{T}_h\longrightarrow \GSp_{2d}$ of algebraic groups over $\Q_p$ as follows.
 For a $\Q_p$-algebra $R$, each element $t\in \mathbf{T}_h(R)$ gives an automorphism of
 the filtered isocrystal $(R\otimes_{\Q_p}N,\Phi,R\otimes_{\Q_p}\Fil_x)$
 by the previous lemma. The induced automorphism on $(R\otimes_{\Q_p}N\otimes_{L_0}B_{\mathrm{crys}})^{\Phi,\Fil}=R\otimes_{\Q_p}V_x$
 defines an element $\iota_{h,x}(t)\in\GSp(R\otimes_{\Q_p}V_x,\psi_x)\cong \GSp_{2d}(R)$.

 By definition, we have $\iota_{h,x}(h)=g_{h,x}$.
 Clearly $\iota_{h,x}$ is independent of the choice of $\Q_p(1)\cong \Q_p$, and the $\GSp_{2d}(\Q_p)$-conjugacy class of $\iota_{h,x}$ is
 independent of the choice of $(V_x,\psi_x)\cong (\Q_p^{2d},\langle\ ,\ \rangle)$. 
\end{defn}

\begin{prop}
 Let $\mathbf{T}_{g_{h,x}}=\mathbf{Z}(g_{h,x})$ be the centralizer of $g_{h,x}$ in $\GSp_{2d}$.
 The homomorphism $\iota_{h,x}$ induces an isomorphism $\mathbf{T}_h\yrightarrow{\cong} \mathbf{T}_{g_{h,x}}$.
\end{prop}

\begin{prf}
 The proof of Proposition \ref{prop:stable-conj} tells us that the base change of $\iota_{h,x}$ to $B_{\mathrm{dR}}$
 can be described as the composite of
 \[
 \mathbf{T}_h\otimes_{\Q_p}B_{\mathrm{dR}}\hooklongrightarrow \mathbf{J}\otimes_{\Q_p}B_{\mathrm{dR}}
 \yrightarrow[\cong]{\xi}\GSp_{2d}\otimes_{\Q_p}B_{\mathrm{dR}}\yrightarrow[\cong]{\Ad(\alpha^{-1})}\GSp_{2d}\otimes_{\Q_p}B_{\mathrm{dR}},
 \]
 where $\alpha\in\GSp_{2d}(B_{\mathrm{dR}})$ is the element defined in the proof of Proposition \ref{prop:stable-conj}.
 Now it is clear that the homomorphism $\mathbf{T}_h\otimes_{\Q_p}B_{\mathrm{dR}}\longrightarrow \mathbf{T}_{g_{h,x}}\otimes_{\Q_p}B_{\mathrm{dR}}$
 induced by $\iota_{h,x}$ is an isomorphism, and thus 
 $\iota_{h,x}\colon \mathbf{T}_h\longrightarrow \mathbf{T}_{g_{h,x}}$ is also an isomorphism.
\end{prf}

\begin{prop}\label{prop:st-conj-mod-conj}
 \begin{enumerate}
  \item For $h'\in \{h\}_{\mathrm{st}}$, take $\gamma\in \mathbf{J}(\overline{\Q}_p)$ such that $h'=\gamma^{-1}h\gamma$.
	Then, for every $\sigma\in\Gal(\overline{\Q}_p/\Q_p)$, $c_{h,h'}(\sigma)=\gamma\sigma(\gamma)^{-1}$ lies in $\mathbf{T}_h(\overline{\Q}_p)$,
	and $c_{h,h'}$ gives an element of $H^1(\Q_p,\mathbf{T}_h)$.
	The map $h'\longmapsto c_{h,h'}$ induces a bijection $\{h\}_{\mathrm{st}}/{\sim}\yrightarrow{\cong}H^1(\Q_p,\mathbf{T}_h)$.
	Similarly, we have a natural bijection $\{g_{h,x}\}_{\mathrm{st}}/{\sim}\yrightarrow{\cong}H^1(\Q_p,\mathbf{T}_{g_{h,x}})$.
  \item Assume that $h$ is elliptic. Then we have a commutative diagram
	\[
	 \xymatrix{%
	\{h\}_{\mathrm{st}}/{\sim}\ar[r]\ar[d]^-{(*)}& H^1(\Q_p,\mathbf{T}_h)\ar[d]^-{\iota_{h,x}}\\
	\{g_{h,x}\}_{\mathrm{st}}/{\sim}\ar[r]& H^1(\Q_p,\mathbf{T}_{g_{h,x}})\lefteqn{,}
	}
	\]
	where the map $(*)$ is given by $h'=\gamma^{-1}h\gamma\longmapsto g_{h',\xi(\gamma)^{-1}x}$.
	(Recall that $\xi$ carries $\gamma\in\mathbf{J}(\overline{\Q}_p)$ to an element of 
	$\GSp(N^\circ\otimes_{\Q_{p^2}}\overline{\Q}_p,\psi)\subset \GSp(N\otimes_{L_0}\overline{L_0},\psi)$.
	Corollary \ref{cor:J-inner-form} ensures that $\xi(\gamma)^{-1}x$ is fixed by $h'$,
	and Proposition \ref{prop:point-counting-period} tells us that $\xi(\gamma)^{-1}x$ lies in $\Omega$.)
 \end{enumerate}
\end{prop}

\begin{prf}
 i) It is well-known that the map $h'\longmapsto c_{h,h'}$ induces a bijection 
 \[
  \{h\}_{\mathrm{st}}/{\sim}\yrightarrow{\cong} \Ker\bigl(H^1(\Q_p,\mathbf{T}_h)\longrightarrow H^1(\Q_p,\mathbf{J})\bigr).
 \]
 Therefore it suffices to show $H^1(\Q_p,\mathbf{J})=1$. Since the derived group $\mathbf{J}_{\mathrm{der}}$ is simply connected,
 we have $H^1(\Q_p,\mathbf{J}_{\mathrm{der}})=1$ by Kneser's theorem. Now the exact sequence
 \[
  1\longrightarrow \mathbf{J}_{\mathrm{der}}\longrightarrow \mathbf{J}\longrightarrow \mathbb{G}_m\longrightarrow 1
 \]
 gives the desired vanishing result. Similarly we can prove $H^1(\Q_p,\GSp_{2d})=1$.

 ii) Take a finite extension $F$ of $\Q_{p^2}$ and $x'\in\boldsymbol{\mathcal{F}}^\circ(F)$ as in the proof of Lemma \ref{lem:torus-stabilizer}.
 Let $E$ be a finite Galois extension of $\Q_{p^2}$ such that $\gamma\in \mathbf{J}(E)$.
 Extending $F$ if necessary, we may assume that $E\subset F$ (although we can take $E=F$, it is better to
 distinguish them in order to avoid confusion).  
 Put $y'=\xi(\gamma)^{-1}(x')\in\boldsymbol{\mathcal{F}}^\circ(F)$.

 First we prove that the subspace $E\otimes_{\Q_p}\Fil_{x'}\subset E\otimes_{\Q_p}(F\otimes_{\Q_{p^2}}N^\circ)$ 
 is mapped to $E\otimes_{\Q_p}\Fil_{y'}$ by $\gamma^{-1}\in\mathbf{J}(E)$.
 Put $W=\gamma^{-1}(E\otimes_{\Q_p}\Fil_{x'})\subset E\otimes_{\Q_p}(F\otimes_{\Q_{p^2}}N^\circ)$.
 As $h'=\gamma^{-1}h\gamma\in\mathbf{J}(\Q_p)$, we have $\gamma\sigma(\gamma)^{-1}\in\mathbf{T}_h(E)$
 for every $\sigma\in\Gal(E/\Q_{p^2})$. By Lemma \ref{lem:torus-stabilizer}, we have
 $W=\sigma(\gamma)^{-1}(E\otimes_{\Q_p}\Fil_{x'})$.
 The commutative diagram  
 \[
  \xymatrix{%
 E\otimes_{\Q_p}(F\otimes_{\Q_{p^2}}N^\circ)\ar[r]^-{\gamma^{-1}}\ar[d]^-{\sigma\otimes\id}&
 E\otimes_{\Q_p}(F\otimes_{\Q_{p^2}}N^\circ)\ar[d]^-{\sigma\otimes\id}\\
 E\otimes_{\Q_p}(F\otimes_{\Q_{p^2}}N^\circ)\ar[r]^-{\sigma(\gamma)^{-1}}&
 E\otimes_{\Q_p}(F\otimes_{\Q_{p^2}}N^\circ)
 }
 \]
 tells us that $(\sigma\otimes\id)(W)=W$. Therefore, by the Galois descent, there exists an $F$-subspace
 $W'\subset F\otimes_{\Q_{p^2}}N^\circ$ such that $W=E\otimes_{\Q_p}W'$.
 It suffices to show that $W'=\Fil_{y'}$.
 By definition, $W'$ is the image of $\Fil_{x'}$ under the map
 $F\otimes_{\Q_{p^2}}N^\circ\yrightarrow{(*)} F\otimes_{\Q_{p^2}}N^\circ$ obtained as the base change
 of $\gamma^{-1}\colon E\otimes_{\Q_p}(F\otimes_{\Q_{p^2}}N^\circ)\longrightarrow E\otimes_{\Q_p}(F\otimes_{\Q_{p^2}}N^\circ)$
 by $E\otimes_{\Q_p}F\longrightarrow F$. 
 Corollary \ref{cor:J-inner-form} tells us that the base change of $\gamma^{-1}\colon E\otimes_{\Q_p}N^\circ\longrightarrow E\otimes_{\Q_p}N^\circ$ by
 $E\otimes_{\Q_p}\Q_{p^2}\longrightarrow E$ coincides with $\xi(\gamma)^{-1}\colon E\otimes_{\Q_{p^2}}N^\circ\longrightarrow E\otimes_{\Q_{p^2}}N^\circ$. Hence $(*)$ is equal to the base change of $\xi(\gamma)^{-1}$ by
 $E\longrightarrow F$, and we conclude that $W'=\xi(\gamma)^{-1}(\Fil_{x'})=\Fil_{y'}$, as desired.
 
 Now we know that $\gamma\colon E\otimes_{\Q_p}N^\circ\yrightarrow{\cong} E\otimes_{\Q_p}N^\circ$ gives an isomorphism of filtered isocrystals
 $(E\otimes_{\Q_p}N^\circ,\Phi,E\otimes_{\Q_p}\Fil_{y'})\yrightarrow{\cong} (E\otimes_{\Q_p}N^\circ,\Phi,E\otimes_{\Q_p}\Fil_{x'})$.
 It induces an isomorphism of corresponding $\Gal(\overline{F}/F)$-representations
 $\widetilde{\gamma}\colon E\otimes_{\Q_p}V_{y'}\yrightarrow{\cong} E\otimes_{\Q_p}V_{x'}$.
 Since $h'=\sigma(\gamma)^{-1}h\sigma(\gamma)$ and $y'=\xi(\sigma(\gamma))^{-1}(x')$ for each $\sigma\in\Gal(E/\Q_p)$,
 we may also define $\widetilde{\sigma(\gamma)}\colon E\otimes_{\Q_p}V_{y'}\yrightarrow{\cong} E\otimes_{\Q_p}V_{x'}$,
 which clearly coincides with $\sigma(\widetilde{\gamma})=(\sigma\otimes \id)\circ \widetilde{\gamma}\circ (\sigma\otimes \id)^{-1}$ by functoriality.
 By the construction, we have $\widetilde{\gamma}\sigma(\widetilde{\gamma})^{-1}=\widetilde{\gamma}\widetilde{\sigma(\gamma})^{-1}=\iota_{h,x}(\gamma\sigma(\gamma)^{-1})$.
 The relation $h'=\gamma^{-1}h\gamma$ is translated into the commutativity of the following diagram:
 \[
 \xymatrix{%
 E\otimes_{\Q_p}V_{y'}\ar[r]^-{g_{h',y'}}\ar[d]^-{\widetilde{\gamma}}& E\otimes_{\Q_p}V_{y'}\ar[d]^-{\widetilde{\gamma}}\\
 E\otimes_{\Q_p}V_{x'}\ar[r]^-{g_{h,x'}}& E\otimes_{\Q_p}V_{x'}\lefteqn{.}
 }
 \]
 In other words, if we fix $(V_{x'},\psi_{x'})\cong (\Q_p^{2d},\langle\ ,\ \rangle)\cong (V_{y'},\psi_{y'})$ and
 regard $\widetilde{\gamma}$ as an element of $\GSp_{2d}(E)$, then we have $g_{h',y'}=\widetilde{\gamma}^{-1}g_{h,x'}\widetilde{\gamma}$.
 Therefore, under the isomorphism $\{g_{h,x}\}/{\sim}\yrightarrow{\cong} H^1(\Q_p,\mathbf{T}_{g_{h,x}})$,
 the conjugacy class of $g_{h',y'}$ corresponds to the cohomology class of the cocycle
 $\sigma\longmapsto \widetilde{\gamma}\sigma(\widetilde{\gamma})^{-1}$.
 This concludes the proof, since we have $\widetilde{\gamma}\sigma(\widetilde{\gamma})^{-1}=\iota_{h,x}(\gamma\sigma(\gamma)^{-1})=\iota_{h,x}(c_{h,h'}(\sigma))$
 as mentioned above. 
\end{prf}

\begin{cor}\label{cor:average-conj-class}
 Assume that $h$ is elliptic. Let $g_h\in G$ be an arbitrary element with $g_h\leftrightarrow h$.
 Then, each $[g]\in \{g_h\}_{\mathrm{st}}/{\sim}$ satisfies the following:
 \[
  \#\bigl\{([h'],x)\bigm| [h']\in\{h\}_{\mathrm{st}}/{\sim}, x\in \Fix(h';\Omega), [g_{h',x}]=[g]\bigr\}=d!2^{d-1}.
 \]
 Here $[-]$ denotes the conjugacy class.
\end{cor}

\begin{prf}
 For $h'\in\{h\}_{\mathrm{st}}$, take $\gamma\in\mathbf{J}(\overline{\Q}_p)$ such that $h'=\gamma^{-1}h\gamma$.
 Then $\xi(\gamma)^{-1}$ induces a bijection $\Fix(h;\Omega)\longrightarrow \Fix(h';\Omega)$.
 Thus by Proposition \ref{prop:st-conj-mod-conj} we have
 \begin{align*}
  \#\{x\in\Fix(h',\Omega)\mid [g_{h',x}]=[g]\}&=\#\{x\in\Fix(h,\Omega)\mid [g_{h',\xi(\gamma)^{-1}x}]=[g]\}\\
  &=\#\{x\in\Fix(h,\Omega)\mid \iota_{h,x}(c_{h'})=c_g\},
 \end{align*}
 where $c_{h'}\in H^1(\Q_p,\mathbf{T}_h)$ (resp.\ $c_g\in H^1(\Q_p,\mathbf{T}_{g_{h,x}})$)
 is the element corresponding to $[h']\in \{h\}_{\mathrm{st}}/{\sim}$
 (resp.\ $[g]\in \{g_h\}_{\mathrm{st}}/{\sim}=\{g_{h,x}\}_{\mathrm{st}}/{\sim}$)
 under the bijection in Proposition \ref{prop:st-conj-mod-conj} i).
 Varying $[h']$ (or equivalently $c_{h'}$), we obtain
 \begin{align*}
  &\#\bigl\{([h'],x)\bigm| [h']\in\{h\}_{\mathrm{st}}/{\sim}, x\in \Fix(h';\Omega), [g_{h',x}]=[g]\bigr\}\\
  &\qquad\qquad=\#\{(c,x)\in H^1(\Q_p,\mathbf{T}_h)\times \Fix(h;\Omega)\mid \iota_{h,x}(c)=c_g\},
 \end{align*}
 whose right hand side is clearly equal to $\#\Fix(h;\Omega)=d!2^{d-1}$
 (\cf Remark \ref{rem:fixed-point-number-period}).
\end{prf}

The following is an analogue of \cite[Theorem 1]{MR1302321} and \cite[Theorem 2.6.8]{MR2383890}:

\begin{prop}\label{prop:fiber-fixed-point-number}
 Let $x$ be as in Definition \ref{defn:transfer}. Assume that $g\in G$ is regular elliptic.
 Then, the number of fixed points in $\wp_{K,p}^{-1}(x)$ under the action of $(g,h)$ is given by
 \[
  \#\{\gamma Kp^\Z\in G/Kp^\Z\mid \gamma Kp^\Z=g_{h,x}^{-1}\gamma gKp^\Z\}.
 \]
 If $h$ is elliptic, this number is equal to the orbital integral
 \[
 O_{g_{h,x}}\Bigl(\frac{\mathbf{1}_{gKp^\Z}}{\vol(K)}\Bigr),
 \]
 where $\mathbf{1}_{gKp^\Z}$ denotes the characteristic function of $gKp^\Z\subset G$
 (for our normalization of the Haar measure, see the last paragraph of Section \ref{subsec:def-RZ}).
\end{prop}

\begin{prf}
 As in Theorem \ref{thm:RZ-period} ii), there exist a finite extension $L'$ of $L$ and
 $\widetilde{x}=(X,\rho,\overline{\eta})\in M_K(L')$ such that $\wp_K(\widetilde{x})=x$. 
 Here, $(X,\rho)$ is an $\mathcal{O}_{L'}$-valued point of $\M$ and
 $\overline{\eta}$ is a $K$-level structure on $X$, namely,
 a $\Gal(\overline{L}/L')$-invariant $K$-orbit of isomorphisms $\Q_p^{2d}\yrightarrow{\cong} V_p(X)$
 preserving symplectic pairings up to $\Q_p^\times$-multiplication.
 Fix $\eta\in \overline{\eta}$.
 Then, by \cite[Proposition 5.37]{MR1393439}, the fiber $\wp_K^{-1}(x)$ of
 $\wp_K\colon M_K\longrightarrow \Omega$ at $x$
 can be identified with $G/K$. The identification is given as follows.
 Let $\widetilde{x}'=(X',\rho',\overline{\eta}')$ be another point in the fiber.
 Then, there exists a unique quasi-isogeny $f\colon X'\longrightarrow X$ satisfying
 $\rho=(f\bmod p)\circ \rho'$. This $f$ automatically preserves polarizations on $X$ and $X'$
 up to $\Q_p^\times$-multiplication.
 Choose a representative $\eta'$ in the level structure $\overline{\eta}'$,
 and define $\gamma\in G$ by the following commutative diagram:
 \[
  \xymatrix{%
 \Q_p^{2d}\ar[r]^-{\eta'}_-{\cong}\ar[d]^-{\gamma}& V_p(X')\ar[d]^-{V_p(f)}\\
 \Q_p^{2d}\ar[r]^-{\eta}_-{\cong}& V_p(X)\lefteqn{.}
 }
 \]
 The class of $\gamma$ in $G/K$ does not depend on the choice of $\eta'$.

 Since $(g,h)$ acts on $\wp^{-1}_K(x)$, it also acts on $G/K$ by the identification above.
 We will show that this action is given by $\gamma K\longmapsto g_{h,x}^{-1}\gamma gK$, where the element
 $g_{h,x}\in G$ is given by the isomorphism $V_x\cong V_p(X)\yrightarrow[\cong]{\eta^{-1}}\Q_p^{2d}$
 (\cf Definition \ref{defn:transfer}). Take $(X',\rho',\overline{\eta}')\in\wp_K^{-1}(x)$ corresponding to
 $\gamma K\in G/K$. Then $(g,h)(X',\rho',\overline{\eta}')=(X',\rho'\circ h,\overline{\eta'\circ g})$.
 Since $\rho_*\circ h\circ (\rho_*)^{-1}\colon D(X)_{\Q}\longrightarrow D(X)_{\Q}$ preserves 
 the Hodge filtration, there exists a quasi-isogeny $\phi\colon X\longrightarrow X$ such that 
 $\rho\circ h=(\phi\bmod p)\circ \rho$.
 For this $\phi$, we have $\rho=((\phi^{-1}\circ f)\bmod p)\circ (\rho'\circ h)$, where $f\colon X'\longrightarrow X$
 is a quasi-isogeny as above.
 Now we have the following commutative diagram:
 \[
  \xymatrix{%
 \Q_p^{2d}\ar[r]^-{\eta'\circ g}_-{\cong}\ar[d]^-{g}& V_p(X')\ar@{=}[d]\\
 \Q_p^{2d}\ar[r]^-{\eta'}_-{\cong}\ar[d]^-{\gamma}& V_p(X')\ar[d]^-{V_p(f)}\\
 \Q_p^{2d}\ar[r]^-{\eta}_-{\cong}\ar[d]^-{g_{h,x}^{-1}}& V_p(X)\ar[d]^-{V_p(\phi)^{-1}}\\
 \Q_p^{2d}\ar[r]^-{\eta}_-{\cong}& V_p(X)\lefteqn{.}
  }
 \]
 This diagram tells us that $(g,h)(X',\rho',\overline{\eta}')$ corresponds to $g_{h,x}^{-1}\gamma gK\in G/K$,
 as desired.

 Now we consider the fiber $\wp_{K,p}^{-1}(x)=\wp_{K}^{-1}(x)/p^\Z$. 
 By \cite[Lemma 5.36]{MR1393439}, it is in bijection with $G/Kp^\Z$. 
 The action of $(g,h)$ is given by $\gamma Kp^\Z\longmapsto g_{h,x}^{-1}\gamma gKp^\Z$.
 Therefore, the number of fixed points under this action is equal to
 \[
  \#\{\gamma Kp^\Z\in G/Kp^\Z\mid \gamma Kp^\Z=g_{h,x}^{-1}\gamma gKp^\Z\}.
 \]
 If $h$ is elliptic, then $g_{h,x}$ is regular elliptic by Proposition \ref{prop:stable-conj}, and thus
 the number above is equal to
 \begin{align*}
  \frac{1}{\vol(Kp^\Z/p^\Z)}\int_{G/p^\Z}\mathbf{1}_{gKp^\Z}(\gamma^{-1}g_{h,x}\gamma)d\gamma
  =O_{g_{h,x}}\Bigl(\frac{\mathbf{1}_{gKp^\Z}}{\vol(K)}\Bigr).
 \end{align*}
\end{prf}

For later use, we will state our result for the action of the inverse element $(g^{-1},h^{-1})$.

\begin{thm}\label{thm:fixed-point-number}
 Let $(g,h)\in G\times J$ be an element and $K$ a compact open subgroup of $K_0$ normalized by $g$.
 \begin{enumerate}
  \item If $h$ is regular elliptic, then we have
	\[
	\#\Fix\bigl((g^{-1},h^{-1});M_K/p^\Z\bigr)=\sum_{x\in\Fix(h;\Omega)}O_{g_{h,x}}\Bigl(\frac{\mathbf{1}_{gKp^\Z}}{\vol(K)}\Bigr),
	\]
	where the Haar measure on $Z(g_{h,x})$ is normalized as in
	the last paragraph of Section \ref{subsec:def-RZ}.
	The left hand side denotes the number of fixed points in the sense of \cite[Definition 2.6]{adicLTF}.
  \item Assume that $h$ is regular non-elliptic, and $gK$ consists of regular elliptic elements.
	Then there is no point on $M_K/p^\Z$ fixed by $(g^{-1},h^{-1})$.
 \end{enumerate}
\end{thm}

\begin{prf}
 i) Since the multiplicity at each fixed point in $\Omega$ is one and $\wp_{K,p}$ is \'etale,
 the multiplicity at each fixed point in $M_K/p^\Z$ is equal to one (\cf \cite[Proposition 2.10, Proposition 2.11]{adicLTF}). Therefore, we have only to count points in $M_K/p^\Z$ fixed by $(g^{-1},h^{-1})$ 
 (or equivalently, fixed by $(g,h)$) as sets. The desired equality immediately follows from
 Proposition \ref{prop:fiber-fixed-point-number}.

 ii) If there exists a fixed point, then we can find $\gamma\in G$ and an integer $n$ such that 
 $p^{-n}\gamma^{-1}g_{h,x}\gamma\in gK$. Since $gK$ consists of regular elliptic elements,
 $g_{h,x}$ is also regular elliptic. By Proposition \ref{prop:stable-conj},
 this contradicts to the assumption that $h$ is non-elliptic.
\end{prf}

\begin{cor}\label{cor:Fix-SO}
 In the setting of Theorem \ref{thm:fixed-point-number} i), we have
 \[
 \sum_{[h']\in \{h\}_{\mathrm{st}}/{\sim}}\#\Fix\bigl((g^{-1},h^{-1});M_K/p^\Z\bigr)
 =d!2^{d-1} \mathit{SO}_{g_h}\Bigl(\frac{\mathbf{1}_{gKp^\Z}}{\vol(K)}\Bigr),
 \]
 where $g_h\in G$ is an arbitrary element with $g_h\leftrightarrow h$.
\end{cor}

\begin{prf}
 By Theorem \ref{thm:fixed-point-number} i) and Corollary \ref{cor:average-conj-class}, we have
 \begin{align*}
  &\sum_{[h']\in \{h\}_{\mathrm{st}}/{\sim}}\#\Fix\bigl((g^{-1},h^{-1});M_K/p^\Z\bigr)
  =\sum_{[h']\in \{h\}_{\mathrm{st}}/{\sim}}\sum_{x\in\Fix(h';\Omega)}O_{g_{h,x}}\Bigl(\frac{\mathbf{1}_{gKp^\Z}}{\vol(K)}\Bigr)\\
  &\qquad =d!2^{d-1}\sum_{[g']\in \{g_h\}_{\mathrm{st}}/{\sim}}O_{g'}\Bigl(\frac{\mathbf{1}_{gKp^\Z}}{\vol(K)}\Bigr)=d!2^{d-1} \mathit{SO}_{g_h}\Bigl(\frac{\mathbf{1}_{gKp^\Z}}{\vol(K)}\Bigr),
 \end{align*}
 as desired.
\end{prf}

\section{Some formal models}
In order to apply the Lefschetz trace formula \cite[Theorem 4.5]{adicLTF},
we need to construct a formal model of $M_K$ for $K$ in a certain family of compact open subgroups of $K_0$.

\subsection{Construction of formal models}
First let us recall the definition of a chain of lattices of $\Q_p^{2d}$ considered in \cite[\S 3]{MR1393439}:

\begin{defn}
 The set $\mathscr{L}$ of $\Z_p$-lattices of $\Q_p^{2d}$ is said to be a chain of lattices
 if the following hold (\cf \cite[Definition 3.1]{MR1393439}):
\begin{itemize}
 \item for $L,L'\in \mathscr{L}$, we have either $L\subset L'$ or
       $L\supset L'$,
 \item and for $L\in\mathscr{L}$ and $a\in\Q_p^\times$, we have $aL\in\mathscr{L}$.
\end{itemize}
 Moreover, $\mathscr{L}$ is said to be self-dual if $L\in\mathscr{L}$ implies $L^\vee\in\mathscr{L}$,
 where $L^\vee$ denotes the dual lattice with respect to the fixed symplectic form $\langle\ ,\ \rangle$
 on $\Q_p^{2d}$.
\end{defn}

For a self-dual chain of lattices $\mathscr{L}$, write $K_{\!\mathscr{L}}$ for the stabilizer of $\mathscr{L}$ in $G$.
It is a parahoric subgroup of $G$.
Moreover, for an integer $m\ge 0$, we put
\[
 K_{\!\mathscr{L},m}=\{g\in K_{\!\mathscr{L}}\mid \text{for every $L\in \mathscr{L}$, $g$ acts trivially on $L/p^mL$}\},
\]
which is an open normal subgroup of $K_{\!\mathscr{L}}$.

Let $N_{\!\mathscr{L}}$ be the set of $g\in G$ such that $g\mathscr{L}=\mathscr{L}$.
It coincides with the normalizer of $K_{\!\mathscr{L}}$ in $G$.
Furthermore, $N_{\mathscr{L}}$ also normalizes
$K_{\!\mathscr{L},m}$ for every $m\ge 0$.
It is known that $N_{\!\mathscr{L}}$ is a compact-mod-center subgroup of $G$.

In this subsection, we will construct a formal model $\M_{\mathscr{L},m}$ of $M_{K_{\!\mathscr{L},m}}$
on which $N_{\mathscr{L}}$ acts as isomorphisms.

In the following, we fix a self-dual chain of lattices $\mathscr{L}$.
First we recall the Rapoport-Zink space with parahoric level $\M_{\mathscr{L}}$
introduced in \cite[Definition 3.21]{MR1393439}:

\begin{defn}
 Let $\M_{\mathscr{L}}$ be the contravariant functor $\Nilp\longrightarrow \mathbf{Set}$ that associates $S$
 with the set of isomorphism classes of $\{(X_L,\rho_L)\}_{L\in\mathscr{L}}$ where
 \begin{itemize}
 \item $X_L$ is a $d$-dimensional $p$-divisible group $X$ over $S$,
 \item and $\rho_L\colon \X\otimes_{\overline{\F}_p}\overline{S}\longrightarrow X_L\otimes_S\overline{S}$ is
       a quasi-isogeny (as in the definition of $\M$, we put $\overline{S}=S\otimes_{\Z_{p^\infty}}\overline{\F}_p$),
\end{itemize}
 satisfying the following conditions:
\begin{itemize}
 \item For $L,L'\in \mathscr{L}$ with $L\subset L'$, the quasi-isogeny $\rho_{L'}\circ\rho_L^{-1}\colon X_L\otimes_S\overline{S}\longrightarrow X_{L'}\otimes_S\overline{S}$ lifts to an isogeny $\widetilde{\rho}_{L',L}\colon X_L\longrightarrow X_{L'}$.
 \item For $L,L'$ as above, $\deg\widetilde{\rho}_{L',L}=\log_p\#(L'/L)$.
 \item For $L\in\mathscr{L}$, the quasi-isogeny $\rho_{pL}\circ [p]\circ\rho_L^{-1}\colon X_L\otimes_S\overline{S}\longrightarrow X_{pL}\otimes_S\overline{S}$ lifts to an isomorphism $\theta_p\colon X_L\yrightarrow{\cong} X_{pL}$
       (here $[p]$ denotes the multiplication by $p$).
 \item There exists a constant $a\in\Q_p^\times$ such that for every $L\in\mathscr{L}$, we can find
       an isomorphism $\lambda_L\colon X_L\longrightarrow (X_{L^\vee})^\vee$ which makes the following
       diagram commutative:
	\[
	 \xymatrix{%
	 \mathbb{X}\otimes_{\overline{\F}_p}\overline{S}\ar[rr]^-{\rho_L}\ar[d]^-{a\lambda_0\otimes\id}&& X_L\otimes_S\overline{S}\ar[d]^-{\lambda_L\otimes \id}\\
	 \mathbb{X}^\vee\otimes_{\overline{\F}_p}\overline{S}&& (X_{L^\vee})^\vee\otimes_S\overline{S}\ar[ll]_-{(\rho_{L^\vee})^\vee}
	}
	\]
       (such $\lambda_L$ automatically satisfies $\lambda_L^\vee=-\lambda_{L^\vee}$).
\end{itemize}
 The functor $\M_{\mathscr{L}}$ is represented by a special formal scheme over $\Spf\Z_{p^\infty}$.

 The group $N_{\mathscr{L}}\times J$ naturally acts on $\M_{\mathscr{L}}$ on the right;
 for $(g,h)\in N_{\mathscr{L}}\times J$, 
 the action $(g,h)\colon \M_{\mathscr{L}}(S)\longrightarrow \M_{\mathscr{L}}(S)$ 
 is given by $\{(X_L,\rho_L)\}\longmapsto \{(X_{gL},\rho_{gL}\circ h)\}$.
 It is easy to see that $(g,1)$ with $g\in K_{\!\mathscr{L}}$ and $(p,p)$
 act trivially on $\M_{\mathscr{L}}$.
\end{defn}

Next we consider level structures on $\M_\mathscr{L}$.

\begin{defn}\label{defn:formal-model}
 For an integer $m\ge 0$, let $\M'_{\mathscr{L},m}$ be the contravariant functor
 $\Nilp\longrightarrow \mathbf{Set}$ that associates $S$
 with the set of isomorphism classes of $\{(X_L,\rho_L,\eta_{L,m})\}_{L\in\mathscr{L}}$ where
 \begin{itemize}
  \item $\{(X_L,\rho_L)\}_{L\in\mathscr{L}}\in\M_{\mathscr{L}}(S)$, and
  \item $\eta_{L,m}\colon L/p^mL\longrightarrow X_L[p^m]$ is a homomorphism
 \end{itemize}
 satisfying the following conditions:
 \begin{itemize}
  \item $\eta_{L,m}$ is a Drinfeld $m$-level structure. Namely, the image of $\eta_{L,m}$ gives
	a full set of sections of $X_L[p^m]$ (\cf \cite[\S 1.8]{MR772569}).
  \item For $L,L'\in\mathscr{L}$ with $L\subset L'$, the following diagrams are commutative:
	\[
	 \xymatrix{
	L/p^mL\ar[d]^{\eta_{L,m}}\ar[r]& L'/p^mL'\ar[d]^{\eta_{L',m}}\\
	X_L[p^m]\ar[r]^-{\widetilde{\rho}_{L',L}}& X_{L'}[p^m]\lefteqn{,}
	}
	\qquad\qquad
	\xymatrix{
	L/p^mL\ar[d]^{\eta_{L,m}}\ar[r]^-{\times p}\ar[r]_-{\cong}& pL/p^{m+1}L\ar[d]^{\eta_{pL,m}}\\
	X_L[p^m]\ar[r]_-{\cong}^-{\theta_p}& X_{pL}[p^m]\lefteqn{.}
	}
	\]
  \item There exists a homomorphism $\Z/p^m\Z\longrightarrow \mu_{p^m}$ such that for every $L\in \mathscr{L}$
	the diagram below is commutative up to constant which is independent of $L$:
	\[
	 \xymatrix{%
	L/p^mL\times L^\vee/p^mL^\vee\ar[d]^-{\eta_{L,m}\times \eta_{L^\vee,m}}\ar[rr]^-{\langle\ ,\ \rangle}&&
	\Z/p^m\Z\ar[d]\\
	X_L[p^m]\times X_{L^\vee}[p^m]\ar[r]^-{\id\times \lambda_{L^\vee}}&
	X_L[p^m]\times X^\vee_L[p^m]\ar[r]& \mu_{p^m}\lefteqn{.}
	}
	\]
 \end{itemize}
 It is easy to see that $\M'_{\mathscr{L},m}$ is represented by a formal scheme which is finite
 over $\M_{\mathscr{L}}$ (\cf \cite[Proposition 1.9.1]{MR772569}).
 The group $N_{\mathscr{L}}\times J$ naturally acts on $\M'_{\mathscr{L},m}$ on the right;
 by $(g,h)\in N_{\mathscr{L}}\times J$, $\{(X_L,\rho_L,\eta_{L,m})\}_{L\in\mathscr{L}}$ is mapped to
 $\{(X_{gL},\rho_{gL}\circ h,\eta_{gL,m}\circ g)\}_{L\in\mathscr{L}}$.
 It is easy to see that $(g,1)$ with $g\in K_{\!\mathscr{L},m}$ and $(p,p)$
 act trivially on $\M'_{\mathscr{L},m}$.
 By \cite[Lemma 7.2]{MR2074715}, $\{\M'_{\mathscr{L},m}\}_{m\ge 0}$ forms a projective system of formal schemes
 equipped with actions of $N_{\mathscr{L}}\times J$.

 Let $\M_{\mathscr{L},m}$ be the intersection of the scheme-theoretic images of
 $\M'_{\mathscr{L},m'}\longrightarrow \M'_{\mathscr{L},m}$ for all $m'\ge m$.
 Obviously $\{\M_{\mathscr{L},m}\}_{m\ge 0}$ again forms a projective system of formal schemes.

 Finally, let $\M^\flat_{\mathscr{L},m}$ be the closed formal subscheme of $\M_{\mathscr{L},m}$
 defined by the quasi-coherent ideal of $\mathcal{O}_{\!\!\M_{\mathscr{L},m}}$ consisting of elements
 killed by $p^l$ for some integer $l\ge 0$. It is flat over $\Spf\Z_{p^\infty}$.
 We obtain a projective system of formal schemes $\{\M^\flat_{\mathscr{L},m}\}_{m\ge 0}$
 equipped with actions of $N_{\mathscr{L}}\times J$. 
 Obviously we can identify $\M^{\flat,\rig}_{\mathscr{L},m}$ with $\M^{\rig}_{\mathscr{L},m}$.
\end{defn}

\begin{rem}\label{rem:m=0}
 By \cite[Proposition A.21]{MR1393439}, we can show that the natural map
 $\M^\rig_{\mathscr{L},m}\longrightarrow \M^\rig_{\mathscr{L}}$ is surjective for every $m\ge 0$.
 In particular, we have $\M^\rig_{\mathscr{L},0}=\M^\rig_{\mathscr{L}}$.
\end{rem}

The formal scheme $\M^\flat_{\mathscr{L},m}$ gives a formal model of the rigid space $M_{K_{\!\mathscr{L},m}}$:

\begin{prop}\label{prop:formal-model-gen-fiber}
 Assume that $K_{\!\mathscr{L},m}\subset K_0$. Then,
 we have a natural isomorphism $\M^{\flat,\rig}_{\mathscr{L},m}\cong M_{K_{\!\mathscr{L},m}}$ of
 rigid spaces over $\Q_{p^\infty}$ which is compatible with actions of $N_{\mathscr{L}}\times J$
 and change of $m$.
\end{prop}

\begin{prf}
 First we construct a morphism $M_{K_{\!\mathscr{L},m}}\longrightarrow \M'^\rig_m$.
 Let $\mathcal{S}$ be a formal scheme of finite type over $\Spf\Z_{p^\infty}$ such that
 $\mathcal{S}^\rig$ is connected, and $\mathcal{S}^\rig\longrightarrow M_{K_{\!\mathscr{L},m}}$ a morphism
 over $\Q_{p^\infty}$. It suffices to construct a morphism $\mathcal{S}^\rig\longrightarrow \M'^\rig_m$.

 By changing $\mathcal{S}$ by its admissible blow-up, we may assume that
 the composite $\mathcal{S}^\rig\longrightarrow M_{K_{\!\mathscr{L},m}}\longrightarrow M$
 extends to $\mathcal{S}\longrightarrow \M$.
 Therefore, we have the following data:
 \begin{itemize}
  \item a $p$-divisible group $X$ on $\mathcal{S}$,
  \item a quasi-isogeny $\rho\colon \X\otimes_{\overline{\F}_p}\overline{\mathcal{S}}\longrightarrow X\times_{\mathcal{S}}\overline{\mathcal{S}}$ satisfying the same condition as in the definition of $\M$,
  \item and a $\pi_1(\mathcal{S}^\rig,\overline{x})$-stable $K_{\!\mathscr{L},m}$-orbit $\overline{\eta}$
	of an isomorphism $\eta\colon \Z_p^{2d}\yrightarrow{\cong}T_pX^{\rig}_{\overline{x}}$
	which preserves polarizations up to multiplication by $\Z_p^\times$.
 \end{itemize}
 Here we put $\overline{\mathcal{S}}=\mathcal{S}\otimes_{\Z_{p^\infty}}\overline{\F}_p$ and
 fix a geometric point $\overline{x}$ of $\mathcal{S}^{\rig}$.

 Fix $\eta\in \overline{\eta}$. It corresponds to a homomorphism 
 $\eta\colon \Q_p^{2d}/\Z_p^{2d}\longrightarrow X^\rig_{\overline{x}}$.
 Choose $L_0\in\mathscr{L}$ with $\Z_p^{2d}\subset L_0$ and consider
 $L\in \mathscr{L}$ with $L_0\subset L\subsetneq p^{-1}L_0$. The image of $L/\Z_p^{2d}\subset \Q_p^{2d}/\Z_p^{2d}$
 under $\eta$ corresponds to a finite \'etale subgroup scheme of $X^\rig$. 
 Since there are only finitely many such $L$,
 by the flattening theorem (\cf \cite{MR1225983}), after replacing $\mathcal{S}$ by an admissible blow-up
 we may assume that for each $L$ there exists a finite flat subgroup scheme $Y_L$ of $X$
 whose rigid generic fiber $Y^\rig_L$ corresponds to the image $\eta(L/\Z_p^{2d})$.
 Put $X_L=X/Y_L$, and write $\varphi_L\colon X\longrightarrow X_L$ for the canonical isogeny.
 The homomorphism $\eta\colon \Q_p^{2d}/\Z_p^{2d}\longrightarrow X^\rig_{\overline{x}}$ induces
 $\eta_L\colon \Q_p^{2d}/L\longrightarrow X^{\rig}_{L,\overline{x}}$, which corresponds to
 $\eta_L\colon L\longrightarrow T_pX^{\rig}_{L,\overline{x}}$. It is easy to see that this homomorphism
 fits into the following commutative diagram:
 \[
  \xymatrix{%
 0\ar[r] &\Z_p^{2d}\ar[r]\ar[d]^{\eta}_{\cong}& L\ar[r]\ar[d]^{\eta_L}& L/\Z_p^{2d}\ar[r]\ar[d]_{\cong}& 0\\
 0\ar[r]&T_pX^\rig_{\overline{x}}\ar[r]& T_pX^\rig_{L,\overline{x}}\ar[r]& Y_L^\rig\ar[r]& 0\lefteqn{.}
 }
 \]
 In particular, $\eta_L$ is an isomorphism.

 For $L,L'\in\mathscr{L}$ with $L_0\subset L\subset L'\subsetneq p^{-1}L_0$, we have a natural closed 
 immersion $Y_L^\rig\longrightarrow Y_{L'}^\rig$. Since $Y_L$ and $Y_{L'}$ are finite flat closed subgroup
 schemes of $X$, this extends to a closed immersion $Y_L\longrightarrow Y_{L'}$.
 Therefore we have a natural isogeny $\widetilde{\rho}_{L',L}\colon X_L\longrightarrow X_{L'}$. 
 Clearly the following diagrams are commutative:
 \[
  \xymatrix{%
 \Q_p^{2d}/\Z_p^{2d}\ar[r]\ar[d]^{\eta}&\Q_p^{2d}/L\ar[r]\ar[d]^{\eta_L}& \Q_p^{2d}/L'\ar[d]^{\eta_{L'}}\\
 X^\rig_{\overline{x}}\ar[r]^{\varphi_L}&X^\rig_{L,\overline{x}}\ar[r]^-{\widetilde{\rho}_{L',L}}& X^\rig_{L',\overline{x}}\lefteqn{,}
 }
 \qquad
  \xymatrix{%
 \Z_p^{2d}\ar[r]\ar[d]^{\eta}_{\cong}&L\ar[r]\ar[d]^{\eta_L}_{\cong}& L'\ar[d]^{\eta_{L'}}_{\cong}\\
 T_pX^\rig_{\overline{x}}\ar[r]^{\varphi_L}&T_pX^\rig_{L,\overline{x}}\ar[r]^-{\widetilde{\rho}_{L',L}}& T_pX^\rig_{L',\overline{x}}\lefteqn{.}
 }
 \]
 For $L\in\mathscr{L}$ with $L_0\subset L\subsetneq p^{-1}L_0$, 
 $X_{L_0}\yrightarrow{[p]}X_{L_0}$ factors as 
 $X_{L_0}\yrightarrow{\widetilde{\rho}_{L,L_0}}X_L\yrightarrow{\theta_{L_0,L}} X_{L_0}$.
 By the construction, we have the following commutative diagrams:
 \[
 \xymatrix{%
 \Q_p^{2d}/L\ar[r]\ar[d]^{\eta_L}& \Q_p^{2d}/p^{-1}L_0\ar[d]^{\eta_{L_0}\circ p}\\
 X^\rig_{L,\overline{x}}\ar[r]^-{\theta_{L_0,L}}& X^\rig_{L_0,\overline{x}}\lefteqn{,}
 }
 \qquad
 \xymatrix{%
 L\ar[r]\ar[d]^{\eta_L}_{\cong}& p^{-1}L_0\ar[d]^{\eta_{L_0}\circ p}_{\cong}\\
 T_pX^\rig_{L,\overline{x}}\ar[r]^-{\theta_{L_0,L}}& T_pX^\rig_{L_0,\overline{x}}\lefteqn{.}
 }
 \]
 Next consider the level structure. For $L\in\mathscr{L}$ with $L_0\subset L\subsetneq p^{-1}L_0$, 
 let $K_{L,m}$ be the kernel of $\GL(L)\longrightarrow \GL(L/p^mL)$.
 As the $K_{\!\mathscr{L},m}$-orbit of $\eta_L\colon L\longrightarrow T_pX^\rig_{L,\overline{x}}$
 is $\pi_1(\mathcal{S}^\rig,\overline{x})$-invariant, so is the $K_{L,m}$-orbit of $\eta_L$.
 Therefore we obtain an isomorphism $L/p^mL\yrightarrow{\cong}X^\rig_L[p^m]$.
 Again by the flattening theorem, we may assume that this isomorphism extends to a homomorphism
 $\eta_{m,L}\colon L/p^mL\longrightarrow X_L[p^m]$ of group schemes. 
 This is a Drinfeld $m$-level structure (\cf \cite[Lemma 1.8.3, Proposition 1.9.1]{MR772569}).
 For $L,L'\in\mathscr{L}$ with $L_0\subset L\subset L'\subsetneq p^{-1}L_0$,
 the following diagrams are clearly commutative:
 \[
 \xymatrix{%
 L/p^mL\ar[r]\ar[d]^{\eta_{L,m}}& L'/p^mL'\ar[d]^{\eta_{L',m}}\\
 X_L[p^m]\ar[r]^-{\widetilde{\rho}_{L',L}}& X_{L'}[p^m]\lefteqn{,}
 }
 \qquad
 \xymatrix{
 L/p^mL\ar[d]^{\eta_{L,m}}\ar[r]& p^{-1}L_0/p^{m-1}L_0\ar[d]^{\eta_{L_0,m}\circ p}\\
 X_L[p^m]\ar[r]^-{\theta_{L_0,L}}& X_{L_0}[p^m]\lefteqn{.}
 }
 \]
 Now we extend the construction above to all $L\in\mathscr{L}$;
 take an integer $n$ such that $L_0\subset p^nL\subsetneq p^{-1}L_0$,
 and put $X_L=X_{p^nL}$, $\varphi_L=p^{-n}\circ \varphi_{p^nL}\colon X\longrightarrow X_L$, 
 $\eta_L=\eta_{p^nL}\circ p^n\colon L\yrightarrow{\cong}T_pX^{\rig}_{L,\overline{x}}$ and 
 $\eta_{L,m}=\eta_{p^nL,m}\circ p^n\colon L/p^mL\longrightarrow X_L[p^m]$.
 Then, for $L,L'\in\mathscr{L}$ with $L\subset L'$, we have a natural isogeny 
 $\widetilde{\rho}_{L',L}\colon X_L\longrightarrow X_{L'}$ which makes the following diagrams commutative:
\[
 \xymatrix{%
 L\ar[r]\ar[d]^{\eta_L}_{\cong}& L'\ar[d]^{\eta_{L'}}_{\cong}\\
 T_pX^\rig_{L,\overline{x}}\ar[r]^-{\widetilde{\rho}_{L',L}}& T_pX^\rig_{L',\overline{x}}\lefteqn{,}
 }
 \qquad
 \xymatrix{%
 L/p^mL\ar[r]\ar[d]^{\eta_{L,m}}& L'/p^mL'\ar[d]^{\eta_{L',m}}\\
 X_L[p^m]\ar[r]^-{\widetilde{\rho}_{L',L}}& X_{L'}[p^m]\lefteqn{.}
 }
\]
 Indeed, take integers $n$, $n'$ such that $L_0\subset p^nL\subsetneq p^{-1}L_0$ and
 $L_0\subset p^{n'}L'\subsetneq p^{-1}L_0$ (the assumption $L\subset L'$ implies $n\le n'$).
 If $p^nL\subset p^{n'}L'$, we set
 \[
  \widetilde{\rho}_{L',L}\colon X_L=X_{p^nL}\yrightarrow{\widetilde{\rho}_{p^{n'}L',p^{n}L}\circ [p^{n'-n}]}X_{p^{n'}L'}=X_{L'}.
 \]
 If $p^nL\supset p^{n'}L'$ (in this case $n<n'$), we set
 \[
 \widetilde{\rho}_{L',L}\colon X_L=X_{p^nL}\yrightarrow{\theta_{L_0,p^nL}\circ [p^{n'-n-1}]}X_{L_0}\yrightarrow{\widetilde{\rho}_{p^{n'}L',L_0}}X_{p^{n'}L'}=X_{L'}.
 \]
 Put $\rho_L=\varphi_L\circ \rho\colon \X\otimes_{\overline{\F}_p}\overline{\mathcal{S}}\longrightarrow X_L\times_{\mathcal{S}}\overline{\mathcal{S}}$.
 Then, it is straightforward to check that $\{(X_L,\rho_L,\eta_{L,m})\}_{L\in\mathscr{L}}$ satisfies
 the conditions in the definition of $\M'_{\mathscr{L},m}$ (for the condition on polarizations, we can work on
 the generic fibers). Moreover, the isomorphism class of $\{(X_L,\rho_L,\eta_{L,m})\}_{L\in\mathscr{L}}$
 is independent of the choice of $\overline{x}$, $\eta$ and $L_0$.
 Hence we get a morphism $\mathcal{S}\longrightarrow \M'_{\mathscr{L},m}$,
 and thus a morphism $M_{K_{\!\mathscr{L},m}}\longrightarrow \M'^\rig_{\mathscr{L},m}$.
 By the construction, it is easy to verify that this morphism is compatible with the actions of $N_{\mathscr{L}}\times J$ and change of $m$.

 For $m'\ge m$, the transition morphism $M_{K_{\!\mathscr{L},m'}}\longrightarrow M_{K_{\!\mathscr{L},m}}$ is
 surjective. 
 Therefore, the morphism $M_{K_{\!\mathscr{L},m}}\longrightarrow \M'^\rig_{\mathscr{L},m}$
 factors through the image of $\M'^\rig_{\mathscr{L},m'}\longrightarrow \M'^\rig_{\mathscr{L},m}$
 for every $m'\ge m$.
 Since $\M'^\rig_{\mathscr{L},m'}\longrightarrow \M'^\rig_{\mathscr{L},m}$ is finite \'etale, 
 $\M^\rig_{\mathscr{L},m}$ coincides with the intersection of the images of
 $\M'^\rig_{\mathscr{L},m'}\longrightarrow \M'^\rig_{\mathscr{L},m}$ for all $m'\ge m$.
 Hence we have a morphism $M_{K_{\!\mathscr{L},m}}\longrightarrow \M^\rig_{\mathscr{L},m}$ of rigid spaces
 over $\Q_{p^\infty}$.

 Next we construct an inverse morphism $\M^\rig_{\mathscr{L},m}\longrightarrow M_{K_{\!\mathscr{L},m}}$.
 Let $\mathcal{S}$ be a formal scheme of finite type over $\Spf\Z_{p^\infty}$ such that
 $\mathcal{S}^\rig$ is connected, and $\mathcal{S}^\rig\longrightarrow \M^\rig_{\mathscr{L},m}$ a morphism
 over $\Q_{p^\infty}$. It suffices to construct a morphism $\mathcal{S}^\rig\longrightarrow M_{K_{\!\mathscr{L},m}}$. By changing $\mathcal{S}$ by its admissible blow-up, we may assume that
 the morphism $\mathcal{S}^\rig\longrightarrow \M^\rig_{\mathscr{L},m}$ comes from a morphism
 $\mathcal{S}\longrightarrow \M_{\mathscr{L},m}$ of formal schemes over $\Z_{p^\infty}$.
 Therefore, we have the data $\{(X_L,\rho_L,\eta_{L,m})\}$ where
 \begin{itemize}
  \item $X_L$ is a $p$-divisible group on $\mathcal{S}$,
  \item $\rho_L\colon \X\otimes_{\overline{\F}_p}\overline{\mathcal{S}}\longrightarrow X_L\times_{\mathcal{S}}\overline{\mathcal{S}}$ is a quasi-isogeny,
  \item and $\eta_{L,m}\colon L/p^mL\longrightarrow X_L[p^m]$
 \end{itemize}
 satisfying suitable conditions (\cf Definition \ref{defn:formal-model}).
 
 Fix a geometric point $\overline{x}$ of $\mathcal{S}^\rig$. By the definition of $\M_{\mathscr{L},m}$,
 we can find a family of isomorphisms
 $\{\eta_L\colon L\yrightarrow{\cong}T_pX^\rig_{L,\overline{x}}\}_{L\in\mathscr{L}}$ satisfying the following:
 \begin{itemize}
  \item[(a)] $\eta_L\bmod p^m=\eta_{L,m,\overline{x}}$.
  \item[(b)] For $L,L'\in\mathscr{L}$ with $L\subset L'$, the following diagrams are commutative:
	\[
	 \xymatrix{
	L\ar[d]^{\eta_{L}}_{\cong}\ar[r]& L'\ar[d]^{\eta_{L'}}_{\cong}\\
	T_pX^\rig_{L,\overline{x}}\ar[r]^-{\widetilde{\rho}_{L',L}}& T_pX^\rig_{L',\overline{x}}\lefteqn{,}
	}
	\qquad\qquad
	\xymatrix{
	L\ar[d]^{\eta_L}_{\cong}\ar[r]^-{\times p}\ar[r]_-{\cong}& pL\ar[d]^{\eta_{pL}}_{\cong}\\
	T_pX^\rig_{L,\overline{x}}\ar[r]_-{\cong}^-{\theta_p}& T_pX^\rig_{pL,\overline{x}}\lefteqn{.}
	}
	\]
  \item[(c)] Fix an isomorphism $\Z_p\yrightarrow{\cong} \Z_p(1)$.
	Then, for every $L\in \mathscr{L}$,
	the diagram below is commutative up to constant which is independent of $L$:
	\[
	 \xymatrix{%
	L\times L^\vee\ar[d]^{\eta_L\times \eta_{L^\vee}}_{\cong}\ar[rr]^-{\langle\ ,\ \rangle}&&
	\Z_p\ar[d]^{\cong}\\
	T_pX^\rig_{L,\overline{x}}\times T_pX^\rig_{L^\vee,\overline{x}}\ar[r]^-{\id\times \lambda_{L^\vee}}&
	T_pX^\rig_{L,\overline{x}}\times T_pX^{\vee\rig}_{L,\overline{x}}\ar[r]& \Z_p(1)\lefteqn{.}
	}
	\]
 \end{itemize}
 Fix $L_0\in\mathscr{L}$ such that $L_0\subset\Z_p^{2d}$.
 The homomorphism $\eta_{L_0}\colon L_0\longrightarrow T_pX^\rig_{L_0,\overline{x}}$ can be identified
 with $\eta_{L_0}\colon \Q_p^{2d}/L_0\longrightarrow X^\rig_{L_0,\overline{x}}$.
 We will see that the image $\eta_{L_0}(\Z_p^{2d}/L_0)$ of $\Z_p^{2d}/L_0$ under $\eta_{L_0}$
 is $\pi_1(\mathcal{S}^\rig,\overline{x})$-stable. Indeed, for each $\sigma\in \pi_1(\mathcal{S}^\rig,\overline{x})$ and $L\in\mathscr{L}$, we can find $g_{L,\sigma}\in \GL(L)$ such that $\sigma\circ\eta_L=\eta_L\circ g_{L,\sigma}$. By (b) above, if we regard $\GL(L)$ as a subset of $\GL_{2d}(\Q_p)$, then $g_{L,\sigma}=g_{L_0,\sigma}$
 for every $L\in\mathscr{L}$. Moreover, (a) tells us that $g_{L_0,\sigma}=g_{L,\sigma}$ lies in
 the kernel of $\GL(L)\longrightarrow \GL(L/p^mL)$,
 and (c) tells us that $g_{L_0,\sigma}\in\GSp_{2d}(\Q_p)$.
 We conclude that $g_{L_0,\sigma}\in K_{\!\mathscr{L},m}$, and thus $g_{L_0,\sigma}\in K_0$
 by the assumption. Hence $\sigma(\eta_{L_0}(\Z_p^{2d}/L_0))=\eta_{L_0}(g_{L_0,\sigma}(\Z_p^{2d}/L_0))=\eta_{L_0}(\Z_p^{2d}/L_0)$, as desired.
 Therefore, the subset $\eta_{L_0}(\Z_p^{2d}/L_0)\subset T_pX^\rig_{L_0,\overline{x}}$
 corresponds to a finite \'etale closed subgroup scheme of $X^\rig_{L_0}$.
 By the flattening theorem, after replacing $\mathcal{S}$, we may assume that this subgroup scheme
 comes from a finite flat closed subgroup scheme $Y$ of $X_{L_0}$. Put $X=X_{L_0}/Y$.
 Then, by the similar argument as above, we have an isomorphism 
 $\eta\colon \Z_p^{2d}\yrightarrow{\cong}T_pX^\rig_{\overline{x}}$ which makes the following diagram
 commutative:
 \[
  \xymatrix{%
 L_0\ar[r]\ar[d]^{\eta_{L_0}}_{\cong}& \Z_p^{2d}\ar[d]^{\eta}_{\cong}\\
 T_pX^\rig_{L_0,\overline{x}}\ar[r]& T_pX^\rig_{\overline{x}}\lefteqn{.}
 }
 \]
 Let $\rho$ be the composite $\X\otimes_{\overline{\F}_p}\overline{\mathcal{S}}\yrightarrow{\rho_{L_0}}X_{L_0}\times_{\mathcal{S}}\overline{\mathcal{S}}\longrightarrow X\times_{\mathcal{S}}\overline{\mathcal{S}}$.
 We will show that the triple $(X,\rho,\overline{\eta})$ gives an $\mathcal{S}^\rig$-valued point of $M_{K_{\!\mathscr{L},m}}$. For existence of a polarization, let $\lambda\colon X\longrightarrow X^\vee$ be the quasi-isogeny
 $X\longrightarrow X_L\yrightarrow{\lambda_L} (X_{L^\vee})^\vee\yrightarrow{(*)}X_L^\vee\longrightarrow X^\vee$, where $(*)$ is the quasi-isogeny $\widetilde{\rho}^{\vee}_{L^\vee,L}$ if $L\subset L^\vee$,
 and $(\widetilde{\rho}^\vee_{L,L^\vee})^{-1}$ if $L^\vee\subset L$.
Then, the following diagram is commutative up to multiplication by $\Z_p^\times$:
 \[
 \xymatrix{%
 \Q_p^{2d}\times \Q_p^{2d}\ar[d]^{\eta\times \eta}_{\cong}\ar[rr]^-{\langle\ ,\ \rangle}&&
 \Q_p\ar[d]^{\cong}\\
 V_pX^\rig_{\overline{x}}\times V_pX^\rig_{\overline{x}}\ar[r]^-{\id\times \lambda}&
 V_pX^\rig_{\overline{x}}\times V_pX^{\vee\rig}_{\overline{x}}\ar[r]& \Q_p(1)\lefteqn{,}
 }
 \]
 where the isomorphism $\Q_p\longrightarrow \Q_p(1)$ is induced from the isomorphism
 $\Z_p\longrightarrow \Z_p(1)$ fixed in (c) above.
 In particular $\lambda(T_pX^\rig_{\overline{x}})=T_pX^{\vee\rig}_{\overline{x}}$, and thus $\lambda$
 is an isomorphism. It obviously satisfies that $\rho^\vee\circ\lambda\circ\rho=a\lambda_0$ for some $a\in\Q_p^\times$.
 The diagram above also tells us that $\eta$ preserves polarizations up to multiplication by $\Z_p^\times$.
 On the other hand, the $K_{\!\mathscr{L},m}$-orbit of $\eta$ is invariant under $\pi_1(\mathcal{S}^\rig,\overline{x})$ simply because the $K_{\!\mathscr{L},m}$-orbit of $\eta_{L_0}$ is invariant
 (recall that we have proved $g_{L_0,\sigma}\in K_{\!\mathscr{L},m}$).
 
 Now we obtain an $\mathcal{S}^\rig$-valued point $(X,\rho,\overline{\eta})$ of $M_{K_{\!\mathscr{L},m}}$.
 It is easy to see that it is independent of the choice of $\overline{x}$ and $L_0$.
 Hence we have a canonical morphism $\mathcal{S}^\rig\longrightarrow M_{K_{\!\mathscr{L},m}}$,
 and thus $\M_{\mathscr{L},m}^\rig\longrightarrow M_{K_{\!\mathscr{L},m}}$.

 It is not difficult to show that the morphisms $M_{K_{\!\mathscr{L},m}}\longrightarrow \M_{\mathscr{L},m}^\rig$
 and $\M_{\mathscr{L},m}^\rig\longrightarrow M_{K_{\!\mathscr{L},m}}$ we have obtained
 are inverse to each other. As $\M_{\mathscr{L},m}^\rig\cong \M_{\mathscr{L},m}^{\flat,\rig}$,
 we have a desired isomorphism $M_{K_{\!\mathscr{L},m}}\cong \M_{\mathscr{L},m}^{\flat,\rig}$.
\end{prf}

\subsection{Boundary strata}
\begin{defn}\label{defn:direct-summand}
 For an integer $h$ with $1\le h\le d$,
 let $\mathcal{S}_{\infty,h}$ be the set of totally isotropic subspaces of dimension $h$ of $\Q_p^{2d}$.
 Recall that a subspace $V$ of $\Q_p^{2d}$ is said to be totally isotropic if $V\subset V^\perp$.
 For a lattice $L$ of $\Z_p^{2d}$ and an integer $m\ge 0$, 
 denote by $\mathcal{S}_{L,m,h}$ the set of direct summands of rank $h$ of $L/p^mL$.

 For a self-dual chain of lattices $\mathscr{L}$ of $\Q_p^{2d}$ and an integer $m\ge 0$, we have a natural map
 $\mathcal{S}_{\infty,h}\longrightarrow \prod_{L\in\mathscr{L}}\mathcal{S}_{L,m,h}$;
 $V\longmapsto V_{\!\mathscr{L},m}=(V\cap L\bmod p^m)_{L\in\mathscr{L}}$.
 We denote its image by $\mathcal{S}_{\mathscr{L},m,h}$.
 The group $N_{\mathscr{L}}$ acts naturally on $\mathcal{S}_{\infty,h}$ and $\prod_{L\in\mathscr{L}}\mathcal{S}_{L,m,h}$,
 and the map above is equivariant with respect to these actions.
 Therefore $N_{\mathscr{L}}$ also acts on $\mathcal{S}_{\mathscr{L},m,h}$.

 Put $\mathcal{S}_\infty=\bigcup_{h=1}^d\mathcal{S}_{\infty,h}$ and $\mathcal{S}_{\mathscr{L},m}=\bigcup_{h=1}^d\mathcal{S}_{\mathscr{L},m,h}$.
 For two elements $\alpha=(\alpha_L)_L$, $\beta=(\beta_L)_L$ in $\mathcal{S}_{\mathscr{L},m}$,
 we write $\alpha\prec \beta$ if $\alpha_L\supset \beta_L$ for every $L\in\mathscr{L}$.
 This gives a partial order on $\mathcal{S}_{\mathscr{L},m}$.
 The action of $N_{\mathscr{L}}$ on $\mathcal{S}_{\mathscr{L},m}$ obviously preserves this partial order.
\end{defn}

\begin{lem}\label{lem:S-finite-set}
 For every $m\ge 0$, $\mathcal{S}_{\mathscr{L},m}$ is a finite set.
\end{lem}

\begin{prf}
 Since $K_{\!\mathscr{L},m}$ acts trivially on $\mathcal{S}_{\mathscr{L},m}$,
 we have a natural surjection $K_{\!\mathscr{L},m}\backslash\mathcal{S}_\infty\longrightarrow \mathcal{S}_{\mathscr{L},m}$. Therefore it suffices to show that $K_{\!\mathscr{L},m}\backslash\mathcal{S}_{\infty,h}$ is a finite set
 for each $1\le h\le d$. For an integer $m'\ge 0$, let $K_{m'}$ be the kernel of
 $K_0\longrightarrow \GSp_{2d}(\Z/p^{m'}\Z)$.
 Then, for a sufficiently large $m'$, we have $K_{m'}\subset K_{\!\mathscr{L},m}$.
 Therefore it suffices to show that $K_{m'}\backslash\mathcal{S}_{\infty,h}$ is a finite set.
 Fix an element $V$ of $\mathcal{S}_{\infty,h}$, and denote by $P$ the stabilizer of $V$ in $G$.
 It is well-known that $G$ acts transitively on $\mathcal{S}_{\infty,h}$ and $P$ is a parabolic subgroup of $G$.
 Therefore, by the Iwasawa decomposition $G=K_0P$ for the hyperspecial subgroup $K_0$, we have
 \[
  K_{m'}\backslash\mathcal{S}_{\infty,h}\cong K_{m'}\backslash G/P\cong K_{m'}\backslash K_0/P\cap K_0,
 \]
 which is obviously a finite set.
\end{prf}

Put $\M^\flat_{\mathscr{L},m,s}=\M^\flat_{\mathscr{L},m}\otimes_{\Z_{p^\infty}}\overline{\F}_p$.
To $\alpha\in\mathcal{S}_{\mathscr{L},m}$, we will attach a closed formal subscheme $\M^\flat_{\mathscr{L},m,\alpha}$
of $\M^\flat_{\mathscr{L},m,s}$.
A family $\{\M^\flat_{\mathscr{L},m,\alpha}\}_{\alpha\in\mathcal{S}_{\mathscr{L},m}}$ of these formal subschemes
plays a role of the ``boundary'' of $M_{K_{\!\mathscr{L},m}}$ inside $\M^\flat_{\mathscr{L},m}$.

\begin{defn}\label{defn:boundary-strata}
 For $\alpha=(\alpha_L)_L\in\mathcal{S}_{\mathscr{L},m}$, we define the subfunctor $\M^\flat_{\mathscr{L},m,\alpha}$
 of $\M^\flat_{\mathscr{L},m,s}$ as follows:
 for an $\overline{\F}_p$-scheme $S$, $\{(X_L,\rho_L,\eta_{L,m})\}_{L\in\mathscr{L}}\in\M^\flat_{\mathscr{L},m}(S)$ lies in
 $\M^\flat_{\mathscr{L},m,\alpha}(S)$ if and only if $\eta_{L,m}\vert_{\alpha_L^\perp}=0$ for every $L\in\mathscr{L}$.
\end{defn}

\begin{lem}\label{lem:boundary-strata}
 The subfunctor $\M^\flat_{\mathscr{L},m,\alpha}$ is represented by a closed formal subscheme of $\M^\flat_{\mathscr{L},m,s}$.
\end{lem}

\begin{prf}
 Obviously, for each $L\in \mathscr{L}$ the condition $\eta_{L,m}\vert_{\alpha_L^\perp}=0$ gives
 a closed formal subscheme of $\M^\flat_{\mathscr{L},m,s}$. Therefore it suffices to show the equivalence of
 $\eta_{L,m}\vert_{\alpha_L^\perp}=0$ and $\eta_{pL,m}\vert_{\alpha_{pL}^\perp}=0$.
 Since $\alpha\in\mathcal{S}_{\mathscr{L},m}$, we have $\alpha_{pL}=p\alpha_L$.
 Therefore the commutative diagram 
\[
 \xymatrix{%
 L/p^mL\ar[d]^{\eta_{L,m}}\ar[r]^-{\times p}\ar[r]_-{\cong}& pL/p^{m+1}L\ar[d]^{\eta_{pL,m}}\\
 X_L[p^m]\ar[r]_-{\cong}^-{\theta_p}& X_{pL}[p^m]
 }
\]
 in the definition of $\M'_{\mathscr{L},m}$ gives the equivalence.
\end{prf}

The following lemma is clear from the definition:

\begin{lem}\label{lem:formal-model-group-action}
 \begin{enumerate}
  \item The action of $J$ on $\M^\flat_{\mathscr{L},m}$ preserves the closed formal subscheme 
	$\M^\flat_{\mathscr{L},m,\alpha}$ for each $\alpha\in\mathcal{S}_{\mathscr{L},m}$.
  \item For $g\in N_{\mathscr{L}}$, 
	the right action $g\colon \M^\flat_{\mathscr{L},m}\longrightarrow \M^\flat_{\mathscr{L},m}$
	induces an isomorphism $\M^\flat_{\mathscr{L},m,\alpha}\longrightarrow \M^\flat_{\mathscr{L},m,g^{-1}\alpha}$ for each $\alpha\in\mathcal{S}_{\mathscr{L},m}$.
 \end{enumerate}
\end{lem}

For $\alpha\in\mathcal{S}_{\mathscr{L},m}$, put $M^\flat_{\mathscr{L},m,\alpha}=t(\M^\flat_{\mathscr{L},m,\alpha})_a$.
It is a closed analytic adic subspace of $M^\flat_{\mathscr{L},m,s}=t(\M^\flat_{\mathscr{L},m,s})_a$.
Moreover, set $M^\flat_{\mathscr{L},m,(\alpha)}=M^\flat_{\mathscr{L},m,\alpha}\setminus \bigcup_{\beta\succ\alpha}M^\flat_{\mathscr{L},m,\beta}$. It is a locally closed subset of $M^\flat_{\mathscr{L},m,s}$.

\begin{prop}\label{prop:formal-model-disjoint}
 \begin{enumerate}
  \item For $\alpha,\beta\in \mathcal{S}_{\mathscr{L},m}$ with $\alpha\neq\beta$, 
	$M^\flat_{\mathscr{L},m,(\alpha)}\cap M^\flat_{\mathscr{L},m,(\beta)}=\varnothing$.
  \item Assume that $d=2$. Then, we have $M^\flat_{\mathscr{L},m,s}=\bigcup_{\alpha\in\mathcal{S}_{\mathscr{L},m}}M^\flat_{\mathscr{L},m,(\alpha)}$.
 \end{enumerate}
\end{prop}

\begin{prf}
 We use the $p$-adic uniformization theorem. Fix a compact open subgroup $K^p$ of $\GSp_{2d}(\A_f^{p})$.
 By \cite[Lemma 4.1]{RZ-GSp4}, we may assume that $(\X,\lambda_0)$ comes from a polarized $d$-dimensional
 abelian variety.
 Let $\Sh_{\mathscr{L},K^p}$ be the moduli space over $\Z_{p^\infty}$ of
 ``$\mathscr{L}$-sets'' of abelian varieties with principal polarizations and $K^p$-level structures 
 introduced in \cite[Definition 6.9]{MR1393439}.
 We shrink $K^p$ so that $\Sh_{\mathscr{L},K^p}$ is represented by a quasi-projective scheme over $\Z_{p^\infty}$.
 For simplicity, we write $\Sh_{\mathscr{L}}$ for $\Sh_{\mathscr{L},K^p}$.
 By using Drinfeld level structures, we can construct towers $\{\Sh'_{\mathscr{L},m}\}_{m\ge 0}$, 
 $\{\Sh_{\mathscr{L},m}\}_{m\ge 0}$ and $\{\Sh^\flat_{\mathscr{L},m}\}_{m\ge 0}$ over $\Sh_{\mathscr{L}}$  
 similarly as in Definition \ref{defn:formal-model}. 
 These towers are endowed with right actions of $N_{\mathscr{L}}$.
 Let $Y^\flat_m$ be a closed subset of $\Sh^\flat_{\mathscr{L},m}$ consisting of supersingular points,
 and $(\Sh^\flat_{\mathscr{L},m})^\wedge_{/Y^\flat_m}$ the formal completion of $\Sh^\flat_{\mathscr{L},m}$ along $Y^\flat_m$.
 The $p$-adic uniformization theorem \cite[Theorem 6.30]{MR1393439}
 tells us that we have a natural $N_{\mathscr{L}}$-equivariant \'etale surjection
 \[
  \theta_m\colon \M^\flat_{\mathscr{L},m}\longrightarrow (\Sh^\flat_{\mathscr{L},m})^\wedge_{/Y^\flat_m}
 \]
 (note that $\GSp_{2d}$ satisfies the Hasse principle; \cf \cite[\S 7]{MR1124982}).

 Put $\Sh^\flat_{\mathscr{L},m,s}=\Sh_{\mathscr{L},m}^\flat\otimes_{\Z_{p^\infty}}\overline{\F}_p$.
 For $\alpha\in \mathcal{S}_{\mathscr{L},m}$, we can define a closed subscheme 
 $\Sh^\flat_{\mathscr{L},m,\alpha}$ of $\Sh^\flat_{\mathscr{L},m,s}$ in the same way as
 $\M^\flat_{\mathscr{L},m,\alpha}$.
 By the construction of $\theta_m$, $\M^\flat_{\mathscr{L},m,\alpha}$ is isomorphic to the fiber product
 $\M^\flat_{\mathscr{L},m,s}\times_{\Sh^\flat_{\mathscr{L},m,s}}\Sh^\flat_{\mathscr{L},m,\alpha}$.
 Thus, as in \cite[Example 4.2]{adicLTF}, we can reduce our problem to an analogue
 for $\{\Sh^\flat_{\mathscr{L},m,\alpha}\}$; namely, for
 $\Sh^\flat_{\mathscr{L},m,(\alpha)}=\Sh^\flat_{\mathscr{L},m,\alpha}\setminus \bigcup_{\beta\succ\alpha}\Sh^\flat_{\mathscr{L},m,\beta}$, it suffices to show the following:
 \begin{enumerate}
  \item[(I)] For $\alpha,\beta\in \mathcal{S}_{\mathscr{L},m}$ with $\alpha\neq\beta$, 
	     $\Sh^\flat_{\mathscr{L},m,(\alpha)}\cap \Sh^\flat_{\mathscr{L},m,(\beta)}=\varnothing$. 
  \item[(II)] Under the assumption $d=2$, we have $\Sh^\flat_{\mathscr{L},m,s}\setminus Y^\flat_m=\bigcup_{\alpha\in\mathcal{S}_{\mathscr{L},m}}\Sh^\flat_{\mathscr{L},m,(\alpha)}$.
 \end{enumerate}
 Let $x=\{(A_L,\eta_{L,m})\}_{L\in\mathscr{L}}\in \Sh^\flat_{\mathscr{L},m,s}(\overline{\F}_p)$.
 We shall prove that $((\Ker\eta_{L,m})^\perp)_L$ lies in $\mathcal{S}_{\mathscr{L},m}$
 provided that $\eta_{L,m}\neq 0$ for some $L$.
 Since $\Sh^\flat_{\mathscr{L},m}$ is flat over $\Z_{p^\infty}$, $x$ can be lifted to a point
 on the generic fiber. Therefore, we can find
 \begin{itemize}
  \item a finite extension $F$ of $\Q_{p^\infty}$,
  \item an ``$\mathscr{L}$-set'' of $p$-divisible groups $\{X_L\}_{L\in\mathscr{L}}$ with a principal
	polarization $\{\lambda_L\}_{L\in\mathscr{L}}$ (\cf \cite[Definition 6.5, Definition 6.6]{MR1393439})
	whose special fiber can be identified with $\{A_L[p^\infty]\}_{L\in\mathscr{L}}$ (with the implicit polarization),
  \item and a family of isomorphisms $\{\widetilde{\eta}_L\colon L\yrightarrow{\cong} T_pX_{L,\overline{\eta}}\}_{L\in\mathscr{L}}$ (here $X_{L,\overline{\eta}}$ denotes the geometric generic fiber of $X_L$) satisfying the following conditions:
	\begin{itemize}
	 \item[(a)] Denote by $\eta_L$ the composite 
		    $L\yrightarrow{\widetilde{\eta}_L}T_pX_{L,\overline{\eta}}\longrightarrow T_pX_{L,s}$,
		    where $X_{L,s}$ denotes the special fiber of $X_L$.
		    Then, $\eta_L\bmod p^m$ coincides with $\eta_{L,m}$ under the identification
		    of $A_L[p^\infty]$ and $X_{L,s}$. 
	 \item[(b)] For $L,L'\in\mathscr{L}$ with $L\subset L'$, the following diagrams are commutative:
		    \[
		    \xymatrix{
		    L\ar[d]^{\widetilde{\eta}_{L}}_{\cong}\ar[r]& L'\ar[d]^{\widetilde{\eta}_{L'}}_{\cong}\\
		    T_pX_{L,\overline{\eta}}\ar[r]& T_pX_{L',\overline{\eta}}\lefteqn{,}
		    }
		    \qquad\qquad
		    \xymatrix{
		    L\ar[d]^{\widetilde{\eta}_L}_{\cong}\ar[r]^-{\times p}\ar[r]_-{\cong}& pL\ar[d]^{\widetilde{\eta}_{pL}}_{\cong}\\
		    T_pX_{L,\overline{\eta}}\ar[r]_-{\cong}^-{\theta_p}& T_pX_{pL,\overline{\eta}}\lefteqn{.}
		    }
		    \]
	 \item[(c)] Fix an isomorphism $\Z_p\yrightarrow{\cong} \Z_p(1)$.
		    Then, for every $L\in \mathscr{L}$,
		    the diagram below is commutative up to constant which is independent of $L$:
		    \[
		    \xymatrix{%
		    L\times L^\vee\ar[d]^{\widetilde{\eta}_L\times \widetilde{\eta}_{L^\vee}}_{\cong}\ar[rr]^-{\langle\ ,\ \rangle}&&
		    \Z_p\ar[d]^{\cong}\\
		    T_pX_{L,\overline{\eta}}\times T_pX_{L^\vee,\overline{\eta}}\ar[r]^-{\id\times \lambda_{L^\vee}}&
		    T_pX_{L,\overline{\eta}}\times T_pX^\vee_{L,\overline{\eta}}\ar[r]& \Z_p(1)\lefteqn{.}
		    }
		    \]
	\end{itemize}
 \end{itemize}
 Let $V$ be the kernel of $\eta_L\otimes_{\Z_p}\Q_p\colon \Q_p^{2d}=L\otimes_{\Z_p}\Q_p\longrightarrow V_pX_{L,s}$.
 By the first diagram in (b) above, $V$ is independent of the choice of $L$.
 Moreover, we have $\Ker\eta_L=V\cap L$,
 and thus $\Ker\eta_{L,m}=V\cap L\bmod p^m$ by (a).
 In particular $V\neq \Q_p^{2d}$, as we are assuming that $\eta_{L,m}\neq 0$ for some $L$.
 Therefore, it suffices to prove that
 $V\subset \Q_p^{2d}$ is coisotropic, that is to say, satisfies $V^\perp\subset V$.
 Denote by $\lambda$ the quasi-isogeny $X_L\yrightarrow{(*)}X_{L^\vee}\yrightarrow{\lambda_{L^\vee}}X^\vee_L$,
 where $(*)$ comes from the $\mathscr{L}$-set structure of $\{X_L\}_{L\in\mathscr{L}}$.
 This induces an alternating pairing
 \[
  \langle\ ,\ \rangle_{\lambda}\colon V_pX_{L,\overline{\eta}}\times V_pX_{L,\overline{\eta}}\yrightarrow{\id\times\lambda}V_pX_{L,\overline{\eta}}\times V_pX^\vee_{L,\overline{\eta}}\longrightarrow \Q_p(1).
 \]
 By the commutative diagrams in (b) and (c) above, $\widetilde{\eta}_L\otimes_{\Z_p}\Q_p\colon \Q_p^{2d}=L\otimes_{\Z_p}\Q_p\yrightarrow{\cong} V_pX_{L,\overline{\eta}}$ maps the pairing $\langle\ ,\ \rangle$ to a scalar multiple of
 $\langle\ ,\ \rangle_\lambda$ under an identification of $\Q_p\cong \Q_p(1)$.
 Thus, it suffices to show that the kernel of the specialization map
 $V_pX_{L,\overline{\eta}}\longrightarrow V_pX_{L,s}$ is a coisotropic subspace of $V_pX_{L,\overline{\eta}}$
 with respect to $\langle\ ,\ \rangle_\lambda$.
 To prove it, we can argue in the same way as in the proof of \cite[Lemma 5.8]{RZ-GSp4}. 
 Let us recall the argument briefly. 
 Take an exact sequence of $p$-divisible groups
 $0\longrightarrow X^0_{L,s}\longrightarrow X_{L,s}\longrightarrow X_{L,s}^{\et}\longrightarrow 0$,
 where $X^0_{L,s}$ is connected and $X_{L,s}^{\et}$ is \'etale.
 It is canonically lifted to an exact sequence $0\longrightarrow X^0_L\longrightarrow X_L\longrightarrow X_L^{\et}\longrightarrow 0$ over $\mathcal{O}_F$, where $X_L^\et$ is an \'etale $p$-divisible group 
 (\cf \cite[p.~76]{MR0347836}). It is easy to see that the kernel of $V_pX_{L,\overline{\eta}}\longrightarrow V_pX_{L,s}$ 
 coincides with $V_pX^0_{L,\overline{\eta}}$.
 Therefore it suffices to prove that the composite of Galois-equivariant homomorphisms
 $(V_pX^0_{L,\overline{\eta}})^\perp\hooklongrightarrow V_pX_{L,\overline{\eta}}\longrightarrow V_pX^\et_{L,\overline{\eta}}$ is zero.
 This follows from the $p$-adic Hodge theory, noting that $\lambda\colon V_pX_{L,\overline{\eta}}\yrightarrow{\cong}V_pX^\vee_{L,\overline{\eta}}=(V_pX_{L,\overline{\eta}})^\vee(1)$ induces a Galois-equivariant
 isomorphism $(V_pX^0_{L,\overline{\eta}})^\perp\cong (V_pX^\et_{L,\overline{\eta}})^\vee(1)$.

 Now we can show (I) and (II) above. For (I), assume that $x=\{(A_L,\eta_{L,m})\}_{L\in\mathscr{L}}\in\Sh^\flat_{\mathscr{L},m,s}(\overline{\F}_p)$ lies in $\Sh^\flat_{\mathscr{L},m,(\alpha)}$.
 Then $\alpha_L^\perp\subset \Ker\eta_{L,m}$ for every $L\in\mathscr{L}$, or equivalently, 
 $\alpha\prec ((\Ker\eta_{L,m})^\perp)_L$. As $x$ belongs to $\Sh^\flat_{\mathscr{L},m,((\Ker\eta_{L,m})^\perp)_L}(\overline{\F}_p)$, we conclude that $\alpha=((\Ker\eta_{L,m})^\perp)_L$.
 In particular, for $x\in\Sh^\flat_{\mathscr{L},m,s}(\overline{\F}_p)$,
 there is at most one $\alpha\in\mathcal{S}_{\mathscr{L},m}$
 with $x\in \Sh^\flat_{\mathscr{L},m,(\alpha)}(\overline{\F}_p)$. This concludes the proof of (I).
 For (II), assume that $\eta_{L,m}=0$ for every $L\in\mathscr{L}$. By \cite[Lemma II.2.1]{MR1876802},
 $A_L[p^\infty]$ has no \'etale part. Since the rational Dieudonn\'e module $D(A_L[p^\infty])_\Q$ is
 polarized, \cite[Lemma 4.1]{RZ-GSp4} tells us that, under the assumption $d=2$,
 $A_L$ is supersingular, namely, $x\in Y^\flat_m(\overline{\F}_p)$.
 Hence we conclude that $\Sh^\flat_{\mathscr{L},m,s}\setminus Y^\flat_m=\bigcup_{\alpha\in\mathcal{S}_{\mathscr{L},m}}\Sh^\flat_{\mathscr{L},m,(\alpha)}$, as desired.
\end{prf}

\begin{rem}\label{rem:locally-algebraizable}
 For every quasi-compact open formal subscheme $\mathscr{U}$ of $\M^\flat_{\mathscr{L},m}$,
 we can choose the level $K^p\subset \GSp_{2d}(\A^p_f)$ such that $\theta_m\vert_{\mathscr{U}}$ is
 an isomorphism from $\mathscr{U}$ onto the image (\cf \cite[Corollaire 3.1.4]{MR2074714}).
 In particular, the formal scheme $\M^\flat_{\mathscr{L},m}$ is locally algebraizable
 in the sense of \cite[Definition 3.18]{formalnearby}.
\end{rem}

\subsection{Group action on $\mathcal{S}_{\mathscr{L},m}$}
In this subsection, we assume that $d=2$. In this case, by the Bruhat-Tits theory,
we can easily classify self-dual chains of lattices:

\begin{lem}\label{lem:lattice-chain-GSp4}
 For every self-dual chain of lattices $\mathscr{L}$, we can find $g\in G$
 so that $g\mathscr{L}$ is one of the following:
\begin{description}
 \item[(hyperspecial type)] $\mathscr{L}_0=\{p^m\Z_p^4\mid m\in \Z\}$.
 \item[(paramodular type)] 
	    \[
	    \mathscr{L}_{\mathrm{para}}=\{p^{m-1}\Z_p\oplus p^m\Z_p\oplus p^m\Z_p\oplus p^m\Z_p, p^m\Z_p\oplus p^m\Z_p\oplus p^m\Z_p\oplus p^{m+1}\Z_p\mid m\in \Z\}.
	    \]
 \item[(Siegel parahoric type)]
	    \[
	     \mathscr{L}_{\mathrm{Siegel}}=\mathscr{L}_0\cup \{p^m\Z_p\oplus p^m\Z_p\oplus p^{m+1}\Z_p\oplus p^{m+1}\Z_p\mid m\in \Z\}.
	    \]
 \item[(Klingen parahoric type)] $\mathscr{L}_{\mathrm{Klingen}}=\mathscr{L}_0\cup \mathscr{L}_{\mathrm{para}}$.
 \item[(Iwahori type)] $\mathscr{L}_{\mathrm{Iw}}=\mathscr{L}_{\mathrm{Siegel}}\cup \mathscr{L}_{\mathrm{para}}$.
\end{description}
\end{lem}

\begin{lem}\label{lem:max-cpt-mod-center}
 Every maximal compact-mod-center subgroup of $G$ is conjugate to one of 
 $N_{\mathscr{L}_0}$, $N_{\mathscr{L}_{\mathrm{para}}}$ or $N_{\mathscr{L}_{\mathrm{Siegel}}}$.
\end{lem}

\begin{prf}
 Since a compact-mod-center subgroup is contained in the normalizer of some parahoric subgroup, it suffices to show that
 $N_{\mathscr{L}_{\mathrm{Klingen}}}=N_{\mathscr{L}_0}\cap N_{\mathscr{L}_{\mathrm{para}}}$ and
 $N_{\mathscr{L}_{\mathrm{Iw}}}=N_{\mathscr{L}_{\mathrm{Siegel}}}\cap N_{\mathscr{L}_{\mathrm{para}}}$.

 For a self-dual chain of lattices $\mathscr{L}$, let $\mathscr{L}_+$ be the subset of $\mathscr{L}$ consisting of lattices $L$
 such that $L^\vee=aL$ for some $a\in\Q_p$, and put $\mathscr{L}_-=\mathscr{L}\setminus \mathscr{L}_+$.
 Then, clearly $g\in N_{\mathscr{L}}$ preserves $\mathscr{L}_+$ and $\mathscr{L}_-$.
 If $\mathscr{L}=\mathscr{L}_{\mathrm{Klingen}}$, then $\mathscr{L}_+=\mathscr{L}_0$ and $\mathscr{L}_-=\mathscr{L}_{\mathrm{para}}$.
 This implies that $N_{\mathscr{L}_{\mathrm{Klingen}}}$ is contained in $N_{\mathscr{L}_0}\cap N_{\mathscr{L}_{\mathrm{para}}}$.
 The other inclusion is obvious.
 If $\mathscr{L}=\mathscr{L}_{\mathrm{Iw}}$, then we have $\mathscr{L}_+=\mathscr{L}_{\mathrm{Siegel}}$ and $\mathscr{L}_-=\mathscr{L}_{\mathrm{para}}$,
 which gives $N_{\mathscr{L}_{\mathrm{Iw}}}=N_{\mathscr{L}_{\mathrm{Siegel}}}\cap N_{\mathscr{L}_{\mathrm{para}}}$.
\end{prf}

The following proposition is very important to control the group action on the ``boundary strata'' of the Rapoport-Zink spaces introduced in the previous subsection:

\begin{prop}\label{prop:linear-algebra}
 Let $\mathscr{L}$ be one of $\mathscr{L}_0$, $\mathscr{L}_{\mathrm{para}}$ or $\mathscr{L}_{\mathrm{Siegel}}$.
 Then, for $1\le h\le 2$ and $m\ge 1$, the natural surjection $\mathcal{S}_{\infty,h}\longrightarrow \mathcal{S}_{\mathscr{L},m,h}$ induces a bijection
 $K_{\!\mathscr{L},m}\backslash\mathcal{S}_{\infty,h}\longrightarrow \mathcal{S}_{\mathscr{L},m,h}$.
\end{prop}

Before proving it, we give its immediate corollary:

\begin{cor}\label{cor:linear-algebra}
 Let $\mathscr{L}$ be as in Proposition \ref{prop:linear-algebra}. For $g\in N_{\mathscr{L}}$ and an integer $m\ge 1$, assume that
 $gK_{\!\mathscr{L},m}=K_{\!\mathscr{L},m}g$ consists of regular elliptic elements. Then the action of $g$ on $\mathcal{S}_{\mathscr{L},m}$ has
 no fixed point.
\end{cor}

\begin{prf}
 Assume that $g$ has a fixed point $V_{\!\mathscr{L},m}\in \mathcal{S}_{\mathscr{L},m}$ with $V\in \mathcal{S}_{\infty}$.
 Then Proposition \ref{prop:linear-algebra} says that $K_{\!\mathscr{L},m}g$ intersects the stabilizer $P$ of $V$ in $G$.
 Since $P$ is a proper parabolic subgroup of $G$, it has no regular elliptic element
 (Lemma \ref{lem:proper-parab-non-elliptic}). This contradicts to the assumption.
\end{prf}

We will give a case-by-case proof of Proposition \ref{prop:linear-algebra}.

\subsubsection{Hyperspecial case}
Here we consider the case $\mathscr{L}=\mathscr{L}_0$. For $V,V'\in \mathcal{S}_{\infty,h}$,
put $I=V\cap \Z_p^4$ and $I'=V'\cap \Z_p^4$. Under the assumption $I\equiv I'\pmod{p^m}$,
we will prove that there exists $g\in K_{\!\mathscr{L}_0,m}$ with $I'=gI$.

First consider the case where $h=2$. Fix a basis $x,y\in I$.
Since the lattice $\Z_p^4$ is self-dual,
we can find $z,w\in \Z_p^4$ such that $\langle x,z\rangle=\langle y,w\rangle=1$,
$\langle x,w\rangle=\langle y,z\rangle=0$ (note that $\Hom_{\Z_p}(\Z_p^4,\Z_p)\longrightarrow \Hom_{\Z_p}(I,\Z_p)$
is surjective, as $I\subset \Z_p^4$ is a direct summand).
Replacing $w$ by $w+\langle z,w\rangle x$, we may assume that $\langle z,w\rangle=0$.
It is easy to see that $x,y,z,w$ spans a self-dual $\Z_p$-lattice of $\Q_p^4$ contained in $\Z_p^4$.
Therefore $x,y,z,w$ form a basis of $\Z_p^4$.
Take $x',y'\in I'$ such that $x\equiv x', y\equiv y'\pmod{p^m}$.
Since $x',y',z,w$ form a basis of $\Z_p^4$, we can take $A,B,C,D\in \Z_p$ so that
$z'=Az+Bw$ and $w'=Cz+Dw$ satisfy 
$\langle x',z'\rangle=\langle y',w'\rangle=1$ and $\langle x',w'\rangle=\langle y',z'\rangle=0$
(note that we have $\langle z',w'\rangle=0$ automatically). 
These are equivalent to the following identity of matrices:
\[
 \begin{pmatrix}A&B\\C&D\end{pmatrix}\begin{pmatrix}\langle x',z\rangle&\langle y',z\rangle\\\langle x',w\rangle&\langle y',w\rangle\end{pmatrix}
=\begin{pmatrix}1&0\\0&1\end{pmatrix}.
\]
By using the fact $x'-x, y'-y\in p^m\Z_p^4$, it is immediate to observe that
   \[
  \begin{pmatrix}\langle x',z\rangle&\langle y',z\rangle\\\langle x',w\rangle&\langle y',w\rangle\end{pmatrix}\in 1+p^m M_2(\Z_p).
  \]
 Therefore we can conclude that $A,D\in 1+p^m\Z_p$ and $B,C\in p^m\Z_p$. 
 In other words, we have $z'\equiv z,w'\equiv w\pmod{p^m}$.
 In particular, $x',y',z',w'$ form a basis of $\Z_p^4$.
 Let $g$ be the automorphism of $\Z_p^4$ that maps $x$ to $x'$, $y$ to $y'$, $z$ to $z'$ and $w$ to $w'$.
 Clearly $g$ is an element of $K_{\!\mathscr{L}_0,m}$ satisfying $gI=I'$.

Next consider the case where $h=1$. Take a basis $x$ of $I$ and a basis $x'$ of $I'$ such that
$x\equiv x'\pmod{p^m}$.
Put $I^\perp=V^\perp\cap \Z_p^4$ and
$I'^\perp=V'^\perp\cap \Z_p^4$, which are direct summands of rank $3$ of $\Z_p^4$
satisfying $I^\perp\equiv I'^\perp \pmod{p^m}$.
There exists $y\in I^\perp$ such that $x$ and $y$ span
a totally isotropic direct summand $I_1$ of rank $2$ of $\Z_p^4$.
Take an arbitrary element $y'\in I'^\perp$ such that $y\equiv y'\pmod{p^m}$.
Then $x'$ and $y'$ span a totally isotropic direct summand $I'_1$ of $\Z_p^4$.
By the argument in the case $h=2$, there exists $g\in K_{\!\mathscr{L}_0,m}$ satisfying
$g(x)=x'$ and $g(y)=y'$. In particular, we have $gI=I'$.

\subsubsection{Paramodular case}
Here assume that $\mathscr{L}=\mathscr{L}_{\mathrm{para}}$. Put $L=p^{-1}\Z_p\oplus \Z_p\oplus \Z_p\oplus \Z_p$.
Then we have $L^\vee=\Z_p\oplus \Z_p\oplus \Z_p\oplus p\Z_p$ and $[L:L^\vee]=p^2$.
In particular, $g\in G$ lies in $K_{\!\mathscr{L}}$ if and only if $gL=L$.
For every $x,y\in L$, we have $\langle x,y\rangle\in p^{-1}\Z_p$.
We introduce temporary terminology:

\begin{defn}
 An element $x$ of $L$ is said to be primitive if $x\notin pL$. A primitive element $x$ of $L$ is said to be type (I) if
 $\langle x,y\rangle\in p^{-1}\Z_p\setminus \Z_p$ for some $y\in L$, and $x$ is said to be type (II)
 if $\langle x,y\rangle\in \Z_p$ for every $y\in L$.
\end{defn}

\begin{rem}\label{rem:prim-elem-explicit}
 A primitive element of type (I) is of the form $(p^{-1}a,b,c,d)$ where $a,b,c,d\in \Z_p$ and either $a\notin p\Z_p$ or $d\notin p\Z_p$.
 A primitive element of type (II) is of the form $(a,b,c,pd)$ where $a,b,c,d\in \Z_p$ and either $b\notin p\Z_p$ or $c\notin p\Z_p$.
\end{rem}

\begin{lem}\label{lem:type-I-II}
 \begin{enumerate}
  \item Let $x\in L$ be a primitive element of type (I), and $y\in L$ an element satisfying $\langle x,y\rangle\in\Z_p$.
	Then there exists $\lambda\in \Z_p$ such that $y-\lambda x$ is not primitive of type (I).
  \item For every primitive element $x\in L$ of type (II), there exists a primitive element $y\in L$ of type (II) such that $\langle x,y\rangle=1$.
  \item Let $m\ge 1$ be an integer. For $x, x'\in L$ with $x\equiv x'\pmod{p^m}$, $x$ is primitive of type (I) (resp.\ (II))
	if and only if $x'$ is primitive of type (I) (resp.\ (II)).
  \item Let $x$ and $y$ be primitive elements of $L$. If $x$ is type (I) and $y$ is type (II), then the $\Z_p$-submodule of $L$
	generated by $x$ and $y$ is a direct summand of rank 2 of $L$.
 \end{enumerate}
\end{lem}

\begin{prf}
 i) Write $x=(p^{-1}a,b,c,d)$ and $y=(p^{-1}a',b',c',d')$ with $a,b,c,d,a',b',c',d'\in \Z_p$. Since $x$ is primitive of type (I),
 $a$ or $d$ is a unit. If $a$ is a unit, by replacing $y$ by $y-(a'/a)x$, we may assume that $a'=0$. Since $\langle x,y\rangle=p^{-1}ad'+bc'-b'c\in \Z_p$,
 $d'$ lies in $p\Z_p$. Therefore $y$ is not primitive of type (I). If $d$ is a unit, 
 by replacing $y$ by $y-(d'/d)x$, we may assume that $d'=0$. A similar argument as above tells us that $a'\in p\Z_p$, namely, $y$ is not primitive of type (I).

 ii) Since $x$ is primitive of type (II), we have $x=(a,b,c,pd)$, where $a,b,c,d\in \Z_p$ and either $b\notin p\Z_p$ or $c\notin p\Z_p$.
 We can find $b',c'\in \Z_p$ such that $bc'-cb'=1$. For $y=(0,b',c',0)\in L$,
 we have $\langle x,y\rangle=1$, as desired.

 iii) It is clear that $x$ is primitive if and only if $x'$ is primitive.
 Moreover, noting that $\langle x,y\rangle-\langle x',y\rangle=\langle x-x',y\rangle\in p^{m-1}\Z_p\subset \Z_p$, it is also immediate to see that
 $x$ is type (I) (resp.\ (II)) if and only if $x'$ is type (I) (resp.\ (II)).

 iv) Write $x=(p^{-1}a,b,c,d)$ and $y=(a',b',c',pd')$ as in Remark \ref{rem:prim-elem-explicit}. To show that they generates a direct summand of rank 2
 of $L$, it suffices to see that there exists a $2\times 2$-minor of the following matrix whose determinant is a unit in $\Z_p$:
 \[
 \begin{pmatrix}
  a& pa'\\ b&b'\\ c&c'\\ d&pd'
 \end{pmatrix}.
 \]
 It is clear from the assumptions that $a$ or $d$ is a unit, and $b'$ or $c'$ is a unit.
\end{prf}

\begin{lem}\label{lem:type-I-II-isotropic}
 Let $I$ be a totally isotropic direct summand of rank 2 of $L$.
 Then there exist a primitive element $x\in I$ of type (I) and a primitive element $y\in I$ of type (II).
 These elements $x,y$ form a basis of $I$.
\end{lem}

\begin{prf}
 Take a basis $x,y$ of $I$. Obviously they are primitive.
 First prove that one of $x$ and $y$ is type (I). Suppose that both $x$ and $y$ are type (II), and write
 $x=(a,b,c,pd)$ and $y=(a',b',c',pd')$ as in Remark \ref{rem:prim-elem-explicit}.
 Then, since $x$ and $y$ generates a direct summand of $L$, there exists a $2\times 2$-minor of the following matrix
 whose determinant is a unit in $\Z_p$:
 \[
 \begin{pmatrix}
  pa& pa'\\ b&b'\\ c&c'\\ pd&pd'
 \end{pmatrix}.
 \]
 Namely, $bc'-b'c\in\Z_p^\times$. Therefore, we have $\langle x,y\rangle=p(ad'-a'd)+bc'-b'c\in \Z_p^\times$.
 This contradicts to $\langle x,y\rangle=0$.

 Next, assume that $x$ is type (I), and find an element of $I$ which is primitive of type (II).
 Since $\langle x,y\rangle=0\in\Z_p$, Lemma \ref{lem:type-I-II} i) tells us that there exists $\lambda\in\Z_p$
 such that $y-\lambda x$ is not primitive of type (I). On the other hand, $y-\lambda x$ is primitive since $x,y-\lambda x$ form a basis of 
 a direct summand $I$. Thus we can conclude that $y-\lambda x$ is primitive of type (II).

 Finally, let $x$ (resp.\ $y$) be an arbitrary primitive element of type (I) (resp.\ (II)) of $I$. 
 Then, Lemma \ref{lem:type-I-II} iv) tells us that $x$ and $y$ form a basis of $I$.
\end{prf}

\begin{lem}\label{lem:z-w}
 Let $I$ be a totally isotropic direct summand of rank 2 of $L$, $x\in I$ a primitive element of type (I) and $y\in I$ a primitive element of type (II).
 Then there exist $z,w\in L$ satisfying the following conditions:
 \begin{itemize}
  \item $x,y,z,w$ form a basis of $L$.
  \item $\langle x,z\rangle=p^{-1}$, $\langle y,w\rangle=1$, $\langle x,w\rangle=\langle y,z\rangle=\langle z,w\rangle=0$.
 \end{itemize}
 Moreover, if we are given an element $u\in L$ satisfying $\langle x,u\rangle\in p^{-1}\Z_p\setminus\Z_p$, then we can find $z$ and $w$
 so that $\langle z,u\rangle=\langle w,u\rangle=0$ holds.
\end{lem}

\begin{prf}
 It suffices to show the latter part, since we can always find such $u\in L$.

 It is easy to see that the images $\overline{x},\overline{y},\overline{u}\in L/pL$ of $x,y,u$ are linearly independent over $\F_p$.
 Therefore, we can find $v\in L$ such that $\overline{x},\overline{y},\overline{u},\overline{v}$ form a basis of $L/pL$.
 Then $x,y,u,v$ form a basis of $L$. Since $\langle u,v\rangle\in p^{-1}\Z_p$ and $\langle u,x\rangle\in p^{-1}\Z_p\setminus \Z_p$,
 there exists $a\in\Z_p$ satisfying $\langle u,v\rangle=a\langle u,x\rangle$.
 Therefore, by replacing $v$ by $v-ax$, we may assume that $\langle u,v\rangle=0$.

 Denote by $I'$ the $\Z_p$-submodule of $L$ generated by $u$, $v$. It is a totally isotropic direct summand of rank 2 of $L$.
 The pairing $\langle\ ,\ \rangle$ induces an isomorphism of $\Q_p$-vector spaces $I'_{\Q_p}\yrightarrow{\cong}(I_{\Q_p})^\vee$,
 where $(-)_{\Q_p}=(-)\otimes_{\Z_p}\Q_p$.
 Thus we obtain a basis $x^\vee,y^\vee$ of $I'_{\Q_p}$ satisfying $\langle x,x^\vee\rangle=\langle y,y^\vee\rangle=1$
 and $\langle x,y^\vee\rangle=\langle y,x^\vee\rangle=0$. Since $I'$ is totally isotropic, $x^\vee$ and $y^\vee$ belongs to $L^\vee$. As $L^\vee\subset L$, we have $x^\vee,y^\vee\in I'_{\Q_p}\cap L=I'$.
 
 Let $I''$ be the $\Z_p$-submodule of $I'$ generated by $x^\vee$, $y^\vee$, and $L''$ the $\Z_p$-lattice of $\Q_p^4$ generated by $x,y,x^\vee,y^\vee$.
 Then we have $L/L''=I'/I''$. Since $L''$ is self-dual, we have $L^\vee\subsetneq (L'')^\vee=L''\subsetneq L$.
 Therefore $[L:L'']=p$, and thus $I'/I''=L/L''\cong \Z/p\Z$.
 In particular, there exist $a,b\in \Q_p$ such that $ax^\vee+by^\vee\in I'\setminus I''$.
 Since $y$ is primitive of type (II), $b=\langle y,ax^\vee+by^\vee\rangle\in\Z_p$. Therefore $a\notin \Z_p$, and thus $p^{-1}x^\vee\in L$.
 Now we can easily observe that $z=p^{-1}x^\vee$ and $w=y^\vee$ satisfy all conditions in the proposition;
 note that $\langle p^{-1}x^\vee,u\rangle=\langle y^\vee,u\rangle=0$ since $u$, $p^{-1}x^\vee$ and $y^\vee$ belong to $I'$, which is 
 a totally isotropic direct summand.
\end{prf}

Now we can prove Proposition \ref{prop:linear-algebra} for $h=2$.

\begin{prop}\label{prop:paramodular-h=2}
 If $\mathscr{L}=\mathscr{L}_{\mathrm{para}}$ and $h=2$, we have $K_{\!\mathscr{L},m}\backslash \mathcal{S}_{\infty,h}\cong \mathcal{S}_{\mathscr{L},m,h}$.
\end{prop}

\begin{prf}
 Let $V$ and $V'$ be elements of $\mathcal{S}_{\infty,2}$ such that $V_{\!\mathscr{L},m}=V'_{\!\mathscr{L},m}$.
 Put $I=V\cap L$ and $I'=V'\cap L$. Then $I$ and $I'$ are totally isotropic and $I\equiv I'\pmod{p^m}$. 
 It suffices to find $g\in K_{\!\mathscr{L},m}$ satisfying $gI=I'$.

 By Lemma \ref{lem:type-I-II-isotropic}, there exist a primitive element $x\in I$ of type (I) and a primitive element $y\in I$ of type (II).
 Take $x',y'\in I'$ such that $x\equiv x', y\equiv y'\pmod{p^m}$. Then $x'$ (resp.\ $y'$) is primitive of type (I) (resp.\ type (II))
 by Lemma \ref{lem:type-I-II} iii), and $x',y'$ form a basis of $I'$ by Lemma \ref{lem:type-I-II-isotropic}.
 As $p^{-m}(y'-y)\in L$, we have $\langle x,p^{-m}(y'-y)\rangle\in p^{-1}\Z_p$. Moreover, if $\langle x,p^{-m}(y'-y)\rangle\in \Z_p$,
 Lemma \ref{lem:type-I-II} i) enables us to find $\lambda\in \Z_p$ such that $p^{-m}(y'-y)-\lambda x$ is not primitive of type (I).
 Replacing $y$ by $y+p^m\lambda x$, we may assume that $p^{-m}(y'-y)$ is not primitive of type (I).

 By Lemma \ref{lem:z-w}, there exist $z,w\in L$ such that
 \begin{itemize}
  \item $x,y,z,w$ form a basis of $L$,
  \item $\langle x,z\rangle=p^{-1}$, $\langle y,w\rangle=1$, $\langle x,w\rangle=\langle y,z\rangle=\langle z,w\rangle=0$,
  \item and $\langle z,p^{-m}(y'-y)\rangle, \langle w,p^{-m}(y'-y)\rangle\in \Z_p$.
 \end{itemize}
 Indeed, if $\langle x,p^{-m}(y'-y)\rangle\notin \Z_p$, we can apply Lemma \ref{lem:z-w} to $u=p^{-m}(y'-y)$; otherwise the third condition is
 automatic since $p^{-m}(y'-y)$ is not primitive of type (I).

 Since $x',y',z,w$ form a basis of $L$, we can take $A,B,C,D\in \Q_p$
 so that $z'=Az+Bw$ and $w'=Cz+Dw$ satisfy $\langle x',z'\rangle=p^{-1}$, 
 $\langle y',w'\rangle=1$ and $\langle x',w'\rangle=\langle y',z'\rangle=0$
 (note that we have $\langle z',w'\rangle=0$ automatically).
 These are equivalent to the following identity of matrices:
 \[
  \begin{pmatrix}A&B\\C&D\end{pmatrix}\begin{pmatrix}p\langle x',z\rangle&\langle y',z\rangle\\p\langle x',w\rangle&\langle y',w\rangle\end{pmatrix}
 =\begin{pmatrix}1&0\\0&1\end{pmatrix}.
 \]
 By using the facts $x'-x\in p^mL$ and $\langle z,p^{-m}(y'-y)\rangle, \langle w,p^{-m}(y'-y)\rangle\in \Z_p$, it is immediate to observe that
  \[
  \begin{pmatrix}p\langle x',z\rangle&\langle y',z\rangle\\p\langle x',w\rangle&\langle y',w\rangle\end{pmatrix}\in 1+p^m M_2(\Z_p).
  \]
 Therefore we can conclude that $A,D\in 1+p^m\Z_p$ and $B,C\in p^m\Z_p$. 
 In other words, we have $z',w'\in L$ and $z'\equiv z,w'\equiv w\pmod{p^m}$.
 In particular, $x',y',z',w'$ form a basis of $L$.
 Let $g$ be the automorphism of $L$ that maps $x$ to $x'$, $y$ to $y'$, $z$ to $z'$ and $w$ to $w'$.
 It is clear that $g$ is an element of $K_{\!\mathscr{L},m}$ satisfying $gI=I'$. This completes the proof.
\end{prf}

Finally we prove Proposition \ref{prop:linear-algebra} for $h=1$.

\begin{prop}
 If $\mathscr{L}=\mathscr{L}_{\mathrm{para}}$ and $h=1$, we have $K_{\!\mathscr{L},m}\backslash \mathcal{S}_{\infty,h}\cong \mathcal{S}_{\mathscr{L},m,h}$.
\end{prop}

\begin{prf}
 Let $V$ and $V'$ be elements of $\mathcal{S}_{\infty,1}$ such that $V_{\!\mathscr{L},m}=V'_{\!\mathscr{L},m}$.
 Put $I=V\cap L$ and $I'=V'\cap L$. It suffices to find $g\in K_{\!\mathscr{L},m}$ satisfying $gI=I'$.

 Take a basis $x$ of $I$ and a basis $x'$ of $I'$ satisfying $x\equiv x'\pmod{p^m}$.
 First consider the case where $x$ is primitive of type (I). Then $x'$ is also primitive of type (I) by Lemma \ref{lem:type-I-II} iii).
 Take $z,z'\in L$ such that $\langle x,z\rangle=\langle x',z'\rangle=p^{-1}$.
 Take an arbitrary primitive element $y\in L$ of type (II). By replacing it by $y-p\langle x,y\rangle z$, we may assume that $\langle x,y\rangle=0$.
 Put $y'=y-p\langle x',y\rangle z'$. Then $y'$ is a primitive element of type (II) satisfying $\langle x',y'\rangle=0$ and $y\equiv y'\pmod{p^m}$.
 By the proof of Proposition \ref{prop:paramodular-h=2}, we can find $g\in K_{\!\mathscr{L},m}$ such that $gx=x'$.

 Next consider the case where $x$ is primitive of type (II). In this case $x'$ is also primitive of type (II).
 By Lemma \ref{lem:type-I-II} ii), there exist primitive elements $z,z'\in L$ of type (II) such that $\langle x,z\rangle=\langle x',z'\rangle=1$.
 Let us prove that we can find a primitive element $y\in L$ of type (I) satisfying $\langle x,y\rangle=0$ and $\langle p^{-m}(x'-x),y\rangle\in\Z_p$.
 If $p^{-m}(x'-x)$ is primitive of type (I), put $y=p^{-m}(x'-x)-az$ where $a=\langle x,p^{-m}(x'-x)\rangle$.
 Then $a\in\Z_p$ as $x$ is primitive of type (II), and $\langle p^{-m}(x'-x),y\rangle=a\langle z,y\rangle\in\Z_p$ since $z$ is
 primitive of type (II). If $p^{-m}(x'-x)$ is not primitive of type (I), take an arbitrary primitive element $w\in L$ of type (I)
 and put $y=w-\langle x,w\rangle z$. In this case, the condition $\langle p^{-m}(x'-x),y\rangle\in\Z_p$ is automatic.
 Put $y'=y-\langle x',y\rangle z'$. Then we have $\langle x',y'\rangle=0$.
 Furthermore we have $y\equiv y'\pmod{p^m}$, for $\langle x',y\rangle=p^m\langle p^{-m}(x'-x),y\rangle\in p^m\Z_p$.
 Hence, again by the proof of Proposition \ref{prop:paramodular-h=2}, we can find $g\in K_{\!\mathscr{L},m}$ such that $gx=x'$.
\end{prf}

Now a proof of Proposition \ref{prop:linear-algebra} for the paramodular case is complete.

\subsubsection{Siegel parahoric case}
Here we assume that $\mathscr{L}=\mathscr{L}_{\mathrm{Siegel}}$.
Put $L_0=\Z_p^4$ and $L_1=\Z_p\oplus \Z_p\oplus p\Z_p\oplus p\Z_p$.
We have $pL_0\subset L_1\subset L_0$ and $L_0/L_1\cong L_1/pL_0\cong \F_p^2$.
An element $g\in G$ lies in $K_{\!\mathscr{L}}$ if and only if $gL_0=L_0$ and $gL_1=L_1$.

\begin{lem}\label{lem:Nakayama}
 Let $x,y\in L_0$ and $z,w\in L_1$ be elements satisfying
 \[
 \langle x,y\rangle=0,\quad \langle x,z\rangle=\langle y,w\rangle=1,\quad \langle x,w\rangle=\langle y,z\rangle=0,
 \quad \langle z,w\rangle=0.
 \]
 Then, $x,y,z,w$ (resp.\ $px,py,z,w$) form a basis of $L_0$ (resp.\ $L_1$).
\end{lem}

\begin{prf}
 Let us prove that the images $\overline{x},\overline{y}$ of $x,y$ in $L_0/L_1$
 form an $\F_p$-basis of $L_0/L_1$. Assume that $u=ax+by\in L_1$ for $a,b\in\Z_p$.
 Then, $a=\langle u,z\rangle\in p\Z_p$ and $b=\langle u,w\rangle\in p\Z_p$, 
 since the pairing $\langle\ ,\ \rangle$ on $L_1$ takes its value in $p\Z_p$.
 This implies that $\overline{x},\overline{y}\in L_0/L_1$ are linearly independent over $\F_p$,
 and thus they form a basis.
 Similarly, we can prove that $\overline{z},\overline{w}\in L_1/pL_0$ form a basis of $L_1/pL_0$.
 From this, it is immediate to conclude the lemma.
\end{prf}

\begin{prop}\label{prop:Siegel-parahoric-basis}
 Let $I_0$ be a totally isotropic direct summand of rank $2$ of $L_0$ and put $I_1=I_0\cap L_1$.
 \begin{enumerate}
  \item Assume that $I_0=I_1$. For a basis $x,y\in I_0=I_1$, there exist $z,w\in L_0$ such that
	\[
	\langle x,z\rangle=\langle y,w\rangle=1,\quad \langle x,w\rangle=\langle y,z\rangle=0,
	\quad \langle z,w\rangle=0.
	\]
  \item Assume that $[I_0:I_1]=p$. For $x\in I_0\setminus I_1$ and $y\in I_1\setminus pI_0$,
	$x,y$ form a basis of $I_0$. Moreover, there exist $z\in L_1$ and $w\in L_0$ such that
	\[
	\langle x,z\rangle=\langle y,w\rangle=1,\quad \langle x,w\rangle=\langle y,z\rangle=0,\quad \langle z,w\rangle=0.
	\]
  \item Assume that $[I_0:I_1]=p^2$. For a basis $x,y\in I_0$,
	there exist $z,w\in L_1$ such that
	\[
	\langle x,z\rangle=\langle y,w\rangle=1,\quad \langle x,w\rangle=\langle y,z\rangle=0,\quad \langle z,w\rangle=0.
	\]
 \end{enumerate}
\end{prop}

\begin{prf}
 i) Take $u,v\in L_0$ such that $x,y,u,v$ form a basis of $L_0$.
 Since $y\notin pL_0$ and $\langle x,y\rangle=0$, either $\langle u,y\rangle$ or $\langle v,y\rangle$ 
 is a unit in $\Z_p$.
 Therefore we may assume that $\langle u,y\rangle\in\Z_p^\times$.
 Then there exists $a\in\Z_p$ such that $\langle u,v-ay\rangle=0$.
 Replacing $v$ by $v-ay$, we may assume that $\langle u,v\rangle=0$.
 Let $J$ be a direct summand of $L_0$ generated by $u,v$, which is totally isotropic.
 The pairing $\langle\ ,\ \rangle$ induces an isomorphism $J\yrightarrow{\cong} I_0^\vee$.
 Therefore, there exists a basis $z,w$ of $J$ satisfying
 \[
 \langle x,z\rangle=\langle y,w\rangle=1,\quad \langle x,w\rangle=\langle y,z\rangle=0.
 \]
 This concludes i).

 ii) It is immediate to see that $x,y$ form a basis of $L_0$. We shall prove the existence of $z,w$.
 Take $u\in L_0$ (resp.\ $v\in L_1$) so that $\overline{x},\overline{u}$ (resp.\ $\overline{y},\overline{v}$)
 form a basis of $L_0/L_1$ (resp.\ $L_1/pL_0$).
 Since $y\notin pL_0$ and $\langle x,y\rangle=0$, either $\langle u,y\rangle$ or $\langle v,y\rangle$ is a
 unit in $\Z_p$. As $\langle v,y\rangle\in p\Z_p$, we have $\langle u,y\rangle\in\Z_p^\times$.
 Replacing $v$ by $v-ay$ for a suitable $a\in\Z_p$, we may assume that $\langle u,v\rangle=0$.
 By the same argument as above, we can find $z,w\in \Z_pu+\Z_pv$ satisfying
 \[
  \langle x,z\rangle=\langle y,w\rangle=1,\quad \langle x,w\rangle=\langle y,z\rangle=0.
 \]
 We should prove that $z\in L_1$. 
 Since $px,y,pu,v$ form a basis of $L_1$,
 it is easy to see that $\langle u,z\rangle\in p\Z_p$ for every $u\in L_1$.
 Therefore $z\in pL_1^\vee=L_1$, as desired. 

 iii) Take $u,v\in L_1$ such that $\overline{u},\overline{v}\in L_1/pL_0$
 form a basis of $L_1/pL_0$.
 Since $L_0/L_1\cong I_0/I_1$, $x,y,u,v$ form a basis of $L_0$.
 As in i) and ii), we may assume that $\langle u,y\rangle\in\Z_p^\times$.
 For $a=\langle u,y\rangle^{-1}\langle u,v\rangle\in \Z_p$,
 we have $\langle u,v-ay\rangle=0$. Since $\langle u,v\rangle\in p\Z_p$, the element $ay$ belongs to $L_1$.
 Therefore, replacing $v$ by $v-ay$, we may assume that $\langle u,v\rangle=0$.
 By the same argument as above, we can find $z,w\in \Z_pu+\Z_pv\subset L_1$ satisfying
 \[
  \langle x,z\rangle=\langle y,w\rangle=1,\quad \langle x,w\rangle=\langle y,z\rangle=0.
 \]
\end{prf}

\begin{prop}\label{prop:Siegel-parahoric-h=2}
 If $\mathscr{L}=\mathscr{L}_{\mathrm{Siegel}}$ and $h=2$, 
 we have $K_{\!\mathscr{L},m}\backslash \mathcal{S}_{\infty,h}\cong \mathcal{S}_{\mathscr{L},m,h}$.
\end{prop}

\begin{prf}
 Let $V$ and $V'$ be elements of $\mathcal{S}_{\infty,2}$ such that $V_{\!\mathscr{L},m}=V'_{\!\mathscr{L},m}$.
 Put $I_i=V\cap L_i$ and $I'_i=V'\cap L_i$ for $i=0,1$. Then $I_i$ and $I'_i$ are totally isotropic
 and $I_i\equiv I'_i\pmod{p^m}$. 
 It suffices to find $g\in K_{\!\mathscr{L},m}$ satisfying $gI_0=I'_0$.

 Let us take $x,y,z,w\in L_0$ as in Proposition \ref{prop:Siegel-parahoric-basis}.
 For $x',y'\in I'_0$ with $x'-x,y'-y\in p^mL_0$, we can take $A,B,C,D\in \Q_p$
 so that $z'=Az+Bw$ and $w'=Cz+Dw$ satisfy
 $\langle x',z'\rangle=\langle y',w'\rangle=1$ and $\langle x',w'\rangle=\langle y',z'\rangle=0$.
 These are equivalent to the following identity of matrices:
 \[
 \begin{pmatrix}A&B\\C&D\end{pmatrix}\begin{pmatrix}\langle x',z\rangle&\langle y',z\rangle\\\langle x',w\rangle&\langle y',w\rangle\end{pmatrix}
  =\begin{pmatrix}1&0\\0&1\end{pmatrix}.
 \]
 By using the assumption $x'-x, y'-y\in p^mL_0$, it is immediate to observe that
  \[
  \begin{pmatrix}\langle x',z\rangle&\langle y',z\rangle\\\langle x',w\rangle&\langle y',w\rangle\end{pmatrix}\in 1+p^m M_2(\Z_p).
  \]
 Therefore we can conclude that $A,D\in 1+p^m\Z_p$ and $B,C\in p^m\Z_p$.
 In other words, $z',w'$ lie in $\Z_pz+\Z_pw$ and $z'-z,w'-w\in p^m\Z_pz+p^m\Z_pw$.

 Now consider the three cases in Proposition \ref{prop:Siegel-parahoric-basis} separately.
 First assume that $I_0=I_1$. Then $I_0'=I_1'$, as $I_0'\subset I_0+p^mL_0\subset L_1+pL_0=L_1$.
 By the assumption $I_1\equiv I_1'\pmod{p^m}$, we may take $x',y'\in I_0'$ above so that $x-x',y-y'\in p^mL_1$.
 Lemma \ref{lem:Nakayama} tells us that
 \begin{itemize}
  \item $x,y,z,w$ and $x',y',z',w'$ are bases of $L_0$,
  \item $x,y,pz,pw$ and $x',y',pz',pw'$ are bases of $L_1$, and
  \item $x',y'$ form a basis of $I_0'$.
 \end{itemize}
 Let $g$ be the automorphism of $L_0$ that maps $x$ to $x'$, $y$ to $y'$, $z$ to $z'$ and $w$ to $w'$.
 It is clear that $g$ is an element of $K_{\!\mathscr{L},m}$ satisfying $gI_0=I_0'$.

 Next consider the case where $[I_0:I_1]=p$. In this case, we can choose $y'$ so that $y'$ belongs to $I'_1$
 and $y'-y\in p^mL_1$, as $I_1\equiv I_1'\pmod{p^m}$.
 Then $z'\in L_1$ and $z'-z\in p^mL_1$; indeed, we have $\langle y',z\rangle=p^m\langle p^{-m}(y'-y),z\rangle\in p^{m+1}\Z_p$ (note that $z\in L_1$) and thus $B\in p^{m+1}\Z_p$.
 Therefore, Lemma \ref{lem:Nakayama} tells us that
 \begin{itemize}
  \item $x,y,z,w$ and $x',y',z',w'$ are bases of $L_0$,
  \item $px,y,z,pw$ and $px',y',z',pw'$ are bases of $L_1$, and
  \item $x',y'$ form a basis of $I_0'$.
 \end{itemize}
 Hence the automorphism of $L_0$ defined as above gives an element $g\in K_{\!\mathscr{L},m}$ satisfying $gI_0=I_0'$.

 Finally consider the case where $[I_0:I_1]=p^2$. Since $z,w\in L_1$, we have
 $z',w'\in L_1$ and $z'-z,w'-w\in p^mL_1$. Lemma \ref{lem:Nakayama} tells us that
 \begin{itemize}
  \item $x,y,z,w$ and $x',y',z',w'$ are bases of $L_0$,
  \item $px,py,z,w$ and $px',py',z',w'$ are bases of $L_1$, and
  \item $x',y'$ form a basis of $I_0'$.
 \end{itemize}
 Hence the automorphism of $L_0$ defined as above gives an element $g\in K_{\!\mathscr{L},m}$ satisfying $gI_0=I_0'$.
\end{prf}

\begin{prop}
 If $\mathscr{L}=\mathscr{L}_{\mathrm{Siegel}}$ and $h=1$, we have $K_{\!\mathscr{L},m}\backslash \mathcal{S}_{\infty,h}\cong \mathcal{S}_{\mathscr{L},m,h}$.
\end{prop}

\begin{prf}
 Let $V$ and $V'$ be elements of $\mathcal{S}_{\infty,1}$ such that $V_{\!\mathscr{L},m}=V'_{\!\mathscr{L},m}$.
 Put $I_i=V\cap L_i$ and $I'_i=V'\cap L_i$ for $i=0,1$. Then we have $I_i\equiv I'_i\pmod{p^m}$. 
 It suffices to find $g\in K_{\!\mathscr{L},m}$ satisfying $gI_0=I'_0$.

 Take a basis $x$ of $I_0$ and an element $x'\in I_0'$ such that
 $x'-x\in p^mL_0$. Then $x'$ is also a basis of $I_0'$, for $x'\notin pL_0$.
 Moreover, since $x'-x\in pL_0\subset L_1$, $x\in L_1$ if and only if $x'\in L_1$.

 First consider the case where $I_0=I_1$. In this case, $I_0'=I_1'$ as $x'\in I_1'$.
 By the assumption $I_1\equiv I_1'\pmod{p^m}$, we may take $x'$ above so that $x-x'\in p^mL_1$.
 Let us observe that there exists $y\in L_0\setminus L_1$ such that $\langle x,y\rangle=0$.
 Consider an $\F_p$-linear map $\langle x,-\rangle\colon L_0/L_1\longrightarrow \F_p$.
 Since $\dim_{\F_p}L_0/L_1=2$, the kernel of this map is non-trivial.
 In other words, there exists $u\in L_0\setminus L_1$ such that $\langle x,u\rangle\in p\Z_p$.
 Take $v\in L_0$ such that $\langle x,v\rangle=1$ and put $y=u-\langle x,u\rangle v$.
 Then we have $\langle x,y\rangle=0$ and $y\in L_0\setminus L_1$
 (note that $\langle x,u\rangle v\in pL_0\subset L_1$).
 Take $z\in L_0$ satisfying $\langle x',z\rangle=1$ and put $y'=y-\langle x',y\rangle z$.
 Since $\langle x',y\rangle=\langle x'-x,y\rangle\in p^m\Z_p$,
 we have $y'\in L_0\setminus L_1$ and $y'-y\in p^mL_0$.
 Let $\widetilde{V}$ (resp.\ $\widetilde{V}'$) be a subspace of $\Q_p^4$ generated by $x,y$ (resp.\ $x',y'$),
 and put $\widetilde{I}_i=\widetilde{V}\cap L_i$, $\widetilde{I}'_i=\widetilde{V}'\cap L_i$ for $i=1,2$.
 Then it is easy to see the following:
 \begin{itemize}
  \item $\widetilde{V}, \widetilde{V}'\in\mathcal{S}_{\infty,2}$.
  \item $x,y$ (resp.\ $x,py$) form a basis of $\widetilde{I}_0$ (resp.\ $\widetilde{I}_1$).
  \item $x',y'$ (resp.\ $x',py'$) form a basis of $\widetilde{I}'_0$ (resp.\ $\widetilde{I}'_1$).
  \item $\widetilde{I}_i\equiv \widetilde{I}'_i \pmod{p^m}$.
 \end{itemize}
 Therefore, by the proof of Proposition \ref{prop:Siegel-parahoric-h=2}, we can find
 $g\in K_{\!\mathscr{L},m}$ such that $gx=x'$.

 Next consider the case where $I_0\neq I_1$. In this case, $I_0'\neq I_1'$ since $x'\notin I_1'$.
 In the same way as above, we can find $y\in L_1\setminus pL_0$ such that $\langle x,y\rangle=0$.
 Since $x'\notin L_1=pL_1^\vee$, there exists $z\in L_1$ such that $\langle x',z\rangle=1$.
 Put $y'=y-\langle x',y\rangle z$. As $\langle x',y\rangle=\langle x'-x,y\rangle\in p^m\Z_p$,
 we have $y'\in L_1\setminus pL_0$ and $y'-y\in p^mL_1$.
 Define $\widetilde{V}$, $\widetilde{V}'$, $\widetilde{I}_i$ and $\widetilde{I}'_i$ as above.
 Then it is easy to see the following:
 \begin{itemize}
  \item $\widetilde{V}, \widetilde{V}'\in\mathcal{S}_{\infty,2}$.
  \item $x,y$ (resp.\ $px,y$) form a basis of $\widetilde{I}_0$ (resp.\ $\widetilde{I}_1$).
  \item $x',y'$ (resp.\ $px',y'$) form a basis of $\widetilde{I}'_0$ (resp.\ $\widetilde{I}'_1$).
  \item $\widetilde{I}_i\equiv \widetilde{I}'_i \pmod{p^m}$.
 \end{itemize}
 Therefore, by the proof of Proposition \ref{prop:Siegel-parahoric-h=2}, we can find
 $g\in K_{\!\mathscr{L},m}$ such that $gx=x'$.
\end{prf}

Now a proof of Proposition \ref{prop:linear-algebra} for the Siegel parahoric case is complete.

\section{Open covering of the Rapoport-Zink space}
From this section, we convert the right action of $G\times J$ on the Rapoport-Zink tower
to the left action by taking inverse.
Therefore, the action of $(g,h)\in G\times J$ on the cohomology $H^i_{\mathrm{RZ}}$
is given by $(g^{-1},h^{-1})^*$.

\subsection{Definition of open covering}
Fix a chain of lattices $\mathscr{L}$ of $\Q_p^{2d}$.
Then we can consider $\M^\flat_{\mathscr{L}}$ (\cf Definition \ref{defn:formal-model})
and its rigid generic fiber $M_{\mathscr{L}}$.
We write $\Irr(\mathscr{L})$ for the set of irreducible components of $\M^{\flat,\red}_{\mathscr{L}}$.
For each $\alpha\in \Irr(\mathscr{L})$, the set $\M^{\flat,\red}_{\mathscr{L}}\setminus \bigcup_{\beta\in \Irr(\mathscr{L}),\alpha\cap \beta=\varnothing}\beta$ is open in $\M^{\flat,\red}_{\mathscr{L}}$, as $\M^{\flat,\red}_{\mathscr{L}}$ is
locally of finite type over $\overline{\F}_p$.
Let $\mathscr{U}_\alpha$ be the open formal subscheme of
$\M^{\flat}_{\mathscr{L}}$ satisfying 
$\mathscr{U}^\red_\alpha=\M^{\flat,\red}_{\mathscr{L}}\setminus \bigcup_{\beta\in \Irr(\mathscr{L}),\alpha\cap \beta=\varnothing}\beta$, and $U_\alpha$ the rigid generic fiber of $\mathscr{U}_\alpha$.
As $\alpha\subset \mathscr{U}^\red_\alpha$, $\{\mathscr{U}_\alpha\}_{\alpha\in\Irr(\mathscr{L})}$
(resp.\ $\{U_\alpha\}_{\alpha\in\Irr(\mathscr{L})}$) gives an open covering of $\M^\flat_{\mathscr{L}}$
(resp.\ $M_{\mathscr{L}}$).
For $\alpha,\beta\in \Irr(\mathscr{L})$, note that $U_\alpha\cap U_\beta\neq\varnothing$ if and only if
$\mathscr{U}^{\mathrm{red}}_\alpha\cap \mathscr{U}^{\mathrm{red}}_\beta\neq\varnothing$,
since $U_\alpha\cap U_\beta=(\mathscr{U}_\alpha\cap \mathscr{U}_\beta)^{\rig}$ and
$\mathscr{U}_\alpha\cap \mathscr{U}_\beta$ is flat over $\Z_{p^\infty}$.

Clearly $J$ naturally acts on $\Irr(\mathscr{L})$, and we have $h \mathscr{U}_\alpha=\mathscr{U}_{h\alpha}$
for $h\in J$. 

\begin{lem}\label{lem:open-cov-properties}
 \begin{enumerate}
  \item Put $\Irr(\mathscr{L})_{\alpha,0}=\{\alpha\}$ and for $m\ge 1$ define 
	the set $\Irr(\mathscr{L})_{\alpha,m}$ inductively as follows;
	$\Irr(\mathscr{L})_{\alpha,m}$ consists of $\beta\in \Irr(\mathscr{L})$ which intersects
	some element in $\Irr(\mathscr{L})_{\alpha,m-1}$. Then $\Irr(\mathscr{L})_{\alpha,m}$ is
	a finite set.
  \item The action of $J$ on $\Irr(\mathscr{L})$ has finite orbits.
  \item For every $\alpha,\beta\in \Irr(\mathscr{L})$,
	the subset $\{h\in J\mid h\alpha\cap \beta\neq \varnothing\}$ of $J$
	is contained in a compact subset of $J$.
  \item For each $\alpha\in\Irr(\mathscr{L})$,
	the subgroup $J_\alpha=\{h\in J\mid h\alpha=\alpha\}$ is open and compact.
 \end{enumerate}
\end{lem}

\begin{prf}
 Since $\M^{\flat,\red}_{\mathscr{L}}$ is locally of finite type
 and every irreducible component of $\M^{\flat,\red}_{\mathscr{L}}$ is projective, 
 $\Irr(\mathscr{L})_{\alpha,1}$ is a finite set. If $\Irr(\mathscr{L})_{\alpha,m-1}$ is a finite set,
 $\Irr(\mathscr{L})_{\alpha,m}=\bigcup_{\beta\in \Irr(\mathscr{L})_{\alpha,m-1}}\Irr(\mathscr{L})_{\beta,1}$
 is also a finite set. This concludes the proof of i).

 We assume that $\mathscr{L}=\mathscr{L}_0=\{p^m\Z_p^{2d}\mid m\in\Z\}$ and prove ii) and iii). 
 In this case, $\M^\flat_{\mathscr{L}_0}=\M$ and the action of $J$ on $\Irr(\mathscr{L}_0)$ is
 transitive (\cite[Theorem 2]{MR2466182}).
 Thus ii) is immediate. 
 Let us prove iii). By the transitivity, we may assume that $\alpha=\beta$.
 Let $\breve{\mathscr{M}}_{\GL}$ be the Rapoport-Zink space for $\GL(2d)$
 corresponding to the $p$-divisible group $\X$ and $J_{\GL}$ the group of self-quasi-isogenies of $\X$
 which naturally acts on $\breve{\mathscr{M}}_{\GL}$. The group $J$ is a closed subgroup of $J_{\GL}$.
 By forgetting the condition on polarization, we have a natural closed immersion $\M\hooklongrightarrow \M_{\GL}$
 and the actions of $J$ and $J_{\GL}$ are compatible. Therefore, it suffices to show that
 $\{h\in J_{\GL_n}\mid h\alpha\cap \alpha\neq\varnothing\}$
 is contained in a compact subset of $J_{\GL}$,
 where $\alpha$ is regarded as a quasi-compact closed subset of $\M_{\GL}^{\red}$. 
 This is essentially proved in \cite[proof of Proposition 2.34]{MR1393439}.

 Next we assume that 
 $\mathscr{L}=\mathscr{L}_{\mathrm{Iw}}=\{p^m\Z_p^i\oplus p^{m+1}\Z_p^{2d-i}\mid m\in\Z, 0\le i\le 2d\}$.
 Then we have a natural morphism
 $\pi\colon \M^\flat_{\mathscr{L}_{\mathrm{Iw}}}\longrightarrow \M^\flat_{\mathscr{L}_0}$, which is proper.
 Fix an arbitrary element $\beta\in\Irr(\mathscr{L}_0)$ and let $\alpha_1,\ldots,\alpha_m$ be
 the collection of elements of $\Irr(\mathscr{L}_{\mathrm{Iw}})$ which are contained in $\pi^{-1}(\beta)$.
 Since $\pi$ is $J$-equivariant, it is easy to observe 
 $\Irr(\mathscr{L}_\mathrm{Iw})=\bigcup_{h\in J}\bigcup_{i=1}^mh\alpha_i$. 
 This concludes the proof of ii) in this case. To prove iii), 
 let $\alpha_0$ (resp.\ $\beta_0$) be an element of
 $\Irr(\mathscr{L}_0)$ which contains $\pi(\alpha)$ (resp.\ $\pi(\beta)$).
 Then, $\{h\in J\mid h\alpha\cap \beta\neq \varnothing\}$ is contained in 
 $\{h\in J\mid h\alpha_0\cap \beta_0\neq \varnothing\}$, and thus we are reduced to the case 
 $\mathscr{L}=\mathscr{L}_0$.

 Finally we consider an arbitrary $\mathscr{L}$. Changing $\mathscr{L}$ by $g\mathscr{L}$ with $g\in G$,
 we may assume that $\mathscr{L}$ is contained in $\mathscr{L}_{\mathrm{Iw}}$.
 Then we have a natural morphism
 $\pi\colon \M^\flat_{\mathscr{L}_{\mathrm{Iw}}}\longrightarrow \M^\flat_{\mathscr{L}}$, which is proper.
 By Remark \ref{rem:m=0} and Proposition \ref{prop:formal-model-gen-fiber},
 the induced morphism $\pi\colon M_{\mathscr{L}_{\mathrm{Iw}}}\longrightarrow M_{\mathscr{L}}$ is
 surjective, and therefore the morphism 
 $\M^{\flat,\red}_{\mathscr{L}_{\mathrm{Iw}}}\longrightarrow \M^{\flat,\red}_{\mathscr{L}}$ is also surjective
 (note that $\M^\flat_{\mathscr{L}}$ is flat over $\Z_{p^\infty}$).
 For ii), let $\beta_1,\ldots,\beta_m\in\Irr(\mathscr{L}_{\mathrm{Iw}})$ be
 elements such that $\Irr(\mathscr{L}_{\mathrm{Iw}})=\bigcup_{h\in J}\bigcup_{i=1}^mh\beta_i$,
 and $\alpha_i$ an element of $\Irr(\mathscr{L})$ which contains $\pi(\beta_i)$.
 Then $\{h\alpha_i\}_{h\in J,1\le i\le m}$ cover $\M^{\flat,\red}_{\mathscr{L}}$,
 and we have $\Irr(\mathscr{L})=\bigcup_{h\in J}\bigcup_{i=1}^mh\alpha_i$.
 To prove iii), take $\alpha'_1,\ldots,\alpha'_k,\beta'_1,\ldots,\beta'_l\in\Irr(\mathscr{L}_{\mathrm{Iw}})$
 so that
 $\pi^{-1}(\alpha)\subset \bigcup_{i=1}^{k}\alpha'_i$ and $\pi^{-1}(\beta)\subset \bigcup_{j=1}^{l}\beta'_j$.
 We have
 \begin{align*}
 &\{h\in J\mid h\alpha\cap \beta\neq\varnothing\}=\{h\in J\mid \pi^{-1}(h\alpha)\cap \pi^{-1}(\beta)\neq\varnothing\}\\
  &\qquad \subset \bigcup_{i=1}^k\bigcup_{j=1}^l\{h\in J\mid h\alpha'_i\cap \beta'_j\neq \varnothing\}.
 \end{align*}
 We already know that the latter set is contained in a compact subset of $J$. 
 This completes the proof of ii) and iii).
 
 Let us prove iv). By \cite[Proposition 2.3.11]{MR2074714}, there exists an open subgroup of $J$
 whose elements stabilize $\alpha$. This implies that $J_{\alpha}$ is an open subgroup of $J$.
 In particular, $J_{\alpha}$ is a closed subgroup of $J$. 
 On the other hand, by iii), $J_\alpha$ is contained in a compact subset of $J$.
 Thus $J_\alpha$ is compact, as desired.
\end{prf}

\begin{cor}\label{cor:open-cov-properties-1}
 For each $\alpha\in \Irr(\mathscr{L})$, there exist only finitely many $\beta\in \Irr(\mathscr{L})$
 with $U_\alpha\cap U_\beta\neq\varnothing$.
\end{cor}

\begin{prf}
 By the construction, it is easy to see that 
 $\mathscr{U}^\red_\alpha\cap \mathscr{U}^\red_\beta=\varnothing$ unless $\beta\in \Irr(\mathscr{L})_{\alpha,2}$.
 Thus the claim follows from Lemma \ref{lem:open-cov-properties} i).
\end{prf}

Consider the quotient $\M^\flat_{\mathscr{L}}/p^\Z$ of $\M^\flat_{\mathscr{L}}$ by the discrete subgroup $p^\Z\subset J$.
For $\overline{\alpha}\in\Irr(\mathscr{L})/p^\Z$, we write $\mathscr{U}_{\overline{\alpha}}$
(resp.\ $U_{\overline{\alpha}}$) for the image of $\mathscr{U}_{\alpha}$ (resp.\ $U_{\alpha}$).
Recall that we have an action of $N_\mathscr{L}$ on $\M^\flat_{\mathscr{L}}/p^\Z$. Therefore,
$N_\mathscr{L}$ also acts on $\Irr(\mathscr{L})/p^\Z$.
This action factors through the quotient $N_\mathscr{L}/p^\Z K_{\!\mathscr{L}}$, which is a finite group.
Therefore, each $N_{\mathscr{L}}$-orbit in $\Irr(\mathscr{L})/p^\Z$ is a finite set.
We denote by $I(\mathscr{L})$ the set of $N_{\mathscr{L}}$-orbits in $\Irr(\mathscr{L})/p^\Z$.
For $\lambda\in I(\mathscr{L})$,
we put $\mathscr{U}_\lambda=\bigcup_{\overline{\alpha}\in\lambda}\mathscr{U}_{\overline{\alpha}}$
and $U_\lambda=\bigcup_{\overline{\alpha}\in\lambda}U_{\overline{\alpha}}$. Then
$\mathscr{U}_\lambda$ and $U_\lambda$ are stable under the action of $N_{\mathscr{L}}$.
Clearly $J/p^\Z$ naturally acts on $I(\mathscr{L})$, and we have $h\mathscr{U}_\lambda=\mathscr{U}_{h\lambda}$
and $hU_\lambda=U_{h\lambda}$ for $h\in J/p^\Z$ and $\lambda\in I(\mathscr{L})$. 

Since $U_\lambda$ is a finite union of $U_{\overline{\alpha}}$'s, the following corollary is
an immediate consequence of Lemma \ref{lem:open-cov-properties} and Corollary \ref{cor:open-cov-properties-1}:

\begin{cor}\label{cor:open-cov-properties-2}
 \begin{enumerate}
  \item For each $\lambda\in I(\mathscr{L})$, there exist only finitely many $\mu\in I(\mathscr{L})$
	with $U_\lambda\cap U_\mu\neq\varnothing$.
  \item The action of $J$ on $I(\mathscr{L})$ has finite orbits. In particular,
	there exist finitely many elements $\lambda_1,\ldots,\lambda_m$ of $I(\mathscr{L})$
	such that $M_{\mathscr{L}}/p^\Z=\bigcup_{h\in J}\bigcup_{i=1}^mhU_{\lambda_i}$.
  \item For every $\lambda,\mu\in I(\mathscr{L})$,
	the subset $\{h\in J/p^\Z\mid hU_\lambda\cap U_\mu\neq \varnothing\}$ of $J/p^\Z$
	is contained in a compact subset of $J/p^\Z$.
  \item For each $\lambda\in I(\mathscr{L})$,
	the subgroup $J_\lambda=\{h\in J\mid h\lambda=\lambda\}$ is open and compact-mod-center.
 \end{enumerate}
\end{cor}

\begin{cor}\label{cor:iii-compactness}
 For $\lambda,\mu\in I(\mathscr{L})$, the set $\{h\in J/p^\Z\mid hU_\lambda\cap U_\mu\neq \varnothing\}$ is
 in fact a compact open subset of $J/p^\Z$.
\end{cor}

\begin{prf}
 Put $C_{\lambda,\mu}=\{h\in J/p^\Z\mid hU_\lambda\cap U_\mu\neq \varnothing\}$ and consider the closure
 $\overline{C_{\lambda,\mu}}$ in $J/p^\Z$. By Corollary \ref{cor:open-cov-properties-2} iii), it is compact.
 The group $J_\lambda$ acts on $C_{\lambda,\mu}$ and $\overline{C_{\lambda,\mu}}$ on the right.
 Since $J_\lambda/p^\Z$ is an open subgroup of $J/p^\Z$, the quotient $\overline{C_{\lambda,\mu}}/J_\lambda$
 is finite. Therefore $C_{\lambda,\mu}/J_\lambda$ is also a finite set.
 Since $J_\lambda/p^\Z$ is compact open, so is $C_{\lambda,\mu}$.
\end{prf}

 For a finite subset $\underline{\lambda}=\{\lambda_1,\ldots,\lambda_m\}$ of $I(\mathscr{L})$,
 we put $\mathscr{U}_{\underline{\lambda}}=\bigcap_{i=1}^m \mathscr{U}_{\lambda_i}$, 
 $U_{\underline{\lambda}}=\bigcap_{i=1}^m U_{\lambda_i}$ and
 $J_{\underline{\lambda}}=\{h\in J\mid h\underline{\lambda}=\underline{\lambda}\}$.
 Since the group $J_{\underline{\lambda}}/\bigcap_{i=1}^mJ_{\lambda_i}$ acts faithfully on $\underline{\lambda}$,
 it is a finite group. Thus, Corollary \ref{cor:open-cov-properties-2} iv) tells us that
 $J_{\underline{\lambda}}$ is an open compact-mod-center subgroup of $J$.

 For an integer $s\ge 1$, let $I(\mathscr{L})_s$ be the set of subsets 
 $\underline{\lambda}\subset I(\mathscr{L})$ such that $\#\underline{\lambda}=s$ and 
 $U_{\underline{\lambda}}\neq\varnothing$. The group $J$ (or $J/p^\Z$) naturally acts on $I(\mathscr{L})_s$.

\begin{lem}\label{lem:finite-orbits}
 For an integer $s\ge 1$, the action of $J$ on $I(\mathscr{L})_s$ has finite orbits.
\end{lem}

\begin{prf}
 Let $I(\mathscr{L})^\sim_s$ be the subset of $I(\mathscr{L})^s$ consisting of elements
 $(\lambda_1,\ldots,\lambda_s)$ such that $\lambda_1,\ldots,\lambda_s$ are mutually disjoint and 
 $U_{\{\lambda_1,\ldots,\lambda_s\}}\neq\varnothing$.
 It is stable under the diagonal action of $J$ on $I(\mathscr{L})^s$, and the natural surjection
 $I(\mathscr{L})^\sim_s\longrightarrow I(\mathscr{L})_s$; 
 $(\lambda_1,\ldots,\lambda_s)\longmapsto \{\lambda_1,\ldots,\lambda_s\}$ is obviously $J$-equivariant.
 Therefore it suffices to see that the action of $J$ on $I(\mathscr{L})^\sim_s$ has finite orbits.

 Consider the first projection $I(\mathscr{L})^\sim_s\longrightarrow I(\mathscr{L})$,
 which is $J$-equivariant. By Corollary \ref{cor:open-cov-properties-2} i), each fiber of this map is finite.
 Thus, Corollary \ref{cor:open-cov-properties-2} ii) tells us the finiteness of  
 $J$-orbits in $I(\mathscr{L})^\sim_s$.
\end{prf}

For an integer $m\ge 0$ and $\lambda\in I(\mathscr{L})$,
let $\mathscr{U}_{\lambda,m}$ (resp.\ $U_{\lambda,m}$) be the inverse image of $\mathscr{U}_\lambda$
(resp.\ $U_\lambda$) under $\M^\flat_{\mathscr{L},m}/p^\Z\longrightarrow \M^\flat_{\mathscr{L}}/p^\Z$
(resp.\ $M_{\mathscr{L},m}/p^\Z\longrightarrow M_{\mathscr{L}}/p^\Z$, where we put 
$M_{\mathscr{L},m}=\M^{\flat,\rig}_{\mathscr{L},m}$).
These are stable under the natural action of $N_{\mathscr{L}}\times J_{\lambda}/p^\Z$
on $\M^\flat_{\mathscr{L},m}/p^\Z$ and $M_{\mathscr{L},m}/p^\Z$.
Similarly we define $\mathscr{U}_{\underline{\lambda},m}$ and $U_{\underline{\lambda},m}$.

As $\M^\flat_{\mathscr{L},m}$ is locally algebraizable (Remark \ref{rem:locally-algebraizable}), 
\cite[Proposition 4.6]{adicLTF} tells us that the compactly supported cohomology 
\[
 H^i_c(U_{\underline{\lambda},m}):=H^i_c(U_{\underline{\lambda},m}\otimes_{\Q_{p^\infty}}\overline{\Q}_{p^\infty},\Q_\ell)\otimes_{\Q_\ell}\overline{\Q}_\ell
\]
is finite-dimensional. Two groups $N_{\mathscr{L}}$ and $J_{\underline{\lambda}}/p^\Z$ naturally act on it.
The action of $K_{\!\mathscr{L},m}\subset N_{\mathscr{L}}$ is trivial, and
the action of $J_{\underline{\lambda}}/p^\Z$ is known to be smooth (\cite[Corollary 7.7]{MR1262943},
\cite[Corollaire 4.4.7]{MR2074714}).
Therefore, for fixed $g\in N_{\mathscr{L}}$, the function
\[
 \eta^i_{\underline{\lambda},gK_{\!\mathscr{L},m}}\colon h\longmapsto \Tr\bigl((g^{-1},h^{-1});H^i_c(U_{\underline{\lambda},m})\bigr)
\]
on $J_{\underline{\lambda}}/p^\Z$ only depends on the coset $gK_{\!\mathscr{L},m}$ and is locally constant.
Extending it by $0$, we regard it as a locally constant compactly supported function on $J/p^\Z$.

The following construction is very important in this work.

\begin{defn}\label{defn:eta-function}
 For each integer $s\ge 0$, take a system of representatives 
 $\underline{\lambda}_{s,1},\ldots,\underline{\lambda}_{s,k_s}$ of the quotient
 $J\backslash I(\mathscr{L})_{s+1}$. For $m\ge 0$ and $g\in N_\mathscr{L}$,
 we define a locally constant compactly supported function
 $\eta_{gK_{\!\mathscr{L},m}}$ on $J/p^\Z$ as follows:
 \[
  \eta_{gK_{\!\mathscr{L},m}}=\sum_{s,t\ge 0}\sum_{i=1}^{k_s}\frac{(-1)^{s+t}}{\vol(J_{\underline{\lambda}_{s,i}}/p^\Z)}\eta^t_{\underline{\lambda}_{s,i},gK_{\!\mathscr{L},m}}.
 \]
\end{defn}

\begin{prop}\label{prop:eta-property}
 Let $m\ge 0$ be an integer such that $K_{\!\mathscr{L},m}\subset K_0$, and $g$ an element of $N_{\mathscr{L}}$.
 For every admissible representation $\rho$ of $J/p^\Z$, we have
 \[
 \Tr(\eta_{gK_{\!\mathscr{L},m}};\rho)=\Tr(g;H_{\mathrm{RZ}}[\rho]^{K_{\!\mathscr{L},m}})
 \]
 (for a definition of $H_{\mathrm{RZ}}[\rho]$, see Definition \ref{defn:H-RZ-Ext}
 and Definition \ref{defn:H-RZ-alternating}).
 Moreover, the image of $\eta_{gK_{\!\mathscr{L},m}}$ in $\overline{\mathcal{H}}(J/p^\Z)$
 is characterized by this property. In particular, the image is independent of the choice of
 $\underline{\lambda}_{s,1},\ldots,\underline{\lambda}_{s,k_s}$.
\end{prop}

\begin{prf}
 By Remark \ref{rem:RZ-rho-K}, Lemma \ref{lem:mod-center} and Proposition \ref{prop:formal-model-gen-fiber},
 we obtain the equality
 \[
  \Tr(g;H_{\mathrm{RZ}}[\rho]^{K_{\!\mathscr{L},m}})=\sum_{i,j\ge 0}(-1)^{i+j}\Tr\bigl(g;\Ext^j_{J/p^\Z}(H_c^i(M_{\mathscr{L},m}/p^\Z),\rho)\bigr).
 \]
 On the other hand, we have an $N_{\mathscr{L}}\times J/p^\Z$-equivariant \v{C}ech spectral sequence
 \[
  E_1^{-s,t}=\bigoplus_{\underline{\lambda}\in I(\mathscr{L})_{s+1}}H^t_c(U_{\underline{\lambda},m})
 \Longrightarrow H^{-s+t}_c(M_{\mathscr{L},m}/p^\Z).
 \]
 It is easy to see that $E_1^{-s,t}\cong \bigoplus_{i=1}^{k_s}\cInd_{J_{\underline{\lambda}_{s,i}}}^J H^t_c(U_{\underline{\lambda}_{s,i},m})$ as $N_{\mathscr{L}}\times J/p^\Z$-representations.
 Since $\cInd_{J_{\underline{\lambda}_{s,i}}}^J H^t_c(U_{\underline{\lambda}_{s,i},m})$ is projective
 in the category of smooth $J/p^\Z$-modules, 
 in the Grothendieck group of finite-dimensional representations of $N_{\mathscr{L}}$
 we can compute as follows:
 \begin{align*}
  &\sum_{i,j\ge 0}(-1)^{i+j}\Ext^j_{J/p^\Z}\bigl(H_c^i(M_{\mathscr{L},m}/p^\Z),\rho\bigr)\\
  &\qquad =\sum_{s,t,j\ge 0}\sum_{i=1}^{k_s}(-1)^{s+t+j}\Ext^j_{J/p^\Z}\bigl(\cInd_{J_{\underline{\lambda}_{s,i}}}^J H_c^t(U_{\underline{\lambda}_{s,i},m}),\rho\bigr)\\
  &\qquad=\sum_{s,t\ge 0}\sum_{i=1}^{k_s}(-1)^{s+t}\Hom_{J_{\underline{\lambda}_{s,i}}/p^\Z}\bigl(H_c^t(U_{\underline{\lambda}_{s,i},m}),\rho\bigr).
 \end{align*}
 Therefore,
 \begin{align*}
  \Tr(g;H_{\mathrm{RZ}}[\rho]^{K_{\!\mathscr{L},m}})&=\sum_{s,t\ge 0}\sum_{i=1}^{k_s}(-1)^{s+t}\Tr\bigl(g;\Hom_{J_{\underline{\lambda}_{s,i}}/p^\Z}(H_c^t(U_{\underline{\lambda}_{s,i},m}),\rho)\bigr)\\
  &=\sum_{s,t\ge 0}\sum_{i=1}^{k_s}\frac{(-1)^{s+t}}{\vol(J_{\underline{\lambda}_{s,i}}/p^\Z)}\int_{J_{\underline{\lambda}_{s,i}}/p^\Z}\eta^t_{\underline{\lambda}_{s,i},gK_{\!\mathscr{L},m}}(h)\Tr\rho(h)dh\\
  &=\int_{J_{\underline{\lambda}_{s,i}}/p^\Z}\eta_{gK_{\!\mathscr{L},m}}(h)\Tr\rho(h)dh\\
  &=\Tr(\eta_{gK_{\!\mathscr{L},m}};\rho),
 \end{align*}
 as desired.
 
 The uniqueness of an element of $\overline{\mathcal{H}}(J/p^\Z)$ with this property 
 follows from \cite[Theorem 0]{MR874042}.
\end{prf}

\subsection{Action of regular element}\label{subsec:action-of-regular-element}
Let $\mathscr{L}$ and $m\ge 0$ be as in the previous section.
We fix a regular element $\gamma\in J$ and consider how $\gamma$ permutes the open subsets
$\{U_\lambda\}_{\lambda\in I(\mathscr{L})}$.

First recall the following well-known lemma:

\begin{lem}\label{lem:centralizer}
 Let $\mathbf{H}$ be a connected reductive group over $\Q_p$, $\mathbf{Z}_{\mathbf{H}}$ its center
 and $\gamma\in \mathbf{H}(\Q_p)$ a regular element.
 Put $H=\mathbf{H}(\Q_p)$, $Z_H=\mathbf{Z}_{\mathbf{H}}(\Q_p)$ and $Z(\gamma)=\mathbf{Z}(\gamma)(\Q_p)$.
 Then the map $Z(\gamma)\backslash H\longrightarrow H/Z_H$; $h\longmapsto h^{-1}\gamma h$ is proper.
\end{lem}

\begin{prf}
 We have only to prove that the inverse image of any compact subset of $H/Z_H$ is compact.
 It follows from \cite[Lemma 18]{MR0414797}.
\end{prf}

\begin{prop}\label{prop:I_gamma}
 Let $I(\mathscr{L})_\gamma$ be the subset $\{\lambda\in I(\mathscr{L})\mid \gamma U_{\lambda}\cap U_\lambda\neq \varnothing\}$ of $I(\mathscr{L})$.
 Then the left action of $Z(\gamma)$ on $I(\mathscr{L})$ preserves $I(\mathscr{L})_\gamma$ and 
 $Z(\gamma)\backslash I(\mathscr{L})_\gamma$ is finite.
 If $\gamma$ is elliptic, $I(\mathscr{L})_\gamma$ is a finite set.
\end{prop}

\begin{prf}
 Take a system of representatives $\lambda_1,\ldots,\lambda_k$ of $J\backslash I(\mathscr{L})$.
 Then, $I(\mathscr{L})$ can be identified with $\coprod_{i=1}^k J/J_{\lambda_i}$.
 Since $J_{\lambda_i}/p^\Z$ is a compact open subgroup of $J/p^\Z$ (Corollary \ref{cor:open-cov-properties-2} iv)),
 to show the finiteness of $Z(\gamma)\backslash I(\mathscr{L})_\gamma$,
 it suffices to show that the set $\{h\in Z(\gamma)\backslash J\mid \gamma U_{h\lambda_i}\cap U_{h\lambda_i}\neq \varnothing\}$
 is compact. The condition 
 $\gamma U_{h\lambda_i}\cap U_{h\lambda_i}\neq \varnothing$ is equivalent to
 $U_{h^{-1}\gamma h\lambda_i}\cap U_{\lambda_i}\neq \varnothing$.
 By Corollary \ref{cor:iii-compactness}, the set
 $\{h'\in J/p^\Z\mid U_{h'\lambda_i}\cap U_{\lambda_i}\neq \varnothing\}$ is
 a compact subset of $J/p^\Z$.
 Therefore, Lemma \ref{lem:centralizer} tells us that the set
 $\{h\in Z(\gamma)\backslash J\mid U_{h^{-1}\gamma h\lambda_i}\cap U_{\lambda_i}\neq \varnothing\}$ is
 compact.

 Assume that $\gamma$ is elliptic, namely, $Z(\gamma)/p^\Z$ is compact.
 Then, for each $\lambda\in I(\mathscr{L})_\gamma$, the $Z(\gamma)$-orbit of $\lambda$ is finite;
 indeed, the $Z(\gamma)$-orbit can be identified with $(Z(\gamma)/p^\Z)/((Z(\gamma)\cap J_\lambda)/p^\Z)$,
 and $(Z(\gamma)\cap J_\lambda)/p^\Z$ is an open subgroup of $Z(\gamma)/p^\Z$.
 In other words, each fiber of the natural projection 
 $I(\mathscr{L})_\gamma\longrightarrow Z(\gamma)\backslash I(\mathscr{L})_\gamma$ is a finite set.
 Hence $I(\mathscr{L})_\gamma$ is a finite set.
\end{prf}

\begin{defn}
 Put 
 \[
  \mathscr{U}_{\gamma,m}=\bigcup_{\lambda\in I(\mathscr{L})_\gamma}\mathscr{U}_{\lambda,m},\quad
  U_{\gamma,m}=\bigcup_{\lambda\in I(\mathscr{L})_\gamma}U_{\lambda,m}.
 \]
 These are stable under the actions of $N_{\mathscr{L}}$ and $Z(\gamma)$.
\end{defn}

\begin{cor}\label{cor:I'}
 Let $I(\mathscr{L})'_\gamma$ be the subset of $I(\mathscr{L})$ consisting of $\lambda$ such that
 $U_{\lambda,m}\cap U_{\gamma,m}\neq \varnothing$.
 Then the left action of $Z(\gamma)$ on $I(\mathscr{L})$ preserves $I(\mathscr{L})'_\gamma$ and 
 $Z(\gamma)\backslash I(\mathscr{L})'_\gamma$ is finite.
 If $\gamma$ is elliptic, $I(\mathscr{L})'_\gamma$ is a finite set.
\end{cor}

\begin{prf}
 It is an easy consequence of Proposition \ref{prop:I_gamma} and Corollary \ref{cor:open-cov-properties-2} i).
 Indeed, for a system of representatives $\lambda_1,\ldots,\lambda_k$ for 
 $Z(\gamma)\backslash I(\mathscr{L})_\gamma$,
 $Z(\gamma)\backslash I(\mathscr{L})'_\gamma$ is contained in the image of a finite set
 $\bigcup_{i=1}^k\{\lambda\in I(\mathscr{L})\mid U_\lambda\cap U_{\lambda_i}\neq\varnothing\}$.
\end{prf}

\begin{defn}
 Put 
 \[
  \mathscr{U}'_{\gamma,m}=\bigcup_{\lambda\in I(\mathscr{L})'_\gamma}\mathscr{U}_{\lambda,m},\quad
  U'_{\gamma,m}=\bigcup_{\lambda\in I(\mathscr{L})'_\gamma}U_{\lambda,m}.
 \]
 These are stable under the actions of $N_{\mathscr{L}}$ and $Z(\gamma)$.
\end{defn}

\begin{rem}\label{rem:bigger-open}
 The closure $\overline{U_{\gamma,m}}$ of $U_{\gamma,m}$ is contained in
 $U'_{\gamma,m}$.
 Indeed, let $x$ be a point in $\overline{U_{\gamma,m}}$ and take $\lambda\in I(\mathscr{L})$
 such that $x\in U_{\lambda,m}$.
 Then $U_{\lambda,m}$ should be intersect $U_{\gamma,m}$ and therefore $\lambda$ lies in $I(\mathscr{L})'_\gamma$.
\end{rem}

If $\gamma$ is elliptic, then $\mathscr{U}_{\gamma,m}$ is quasi-compact,
and thus the cohomology $H^i_c(U_{\gamma,m})$ is finite-dimensional.
For $g\in N_{\mathscr{L}}$, the alternating sum of the traces of $(g^{-1},\gamma^{-1})$ on $H^i_c(U_{\gamma,m})$
can be computed by the function $\eta_{gK_{\!\mathscr{L},m}}$ introduced in the previous subsection:

\begin{prop}\label{prop:trace-orbital-integral}
 Assume that $\gamma$ is elliptic. Then, for $g\in N_{\mathscr{L}}$ we have
 \[
 \sum_i(-1)^i\Tr\bigl((g^{-1},\gamma^{-1});H^i_c(U_{\gamma,m})\bigr)
 =O_\gamma(\eta_{gK_{\!\mathscr{L},m}}).
 \]
 For our normalization of the Haar measure, see the last paragraph of Section \ref{subsec:def-RZ}.
\end{prop}

\begin{prf}
 For an integer $s\ge 1$, put $I(\mathscr{L})_{\gamma,s}=\{\underline{\lambda}\in I(\mathscr{L})_s\mid \underline{\lambda}\subset I(\mathscr{L})_\gamma\}$. Then we have the \v{C}ech spectral sequence
 \[
  E_1^{-s,t}=\bigoplus_{\underline{\lambda}\in I(\mathscr{L})_{\gamma,s+1}}H^t_c(U_{\underline{\lambda},m})
 \Longrightarrow H^{-s+t}_c(U_{\gamma,m}).
 \]
 Therefore, we can compute
 \begin{align*}
  &\sum_i(-1)^i\Tr\bigl((g^{-1},\gamma^{-1});H^i_c(U_{\gamma,m})\bigr)
  =\sum_{s,t\ge 0}\sum_{\underline{\lambda}\in I(\mathscr{L})_{\gamma,s+1}}(-1)^{s+t}\Tr\bigl((g^{-1},\gamma^{-1});H^t_c(U_{\underline{\lambda},m})\bigr)\\
  &\qquad=\sum_{s,t\ge 0}\sum_{i=1}^{k_s}(-1)^{s+t}\sum_{\substack{h\in J/J_{\underline{\lambda}_{s,i}}\\ h\underline{\lambda}_{s,i}\subset I(\mathscr{L})_\gamma}}\Tr\bigl((g^{-1},\gamma^{-1});H^t_c(U_{h\underline{\lambda}_{s,i},m})\bigr)\\
  &\qquad=\sum_{s,t\ge 0}\sum_{i=1}^{k_s}(-1)^{s+t}\sum_{\substack{h\in J/J_{\underline{\lambda}_{s,i}}\\ h\underline{\lambda}_{s,i}\subset I(\mathscr{L})_\gamma}}\eta^t_{\underline{\lambda}_{s,i},gK_{\!\mathscr{L},m}}(g,h^{-1}\gamma h)\\
  &\qquad\stackrel{(*)}{=}\sum_{s,t\ge 0}\sum_{i=1}^{k_s}\frac{(-1)^{s+t}}{\vol(J_{\underline{\lambda}_{s,i}}/p^\Z)}\int_{J/p^\Z}\eta^t_{\underline{\lambda}_{s,i},gK_{\!\mathscr{L},m}}(g,h^{-1}\gamma h)dh=O_\gamma(\eta_{gK_{\!\mathscr{L},m}}).
 \end{align*}
 For $(*)$, note that if $h\underline{\lambda}_{s,i}\nsubset I(\mathscr{L})_\gamma$, then
 $\gamma h\underline{\lambda}_{s,i}\neq h\underline{\lambda}_{s,i}$; indeed, 
 if an element $\lambda$ in $h\underline{\lambda}_{s,i}\setminus I(\mathscr{L})_\gamma$ satisfies
 $\gamma\lambda\in h\underline{\lambda}_{s,i}$, we have $U_{h\underline{\lambda}_{s,i}}\subset U_\lambda\cap U_{\gamma\lambda}=\varnothing$, which contradicts to the fact $h\underline{\lambda}_{s,i}\in I(\mathscr{L})_{s+1}$.
 Thus $h\underline{\lambda}_{s,i}\nsubset I(\mathscr{L})_\gamma$ implies
 $h^{-1}\gamma h\notin J_{\underline{\lambda}_{s,i}}$ and
 $\eta^t_{\underline{\lambda}_{s,i},gK_{\!\mathscr{L},m}}(g,h^{-1}\gamma h)=0$.
\end{prf}

Next we consider the case where $\gamma$ is not elliptic.
In this case, the centralizer $\mathbf{Z}(\gamma)$ is a maximal torus which is not anisotropic modulo the center
of $\mathbf{J}$. 
We can take a discrete torsion-free cocompact subgroup $\Gamma$ of $Z(\gamma)/p^\Z$.

Take a system of representatives $\lambda'_1,\ldots,\lambda'_l$ for $Z(\gamma)\backslash I(\mathscr{L})'_\gamma$.
For each $i$ with $1\le i\le l$, let $C_i$ be the subset of $J/p^\Z$ consisting of $h\in J/p^\Z$ such that
there exist $\mu_1,\mu_2,\mu_3\in I(\mathscr{L})$ with 
$U_{\lambda'_i}\cap U_{\mu_1}\neq \varnothing$, $U_{\mu_1}\cap U_{\mu_2}\neq\varnothing$,
$U_{\mu_2}\cap U_{\mu_3}\neq\varnothing$ and $U_{\mu_3}\cap U_{h\lambda'_i}\neq \varnothing$.
By Corollary \ref{cor:open-cov-properties-2} i) and Corollary \ref{cor:iii-compactness}, 
$C_i$ is compact. Put $C=\bigcup_{i=1}^l C_i$.

Since $\Gamma$ is discrete, $C\cap \Gamma$ is finite. Therefore, by shrinking $\Gamma$ if necessary,
we may assume that $C\cap \Gamma=\{1\}$. 

\begin{lem}\label{lem:I'-finite}
 The quotient $\Gamma\backslash I(\mathscr{L})'_\gamma$ is a finite set.
\end{lem}

\begin{prf}
 As in the proof of Proposition \ref{prop:I_gamma}, it follows from compactness of $Z(\gamma)/\Gamma$
 and Corollary \ref{cor:I'}.
\end{prf}

\begin{lem}\label{lem:Gamma-tot-discont}
 For $\lambda\in I(\mathscr{L})'_\gamma$ and $\gamma'\in\Gamma$, assume that
 there exist $\mu_1,\mu_2,\mu_3\in I(\mathscr{L})$ such that
 $U_\lambda\cap U_{\mu_1}\neq \varnothing$, $U_{\mu_1}\cap U_{\mu_2}\neq\varnothing$,
 $U_{\mu_2}\cap U_{\mu_3}\neq\varnothing$ and $U_{\mu_3}\cap U_{\gamma'\lambda}\neq \varnothing$.
 Then $\gamma'=1$.
 
 In particular, if $\gamma'\neq 1$ then we have $U_{\gamma'\lambda}\cap U_\lambda=\varnothing$
 for every $\lambda\in I(\mathscr{L})'_\gamma$.
\end{lem}

\begin{prf}
 We may write $\lambda=h\lambda'_i$ with $h\in Z(\gamma)$ and $1\le i\le l$.
 If there exist $\mu_1,\mu_2,\mu_3$ as above, then
 $U_{\lambda'_i}\cap U_{h^{-1}\mu_1}\neq \varnothing$, $U_{h^{-1}\mu_1}\cap U_{h^{-1}\mu_2}\neq\varnothing$,
 $U_{h^{-1}\mu_2}\cap U_{h^{-1}\mu_3}\neq\varnothing$ and $U_{h^{-1}\mu_3}\cap U_{h^{-1}\gamma'h\lambda'_i}\neq \varnothing$, and thus $h^{-1}\gamma'h\in C_i$. Since $\gamma',h\in Z(\gamma)$ and $Z(\gamma)$ is abelian,
 we have $\gamma'=h^{-1}\gamma'h$. Therefore $\gamma'\in C\cap \Gamma=\{1\}$, as desired.
\end{prf}

By the lemma above, we can take quotient of $\mathscr{U}_{\gamma,m}$, $\mathscr{U}'_{\gamma,m}$,
$U_{\gamma,m}$ and $U'_{\gamma,m}$ by $\Gamma$.
Lemma \ref{lem:I'-finite} tells us that $\Gamma\backslash\mathscr{U}_{\gamma,m}$ is quasi-compact,
and thus the cohomology $H^i_c(\Gamma\backslash U_{\gamma,m})$ is finite-dimensional.
We write $p_\Gamma$ for the natural maps
$\mathscr{U}'_{\gamma,m}\longrightarrow \Gamma\backslash\mathscr{U}'_{\gamma,m}$ and
$U'_{\gamma,m}\longrightarrow \Gamma\backslash U'_{\gamma,m}$.

\begin{prop}\label{prop:quotient-property}
 \begin{enumerate}
  \item Let $g$ be an element of $N_\mathscr{L}$.
	Assume that a point $x\in \Gamma\backslash U_{\gamma,m}$ is fixed by $(g,\gamma)$.
	Then, every point $y\in U_{\gamma,m}$ satisfying $p_\Gamma(y)=x$ is fixed by $(g,\gamma)$.
  \item For $\lambda\in I(\mathscr{L})'_\gamma\setminus I(\mathscr{L})_\gamma$,
	we have $p_\Gamma(U_{\gamma\lambda,m})\cap p_\Gamma(U_{\lambda,m})=\varnothing$.
 \end{enumerate}
\end{prop}

\begin{prf}
 i) Take $\lambda\in I(\mathscr{L})_\gamma$ such that
 $y\in U_{\lambda,m}$. As $x$ is fixed by $(g,\gamma)$, there exists $\gamma'\in\Gamma$
 such that $(g,\gamma'\gamma)y=y$. In particular $U_{\gamma'\gamma\lambda,m}\cap U_{\lambda,m}\neq \varnothing$.
 Since $\lambda\in I(\mathscr{L})_\gamma$, we have $U_{\gamma\lambda,m}\cap U_{\lambda,m}\neq \varnothing$,
 and thus $U_{\gamma'\gamma\lambda,m}\cap U_{\gamma'\lambda,m}\neq \varnothing$.
 Therefore, by Lemma \ref{lem:Gamma-tot-discont},
 we can conclude that $\gamma'=1$. Hence $y$ is fixed by $(g,\gamma)$.

 ii) Assume that $p_\Gamma(U_{\gamma\lambda,m})\cap p_\Gamma(U_{\lambda,m})$ is non-empty.
 Then there exists $\gamma'\in\Gamma$ such that $U_{\gamma\lambda,m}\cap U_{\gamma'\lambda,m}\neq\varnothing$.
 Since $\lambda\in I(\mathscr{L})'_\gamma$, we can find $\mu\in I(\mathscr{L})_\gamma$ such that
 $U_{\lambda,m}\cap U_{\mu,m}\neq\varnothing$. By the definition of $I(\mathscr{L})_\gamma$,
 we have $U_{\mu,m}\cap U_{\gamma\mu,m}\neq\varnothing$.
 Therefore, all of $U_{\lambda,m}\cap U_{\mu,m}$, $U_{\mu,m}\cap U_{\gamma\mu,m}$,
 $U_{\gamma\mu,m}\cap U_{\gamma\lambda,m}$ and $U_{\gamma\lambda,m}\cap U_{\gamma'\lambda,m}$ are
 non-empty. Lemma \ref{lem:Gamma-tot-discont} tells us that $\gamma'=1$.
 Hence $U_{\gamma\lambda,m}\cap U_{\lambda,m}\neq\varnothing$, which contradicts to the assumption
 $\lambda\notin I(\mathscr{L})_\gamma$.
\end{prf}

By the same argument as in the proof of Proposition \ref{prop:trace-orbital-integral},
we can obtain the following:

\begin{prop}\label{prop:trace-orbital-integral-non-ell}
 For $g\in N_{\mathscr{L}}$, we have
 \[
 \sum_i(-1)^i\Tr\bigl((g^{-1},\gamma^{-1});H^i_c(\Gamma\backslash U_{\gamma,m})\bigr)
 =\vol\bigl(Z(\gamma)/\widetilde{\Gamma}\bigr) O_\gamma(\eta_{gK_{\!\mathscr{L},m}}),
 \]
 where $\widetilde{\Gamma}$ denotes the inverse image of $\Gamma$ under $J\longrightarrow J/p^\Z$. 
\end{prop}

\section{Application of Lefschetz trace formula}\label{section:LTF}
From now on, we assume that $d=2$.
In this section, we will apply the Lefschetz trace formula 
to compute the left hand side of the identities in Proposition \ref{prop:trace-orbital-integral}
and Proposition \ref{prop:trace-orbital-integral-non-ell}.
Let $\mathscr{L}$ be one of $\mathscr{L}_0$, $\mathscr{L}_{\mathrm{para}}$ or $\mathscr{L}_{\mathrm{Siegel}}$
in Lemma \ref{lem:lattice-chain-GSp4}, and $m\ge 1$ be an integer.
In this case, $K_{\!\mathscr{L},m}$ is contained in $K_0$.
We put $K=K_{\!\mathscr{L},m}$ for simplicity.
The goal of this section is as follows:

\begin{thm}\label{thm:LTF-consequence}
 For $g\in N_{\mathscr{L}}$, assume that $gK$ consists of regular elliptic elements.
 \begin{enumerate}
  \item For every regular elliptic element $\gamma\in J$, we have
	\begin{align*}
	 O_\gamma(\eta_{gK})&=\#\Fix\bigl((g,\gamma);M_{\mathscr{L},m}/p^\Z\bigr)=\sum_{x\in\Fix(\gamma;\Omega)}O_{g_{\gamma,x}}\Bigl(\frac{\mathbf{1}_{gKp^\Z}}{\vol(K)}\Bigr),
	\end{align*}
	(The latter equality has been proved in Theorem \ref{thm:fixed-point-number} and .)
  \item For every regular non-elliptic element $\gamma\in J$, we have $O_\gamma(\eta_{gK})=0$.
 \end{enumerate}
\end{thm}

Before proving the theorem, we record the following corollary, which is an immediate consequence of
Theorem \ref{thm:LTF-consequence} i) and Corollary \ref{cor:Fix-SO}:

\begin{cor}\label{cor:LTF-SO}
 Let $g$ be an element of $N_{\mathscr{L}}$ such that
 $gK$ consists of regular elliptic elements.
 For $\gamma\in J^{\mathrm{ell}}$,
 let $g_\gamma$ be an element of $G$ with $g_\gamma\leftrightarrow \gamma$.
 Then we have
 \[
  \textit{SO}_\gamma(\eta_{gK})=4\textit{SO}_{g_\gamma}\Bigl(\frac{\mathbf{1}_{gKp^\Z}}{\vol(K)}\Bigr).
 \]
\end{cor}

\subsection{Complement on \cite{adicLTF}}
We would like to apply \cite[Theorem 4.5]{adicLTF} to $\mathscr{U}_{\gamma,m}$ and
$\Gamma\backslash \mathscr{U}_{\gamma,m}$, but the rigid generic fiber of these formal schemes
are not partially proper over $\Q_{p^\infty}$.
Here we will give a slightly stronger version of \cite[Theorem 4.5]{adicLTF} which is applicable to our cases.
All techniques we need are included in \cite[\S 4]{adicLTF}.

We use the same notation as in \cite[\S 4]{adicLTF}.
Let $R$ be a complete discrete valuation ring and $k$ an algebraic closure of the fraction field $F$ of $R$.
Put $\mathcal{S}=\Spf R$ and $S=\Spa(k,k^+)$, where $k^+$ is the valuation ring of $k$.
Let $\mathcal{X}$ be a quasi-compact special formal scheme which is separated over $\mathcal{S}$.
Then we can associate $\mathcal{X}$ with the adic spaces $X=t(\mathcal{X})_a$, $X_\eta=t(\mathcal{X})_{\eta}$
and $X_{\overline{\eta}}=t(\mathcal{X})_{\overline{\eta}}$. We denote the special fiber of $\mathcal{X}$
(resp.\ $X$) by $\mathcal{X}_s$ (resp.\ $X_s$).

Let $\mathcal{X}'$ be a quasi-compact special formal scheme separated over $\mathcal{S}$ which contains
$\mathcal{X}$ as an open formal subscheme. We can define $X'$, $X'_\eta$, $X'_{\overline{\eta}}$,
$\mathcal{X}'_s$ and $X'_s$ similarly.
Let $\mathcal{T}$ be a finite set equipped with a partial order
and $\{\mathcal{Y}_\alpha\}_{\alpha\in\mathcal{T}}$ a family of closed formal subschemes
of $\mathcal{X}_s$ indexed by $\mathcal{T}$. We put 
$Y_\alpha=t(\mathcal{Y}_{\alpha})_a=t(\mathcal{Y}_\alpha)\times_{t(\mathcal{X})}X$,
and assume the same condition as in \cite[Assumption 4.1]{adicLTF};

\begin{assump}\label{assump:4.1}
  \begin{enumerate}
  \item $X_s=\bigcup_{\alpha\in\mathcal{T}}Y_{\alpha}$.
  \item For $\alpha\in\mathcal{T}$, put $Y_{(\alpha)}=Y_\alpha\setminus\bigcup_{\beta>\alpha}Y_\beta$. Then, for $\alpha,\beta\in\mathcal{T}$ with $\alpha\neq \beta$,
	$Y_{(\alpha)}\cap Y_{(\beta)}=\varnothing$.
 \end{enumerate}
\end{assump}

Let $f\colon \mathcal{X}'\longrightarrow \mathcal{X}'$ be an isomorphism preserving $\mathcal{X}$.
We also denote the induced isomorphism $X'\yrightarrow{\cong}X'$ by the same symbol $f$.
The induced isomorphisms $X'_{\eta}\yrightarrow{\cong}X'_{\eta}$ and 
$X'_{\overline{\eta}}\yrightarrow{\cong}X'_{\overline{\eta}}$ are
denoted by $f_{\eta}$ and $f_{\overline{\eta}}$,
respectively. We will make the same assumption as in \cite[Assumption 4.3]{adicLTF};

\begin{assump}\label{assump:4.3}
 There exist an order-preserving bijection $f\colon \mathcal{T}\yrightarrow{\cong}\mathcal{T}$
 and a system of closed constructible subsets $\{Y_\alpha(n)\}_{n\ge 1}$
 of $X$ for each $\alpha\in \mathcal{T}$ satisfying the following:
 \begin{enumerate}
  \item $Y_\alpha(n+1)\subset Y_\alpha(n)$ for every $n\ge 1$. 
  \item $\bigcap_{n\ge 1}Y_\alpha(n)=Y_\alpha$.
  \item $f(Y_\alpha(n))=Y_{f(\alpha)}(n)$ for every $\alpha\in\mathcal{T}$ and $n\ge 1$.
  \item $f(\alpha)\neq \alpha$ for every $\alpha\in\mathcal{T}$.
 \end{enumerate}
\end{assump}

\begin{rem}
 As proved in \cite[Proposition 4.18]{adicLTF}, if $f$ induces an isomorphism of formal schemes $\mathcal{Y}_\alpha\yrightarrow{\cong} \mathcal{Y}_{f(\alpha)}$ for each $\alpha\in\mathcal{T}$,
then we can find a system of closed constructible subsets $\{Y_\alpha(n)\}_{n\ge 1}$
 of $X$ satisfying Assumption \ref{assump:4.3} i), ii), iii).
\end{rem}

\begin{thm}\label{thm:LTF-formal-stronger}
 In addition to Assumption \ref{assump:4.1} and Assumption \ref{assump:4.3}, assume the following:
 \begin{enumerate}
  \item[(a)] $\mathcal{X}'$ is locally algebraizable (\cite[Definition 3.18]{formalnearby}) and
	$X_{\eta}$ is smooth over $\Spa(F,R)$. 
  \item[(b)] The closure $\overline{X_\eta}$ of $X_\eta$ in $X'_\eta$ is partially proper over $\Spa(F,R)$.
  \item[(c)] For every $x\in X'_{\overline{\eta}}\setminus X_{\overline{\eta}}$, $f_{\overline{\eta}}(x)\neq x$.
 \end{enumerate}
 Let $\Lambda=\Z/\ell^n\Z$ for a prime $\ell$ which is invertible in $k^+$ and some integer $n\ge 1$.
 Then, $\Fix (f_{\overline{\eta}}\vert_{X_{\overline{\eta}}})$ (\cf \cite[Example 2.9]{adicLTF}) is proper over $S$ and we have
 \[
  \Tr\bigl(f_{\overline{\eta}}^*;R\Gamma_c(X_{\overline{\eta}},\Lambda)\bigr)=\#\Fix (f_{\overline{\eta}}\vert_{X_{\overline{\eta}}}).
 \]
 For the definition of the right hand side, see \cite[Definition 2.6]{adicLTF}.
\end{thm}

\begin{prf}
 First of all, \cite[Proposition 4.6]{adicLTF} tells us that $R\Gamma_c(X_{\overline{\eta}},\Lambda)$
 is a perfect $\Lambda$-complex, and therefore we can consider the trace 
 $\Tr(f_{\overline{\eta}}^*;R\Gamma_c(X_{\overline{\eta}},\Lambda))$.

 Note that the proof of \cite[Lemma 4.8]{adicLTF} does not require partially properness of $X_\eta$;
 thus we can find an integer $n_\alpha\ge 1$ for each $\alpha\in\mathcal{T}$
 satisfying the following conditions:
\begin{itemize}
 \item For every $\alpha\in\mathcal{T}$, $n_\alpha=n_{f(\alpha)}$.
 \item For $\alpha\in \mathcal{T}$, 
	put $U_\alpha=Y_\alpha(n_\alpha)\setminus \bigcup_{\beta>\alpha}Y_{\beta}(n_\beta)$.
	Then we have $U_\alpha\cap U_\beta=\varnothing$ for every $\alpha,\beta\in\mathcal{T}$
	with $\alpha\neq \beta$.
\end{itemize}
Put $W=\bigcup_{\alpha\in\mathcal{T}}Y_\alpha(n_\alpha)$ and
$X_0=X\setminus W$. As in the proof of \cite[Lemma 4.12 i)]{adicLTF},
we can show that $X_0$ is a quasi-compact open adic subspace of $X_\eta$.
Moreover, by exactly the same method as in the proof of \cite[Proposition 4.10]{adicLTF},
we can obtain the equality 
 $\Tr(f^*_{\overline{\eta}};R\Gamma_c(X_{\overline{\eta}},\Lambda))=\Tr(f^*_{\overline{\eta}};R\Gamma_c(X_{0,\overline{\eta}},\Lambda))$.

 Consider the closure $\overline{X_{0,\overline{\eta}}}$ of $X_{0,\overline{\eta}}$ in $X'_{\overline{\eta}}$.
 Since $X'_{\overline{\eta}}$ is taut (\cf \cite[Lemma 4.14]{formalnearby}), $\overline{X_{0,\overline{\eta}}}$
 is quasi-compact. On the other hand, by the assumption (b), $\overline{X_{0,\overline{\eta}}}$
 is partially proper over $S$. Thus $\overline{X_{0,\overline{\eta}}}$ is proper over $S$.
 Let us observe that for $x\in X'_{\overline{\eta}}\setminus X_{0,\overline{\eta}}$ we have
 $f_{\overline{\eta}}(x)\neq x$. Indeed, if $x\notin X_{\overline{\eta}}$, by the assumption (c)
 we have $f_{\overline{\eta}}(x)\neq x$.
 If $x\in X_{\overline{\eta}}\setminus X_{0,\overline{\eta}}$, 
 then $x\in W_{\overline{\eta}}=\coprod_{\alpha\in\mathcal{T}}U_{\alpha,\overline{\eta}}$, and thus
 we can find $\alpha\in\mathcal{T}$ such that $x\in U_{\alpha,\overline{\eta}}$.
 Since $U_{\alpha,\overline{\eta}}\cap f(U_{\alpha,\overline{\eta}})=U_{\alpha,\overline{\eta}}\cap U_{f(\alpha),\overline{\eta}}=\varnothing$, $x$ and $f_{\overline{\eta}}(x)$ are distinct.
 Note that this implies that $\Fix(f_{\overline{\eta}}\vert_{X_{0,\overline{\eta}}})=\Fix(f_{\overline{\eta}}\vert_{\overline{X_{0,\overline{\eta}}}})=\Fix(f_{\overline{\eta}}\vert_{X_{\overline{\eta}}})$.
 In particular $\Fix(f_{\overline{\eta}}\vert_{X_{\overline{\eta}}})$ is proper over $S$.

 Now we can apply the Lefschetz trace formula \cite[Theorem 3.13]{adicLTF} to 
 $X_{0,\overline{\eta}}\hooklongrightarrow \overline{X_{0,\overline{\eta}}}$, and obtain
 \[
 \Tr\bigl(f^*_{\overline{\eta}};R\Gamma_c(X_{\overline{\eta}},\Lambda)\bigr)=\Tr\bigl(f_{\overline{\eta}}^*;R\Gamma_c(X_{0,\overline{\eta}},\Lambda)\bigr)=\#\Fix(f_{\overline{\eta}}\vert_{X_{0,\overline{\eta}}})=
 \#\Fix(f_{\overline{\eta}}\vert_{X_{\overline{\eta}}}).
 \]
 For the final equality, we use \cite[Proposition 2.10]{adicLTF}.
\end{prf}

\begin{rem}
 At least when the characteristic of $k$ is $0$,
 we can deduce from Theorem \ref{thm:LTF-formal-stronger} the analogous result for $\ell$-adic coefficient 
 simply by taking projective limit (\cf \cite[proof of Corollary 4.40]{formalnearby}).
\end{rem}

\subsection{Proof of Theorem \ref{thm:LTF-consequence}}

Here we give a proof of Theorem \ref{thm:LTF-consequence}.

\bigbreak

First we consider the case where $\gamma\in J$ is elliptic. 
We will apply Theorem \ref{thm:LTF-formal-stronger} to $\mathscr{U}_{\gamma,m}\subset \mathscr{U}'_{\gamma,m}$
and $(g,\gamma)\colon \mathscr{U}'_{\gamma,m}\longrightarrow \mathscr{U}'_{\gamma,m}$.
By Remark \ref{rem:locally-algebraizable}, $\mathscr{U}'_{\mathscr{L},m}$ is locally algebraizable
(actually it is algebraizable).
We know that the generic fiber $U_{\gamma,m}$ is smooth over $\Q_{p^\infty}$.
Moreover, Remark \ref{rem:bigger-open} tells us that the closure $\overline{U_{\gamma,m}}$ of $U_{\gamma,m}$
inside $M_{\mathscr{L},m}/p^\Z$ is the same as that inside $U'_{\gamma,m}$.
Since $M_{\mathscr{L},m}/p^\Z$ is partially proper over $\Q_{p^\infty}$, 
so is $\overline{U_{\gamma,m}}$.
Let $x\in U'_{\gamma,m}\setminus U_{\gamma,m}$. 
Then we can find $\lambda\in I(\mathscr{L})_\gamma'\setminus I(\mathscr{L})_\gamma$
such that $x\in U_{\lambda,m}$. 
Since $U_{\lambda,m}\cap (g,\gamma)U_{\lambda,m}=U_{\lambda,m}\cap U_{\gamma\lambda,m}=\varnothing$,
we have $x\neq (g,\gamma)x$.
Thus the assumptions (a), (b), (c) in Theorem \ref{thm:LTF-formal-stronger} are satisfied.

Recall that for each $\alpha\in \mathcal{S}_{\mathscr{L},m}$, 
a closed formal subscheme $\M^\flat_{\mathscr{L},m,\alpha}$
is attached (Definition \ref{defn:boundary-strata}, Lemma \ref{lem:boundary-strata}).
We denote by $\mathcal{Y}_\alpha$ the restriction of $\M^\flat_{\mathscr{L},m,\alpha}/p^\Z$ to
$\mathscr{U}_{\gamma,m}$. 
In Proposition \ref{prop:formal-model-disjoint}, we have checked that
$\{\mathcal{Y}_\alpha\}_{\alpha\in\mathcal{S}_{\mathscr{L},m}}$ satisfies Assumption \ref{assump:4.1}. 
By Lemma \ref{lem:formal-model-group-action}, the isomorphism $(g,\gamma)\colon \mathscr{U}_{\gamma,m}\longrightarrow \mathscr{U}_{\gamma,m}$
induces an isomorphism $\mathcal{Y}_\alpha\longrightarrow \mathcal{Y}_{g\alpha}$.
Therefore, Assumption \ref{assump:4.3} is also satisfied; 
note that Assumption \ref{assump:4.3} iv) is nothing but Corollary \ref{cor:linear-algebra}.

Therefore all the conditions of Theorem \ref{thm:LTF-formal-stronger} are verified,
and we have
\[
 \sum_i(-1)^i\Tr\bigl((g^{-1},h^{-1});H^i_c(U_{\gamma,m})\bigr)=\#\Fix\bigl((g,\gamma);U_{\gamma,m}\bigr).
\]
By Proposition \ref{prop:trace-orbital-integral}, the left hand side is equal to
$O_\gamma(\eta_{gK})$. As there is no fixed point under $(g,\gamma)$ in $(M_{\mathscr{L},m}/p^\Z)\setminus U_{\gamma,m}$, the right hand side is equal to $\#\Fix((g,\gamma);M_{\mathscr{L},m}/p^\Z)$
(\cf \cite[Proposition 2.10]{adicLTF}). This completes the proof of Theorem \ref{thm:LTF-consequence} i).

\bigbreak

The case where $\gamma$ is non-elliptic is similar. Let $\Gamma\subset Z(\gamma)/p^\Z$ be as in
Section \ref{subsec:action-of-regular-element}, and apply Theorem \ref{thm:LTF-formal-stronger}
to $\Gamma\backslash\mathscr{U}_{\gamma,m}\subset \Gamma\backslash\mathscr{U}'_{\gamma,m}$
and $(g,\gamma)\colon \Gamma\backslash\mathscr{U}'_{\gamma,m}\longrightarrow \Gamma\backslash\mathscr{U}'_{\gamma,m}$.
The conditions (a), (b), (c) in Theorem \ref{thm:LTF-formal-stronger} can be proved in the same manner.
(a) is clear. For (b), note that $\Gamma\backslash \overline{U_{\gamma,m}}$ is a closed subset of
$\Gamma\backslash U'_{\gamma,m}$ which is partially proper over $\Q_{p^\infty}$.
Thus the closure of $\Gamma\backslash U_{\gamma,m}$ in $\Gamma\backslash U'_{\gamma,m}$,
being a closed subset of $\Gamma\backslash \overline{U_{\gamma,m}}$, is also partially proper over
$\Q_{p^\infty}$.
(c) can be proved in the same way as above by using Proposition \ref{prop:quotient-property} ii).

It is easy to show that $\mathcal{Y}_\alpha$ for $\alpha\in\mathcal{S}_{\mathscr{L},m}$ induces a
closed formal subscheme $\Gamma\backslash \mathcal{Y}_\alpha$ of $\Gamma\backslash \mathscr{U}_{\gamma,m}$,
and a family $\{\Gamma\backslash\mathcal{Y_\alpha}\}_{\alpha\in \mathcal{S}_{\mathscr{L},m}}$
satisfies Assumption \ref{assump:4.1} and Assumption \ref{assump:4.3}.

Therefore we can apply Theorem \ref{thm:LTF-formal-stronger},
and we have
\[
 \sum_i(-1)^i\Tr\bigl((g^{-1},h^{-1});H^i_c(\Gamma\backslash U_{\gamma,m})\bigr)=\#\Fix\bigl((g,\gamma);\Gamma\backslash U_{\gamma,m}\bigr).
\]
By Proposition \ref{prop:trace-orbital-integral-non-ell}, the left hand side is equal to
$\vol(Z(\gamma)/\widetilde{\Gamma})O_\gamma(\eta_{gK})$.
By Proposition \ref{prop:quotient-property} i), the right hand side is equal to 
$\#\Fix((g,\gamma);U_{\gamma,m})$, which is zero by Theorem \ref{thm:fixed-point-number} ii).
This completes the proof of Theorem \ref{thm:LTF-consequence} ii).

\section{Computation of the character}\label{sec:computation-of-char}

\subsection{Local Langlands correspondence for $G$ and $J$}\label{subsec:LLC-G-J}
A candidate of the local Langlands correspondence for $G$ and $J$ has been constructed by
Gan-Takeda \cite{MR2800725} and Gan-Tantono \cite{Gan-Tantono}. Here we review their results briefly.
In this subsection, let $F$ be a $p$-adic field and put $G=\GSp_4(F)$, $J=\mathrm{GU}(2,D)$ where
$D$ is a quaternion division algebra over $F$.

First we recall basic definitions on $L$-parameters. 

\begin{defn}
 \begin{enumerate}
  \item An $L$-parameter for $G$ is a homomorphism $\phi\colon W_F\times \SL_2(\C)\longrightarrow \GSp_4(\C)$
	satisfying the following conditions:
	\begin{itemize}
	 \item The restriction to the first factor $\phi\vert_{W_F}\colon W_F\longrightarrow \GSp_4(\C)$
	       is continuous, where $\GSp_4(\C)$ is endowed with the discrete topology.
	       The restriction to the second factor $\phi\vert_{\SL_2(\C)}\colon \SL_2(\C)\longrightarrow \GSp_4(\C)$
	       comes from a homomorphism of algebraic groups.
	 \item For $w\in W_F$, $\phi(w,1)\in \GSp_4(\C)$ is semisimple.
	\end{itemize}
	We denote by $\Phi(G)$ the set of $\GSp_4(\C)$-conjugacy classes of $L$-parameters for $G$.
  \item An $L$-parameter $\phi\colon W_F\times \SL_2(\C)\longrightarrow \GSp_4(\C)$ is said to be
	discrete if $\Imm \phi$ is not contained in any proper parabolic subgroup of $\GSp_4(\C)$.
  \item Let $P$ be a proper parabolic subgroup of $J$, which is unique up to conjugacy.
	It determines a conjugacy class $\widehat{P}$ of parabolic subgroups of $\GSp_4(\C)=\widehat{\GSp}_4$
	(in fact, it is the conjugacy class containing the Siegel parabolic subgroup).
	A proper parabolic subgroup of $\GSp_4(\C)$ belonging to $\widehat{P}$ is said to be relevant for $J$.
	An $L$-parameter $\phi\colon W_F\times \SL_2(\C)\longrightarrow \GSp_4(\C)$ is said to be
	relevant for $J$ if $\Imm \phi$ is not contained in any proper parabolic subgroup of $\GSp_4(\C)$
	which is not relevant for $J$.
	We denote by $\Phi(J)$ the set of $\GSp_4(\C)$-conjugacy classes of $L$-parameters for $G$
	which are relevant for $J$.
 \end{enumerate}
\end{defn}

The main theorem of \cite{MR2800725} can be summarized as follows:

\begin{thm}[Gan-Takeda \cite{MR2800725}]\label{thm:Gan-Takeda}
 Let $\Irr(G)$ denote the set of isomorphic classes of irreducible admissible representations of $G$.
 \begin{enumerate}
  \item There exists a natural surjection $\Irr(G)\longrightarrow \Phi(G)$ with finite fibers.
	For $\phi\in \Phi(G)$, we denote the fiber at $\phi$ by $\Pi_\phi^G$ and call it the $L$-packet
	corresponding to $\phi$.
  \item An irreducible representation $\pi\in\Irr(G)$ is a (essentially) discrete series
	if and only if $\pi\in \Pi^G_\phi$ for a discrete $L$-parameter $\phi\in\Phi(G)$.
  \item For $\phi\in \Phi(G)$, we put $A_\phi=\pi_0(Z_{\GSp_4}(\Imm\phi))$ and write $\widehat{A}_\phi$
	the set of irreducible characters of $A_\phi$. Then, there exists a natural bijection
	between $\Pi^G_\phi$ and $\widehat{A}_\phi$.
  \item For $\phi\in\Phi(G)$ and $\pi\in\Pi^G_\phi$, the central character of $\pi$ is equal to
	$(\simil\circ\phi)\circ\Art$, where $\simil\colon \GSp_4(\C)\longrightarrow \C^\times$
	denotes the similitude character and 
	$\Art\colon F^\times\longrightarrow W_F^{\mathrm{ab}}$ denotes the isomorphism
	of local class field theory (normalized so that a uniformizer is mapped to a lift of
	the geometric Frobenius element). 
  \item For $\phi\in\Phi(G)$ and a smooth character $\chi$ of $F^\times$, let $\phi\otimes\chi$ be
	the $L$-parameter given by
	$(\phi\otimes\chi)(w,v)\longmapsto \chi(\Art^{-1}(w))\phi(w,v)$
	for $(w,v)\in W_F\times\SL_2(\C)$.
	Then we have $\Pi^G_{\phi\otimes\chi}=\{\pi\otimes \chi_G\mid \pi\in\Pi^G_\phi\}$,
	where $\chi_G=\chi\circ\simil$ as in Lemma \ref{lem:RZ-twist}.
  \item There is a way to characterize the map $\Irr(G)\longrightarrow \Phi(G)$
	by means of local factors and Plancherel measures.
 \end{enumerate}
\end{thm}

Moreover, from the construction in \cite{MR2800725}, we have the following:

\begin{thm}\label{thm:sc-packet-G}
 If $\phi\in \Phi(G)$ is discrete and trivial on $\SL_2(\C)$, then $\Pi^G_\phi$ consists of supercuspidal
 representations.
\end{thm}

\begin{prf}
 We freely use the notation in \cite{MR2800725} and \cite{MR2846304}.
 First assume that $\phi$ is not irreducible as a 4-dimensional representation of $W_F$.
 Then $\phi=\phi_1\oplus \phi_2$ where $\phi_1$ and $\phi_2$ are two-dimensional irreducible representations
 of $W_F$ with $\det\phi_1=\det\phi_2$ and $\phi_1\ncong\phi_2$
 (\cf \cite[Lemma 6.2 (ii)]{MR2800725}).
 Let $\tau_i$ be the irreducible representation of $\GL_2(F)$ corresponding to $\phi_i$
 by the local Langlands correspondence. By the assumption, $\tau_1$ and $\tau_2$ are supercuspidal and
 $\tau_1\ncong \tau_2$.
 By the construction of $L$-packets (\cf \cite[\S 7]{MR2800725}), we have
 $\Pi^G_\phi=\{\theta_{(2,2)}(\tau_1\boxtimes\tau_2), \theta_{(4,0)}(\tau^D_1\boxtimes\tau^D_2)\}$,
  where $\theta_{(2,2)}$ (resp.\ $\theta_{(4,0)}$) denotes the theta correspondence between
 $\GSp_4$ and $\mathrm{GSO}_{2,2}$ (resp.\ $\mathrm{GSO}_{4,0}$) and $\tau_i^D$ denotes
 the representation of $D^\times$ corresponding to $\tau_i$ by the local Jacquet-Langlands correspondence.
 \cite[Theorem 8.2 (ii)]{MR2846304} tells us that $\theta_{(2,2)}(\tau_1\boxtimes\tau_2)$ is
 supercuspidal. 
 \cite[Theorem 8.1 (iii)]{MR2846304} tells us that $\theta_{(4,0)}(\tau^D_1\boxtimes\tau^D_2)$ is
 supercuspidal. Therefore $\Pi^G_\phi$ consists of two supercuspidal representations.

 Next assume that $\phi$ is irreducible as a 4-dimensional representation.
 Then, $\Pi_\phi^G$ consists of a single representation $\pi\in\Irr(G)$ such that 
 $\theta_{(3,3)}(\pi)=\Pi\boxtimes\mu$, where $\theta_{(3,3)}$ denotes the theta correspondence 
 between $\GSp_4$ and $\mathrm{GSO}_{3,3}$, $\Pi$ is the irreducible representation of $\GL_4(F)$
 corresponding to the 4-dimensional representation $\phi$, and $\mu=(\simil\circ\phi)\circ\Art$
 is the character of $F^\times$ corresponding to $\simil\circ\phi$.
 By the assumption on $\phi$, $\Pi$ is supercuspidal, and thus $\theta_{(3,3)}(\pi)$ is supercuspidal.
 Therefore, \cite[Section 14, Table 1]{MR2846304} tells us that $\pi$ is supercuspidal, as desired.
\end{prf}

The main theorem of \cite{Gan-Tantono} is as follows:

\begin{thm}[Gan-Tantono \cite{Gan-Tantono}]\label{thm:Gan-Tantono}
  Let $\Irr(J)$ denote the set of isomorphic classes of irreducible admissible representations of $J$.
 \begin{enumerate}
  \item There exists a natural surjection $\Irr(J)\longrightarrow \Phi(J)$ with finite fibers.
	For $\phi\in \Phi(J)$, we denote the fiber at $\phi$ by $\Pi_\phi^J$ and call it the $L$-packet
	corresponding to $\phi$.
  \item An irreducible representation $\rho\in\Irr(G)$ is a (essentially) discrete series
	if and only if $\rho\in \Pi^J_\phi$ for a discrete $L$-parameter $\phi\in\Phi(J)$.
  \item For $\phi\in \Phi(J)$, we put $B_\phi=\pi_0(Z_{\Sp_4}(\Imm\phi))$ and write $\widehat{B}_\phi$
	the set of irreducible characters of $B_\phi$. The natural map $B_\phi\longrightarrow A_\phi$ is
	surjective, and thus $\widehat{A}_\phi$ can be regarded as a subgroup of $\widehat{B}_\phi$.
	There exists a natural bijection
	between $\Pi^J_\phi$ and $\widehat{B}_\phi\setminus \widehat{A}_\phi$.
  \item For $\phi\in\Phi(J)$ and $\rho\in\Pi^J_\phi$, the central character of $\rho$ is equal to
	$(\simil\circ\phi)\circ\Art$. 
  \item For $\phi\in\Phi(J)$, 
	we have $\Pi^J_{\phi\otimes\chi}=\{\rho\otimes \chi_J\mid \rho\in\Pi^J_\phi\}$,
	where $\chi_J=\chi\circ\simil$ as in Lemma \ref{lem:RZ-twist}.
  \item There is a way to characterize the map $\Irr(J)\longrightarrow \Phi(J)$
	by means of local factors and Plancherel measures.
 \end{enumerate}
\end{thm}

We also have an analogous result as Theorem \ref{thm:sc-packet-G}:

\begin{thm}\label{thm:sc-packet-J}
 If $\phi\in \Phi(J)$ is discrete and trivial on $\SL_2(\C)$, then $\Pi^J_\phi$ consists of supercuspidal
 representations.
\end{thm}

\begin{prf}
 We freely use the notation in \cite{Gan-Tantono}.
 First assume that $\phi$ is not irreducible as a 4-dimensional representation of $W_F$.
 Then $\phi=\phi_1\oplus \phi_2$ as in the proof of Theorem \ref{thm:sc-packet-G}.
 Define $\tau_1$ and $\tau_2$ similarly. They are supercuspidal and not isomorphic to each other.
 By the construction of $L$-packets (\cf \cite[\S 7]{Gan-Tantono}), we have
 $\Pi^J_\phi=\{\theta_{(1,1)}(\tau^D_1\boxtimes\tau_2), \theta_{(1,1)}(\tau^D_2\boxtimes\tau_1)\}$,
 where $\theta_{(1,1)}$ denotes the theta correspondence between
 $J=\GSp(1,1)$ and $\mathrm{GO}^*(1,1)$.
 Therefore, \cite[Proposition 5.4 (iv)]{Gan-Tantono} tells us that
 $\Pi^J_\phi$ consists of two supercuspidal representations.

 Next assume that $\phi$ is irreducible as a 4-dimensional representation.
 Then, $\Pi_\phi^J$ consists of a single representation $\rho\in\Irr(J)$ such that 
 $\theta_{(3,0)}(\rho)=\Pi\boxtimes\mu$, where $\theta_{(3,0)}$ denotes the theta correspondence 
 between $J=\GSp(1,1)$ and $\mathrm{GO}^*(3,0)$, $\Pi$ is the irreducible representation of $D_4^\times$
 corresponding to the 4-dimensional representation $\phi$
 (here $D_4$ is the central division algebra over $\Q_p$ with invariant $1/4$),
 and $\mu=(\simil\circ\phi)\circ\Art$
 is the character of $F^\times$ corresponding to $\simil\circ\phi$.
 By the assumption on $\phi$, $\Pi$ is not a character.
 Therefore, \cite[Proposition 5.7]{Gan-Tantono} tells us that $\rho$ is supercuspidal;
 indeed, if $\rho$ is neither supercuspidal nor a twist of the Steinberg representation,
 then $\theta_{(3,0)}(\rho)=0$, and if $\rho$ is a twist of the Steinberg representation,
 then $\theta_{(3,0)}(\rho)$ is a character.
\end{prf}

For an irreducible admissible representation $\pi$ of $G$, we denote the character of $\pi$ by $\theta_\pi$.
It is locally constant on $G^{\mathrm{reg}}$; namely, $\theta_\pi$ is a unique locally constant function
on $G^{\mathrm{reg}}$ such that $\Tr(f;\pi)=\int_G f(g)\theta_\pi(g)dg$ for every $f\in\mathcal{H}(G)$
with $\supp f\subset G^{\mathrm{reg}}$.
For $\phi\in\Phi(G)$, put $\theta_{\Pi^G_\phi}=\sum_{\pi\in \Pi^G_\phi}\theta_\pi$.

Similarly, we define $\theta_{\rho}$ and $\theta_{\Pi^J_\phi}$ for $\rho\in\Irr(J)$ and $\phi\in\Phi(J)$.

\begin{defn}
 \begin{enumerate}
  \item	For $\phi\in\Phi(G)$, the corresponding $L$-packet $\Pi_\phi^G$ is said to be stable
	if $\theta_{\Pi^G_\phi}$ is a stable function on $G^{\mathrm{reg}}$, that is,
	$\theta_{\Pi^G_\phi}(g)=\theta_{\Pi^G_\phi}(g')$ for every $g,g'\in G^{\mathrm{reg}}$ which
	are stably conjugate.
	Similarly we can define the stability of the $L$-packet $\Pi_\phi^J$ for $\phi\in\Phi(J)$.
  \item For $\phi\in\Phi(J)$, we say that $\Pi^G_\phi$ and $\Pi^J_\phi$ satisfy the character relation
	if $\theta_{\Pi^G_\phi}(g)=-\theta_{\Pi^J_\phi}(h)$ for every $g\in G^{\mathrm{ell}}$, 
	$h\in J^{\mathrm{ell}}$ with $g\leftrightarrow h$ (\cf Section \ref{sec:notation}).
 \end{enumerate}
\end{defn}

\begin{rem}
 For every discrete (or more generally, tempered) $L$-parameter $\phi$,
 it is expected that $\Pi^G_\phi$ and $\Pi^J_\phi$ are stable and satisfy the character relation.
 It is plausible that one can deduce these properties from the stable trace formula.

 If $\phi$ is a TRSELP in the sense of \cite{MR2480618}, then the stability and the character relation
 for $\Pi_\phi^G$ and $\Pi^J_\phi$ are already known due to \cite{MR2480618}, \cite{Kaletha-isocrystal}
 and \cite{Lust}.
\end{rem}

\subsection{Computation of the character}

\begin{thm}\label{thm:main}
 Let $\phi\in \Phi(J)$ be an $L$-parameter such that $\Pi^G_\phi$, $\Pi^J_\phi$ are stable
 and satisfy the character relation.
 Then, for every $f\in\mathcal{H}(G)$ with $\supp f\subset G^{\mathrm{ell}}$, we have
 \[
 \sum_{\rho\in\Pi_\phi^J}\Tr(f;H_{\mathrm{RZ}}[\rho])=-4\sum_{\pi\in\Pi^G_\phi}\Tr(f;\pi).
 \]
\end{thm}

\begin{prf}
 By Theorem \ref{thm:Gan-Takeda} iv) and Theorem \ref{thm:Gan-Tantono} iv),
 all representations in $\Pi^G_\phi$ and $\Pi^J_\phi$ share the same central character
 $\omega=(\simil\circ\phi)\circ\Art\colon \Q_p^\times\longrightarrow \overline{\Q}_\ell^\times$.
 First we will reduce the theorem to the case where $\omega\vert_{p^\Z}$ is trivial.
 The method is similar to the proof of Corollary \ref{cor:RZ-adm-vanish}.
 Take $c\in\overline{\Q}^\times_\ell$ such that $c^2=\omega(p)$, and $\chi\colon \Q_p^\times\longrightarrow \overline{\Q}_\ell^\times$ the character given by $\chi(a)=c^{-v_p(a)}$.
 Consider the $L$-packets
 \[
  \Pi^G_{\phi\otimes\chi}=\{\pi\otimes\chi_G\mid \pi\in\Pi^G_\phi\},\quad
 \Pi^J_{\phi\otimes\chi}=\{\rho\otimes\chi_J\mid \rho\in\Pi^J_\phi\}
 \]
  corresponding to $\phi\otimes\chi$
 (\cf Theorem \ref{thm:Gan-Takeda} v), Theorem \ref{thm:Gan-Tantono} v)).
 It is clear that these $L$-packets are stable and satisfy the character relation.
 Moreover, every representation belonging to these $L$-packets has the central character trivial on $p^\Z$.

 We have 
 $H_{\mathrm{RZ}}[\rho\otimes\chi_J]=H_{\mathrm{RZ}}[\rho]\otimes\chi_G$ by Lemma \ref{lem:RZ-twist}.
 Therefore, if we have the theorem for the $L$-parameter $\phi\otimes\chi\in\Phi(J)$, then
 \begin{align*}
  \sum_{\rho\in\Pi_\phi^J}\Tr(f;H_{\mathrm{RZ}}[\rho])
  &=\sum_{\rho\otimes\chi_J\in\Pi_{\phi\otimes\chi}^J}\Tr(f;H_{\mathrm{RZ}}[\rho\otimes\chi_J]\otimes\chi_G^{-1})\\
  &=\sum_{\rho\otimes\chi_J\in\Pi_{\phi\otimes\chi}^J}\Tr(f\cdot \chi_G^{-1};H_{\mathrm{RZ}}[\rho\otimes\chi_J])\\
  &=-4\sum_{\pi\otimes\chi_G\in\Pi_{\phi\otimes\chi}^G}\Tr(f\cdot\chi_G^{-1};\pi\otimes\chi_G)
  =-4\sum_{\pi\in\Pi_\phi^G}\Tr(f;\pi),
 \end{align*}
 and thus the theorem also holds for $\phi$.

 In the following, we assume that the central character $\omega$ is trivial on $p^\Z$.

 By the similar way as in \cite[Lemma 3.5]{LT-LTF}, we can prove that $G^{\mathrm{ell}}$ is contained in
 the union of all open compact-mod-center subgroups of $G$. Therefore, we may assume that $\supp f$ is contained
 in an open compact-mod-center subgroup of $G$. Recall that a maximal open compact-mod-center subgroup
 is conjugate to $N_{\mathscr{L}}$ where $\mathscr{L}$ is one of $\mathscr{L}_0$, $\mathscr{L}_{\mathrm{para}}$
 or $\mathscr{L}_{\mathrm{Siegel}}$ (Lemma \ref{lem:max-cpt-mod-center}).
 Note that both sides of the identity in the theorem do not change
 if we replace $f$ by its conjugate.  Therefore, we may assume that $\supp f$ is contained in
 $N_{\mathscr{L}}$ where $\mathscr{L}$ is one of $\mathscr{L}_0$, $\mathscr{L}_{\mathrm{para}}$
 or $\mathscr{L}_{\mathrm{Siegel}}$. Since $\{K_{\!\mathscr{L},m}\}_{m\ge 1}$ form a fundamental system of
 neighborhoods of $1\in N_{\mathscr{L}}$ consisting of normal subgroups of $N_{\mathscr{L}}$,
 $f$ can be written as a linear combination of $\mathbf{1}_{gK_{\!\mathscr{L},m}}$ with $g\in N_{\mathscr{L}}$.
 Hence we are reduced to the case where $f=\vol(K_{\!\mathscr{L},m})^{-1}\mathbf{1}_{gK_{\!\mathscr{L},m}}$.
 For simplicity, put $K=K_{\!\mathscr{L},m}$. Note that $gK=\supp f\subset G^{\mathrm{ell}}$.

 By Proposition \ref{prop:eta-property}, Theorem \ref{thm:LTF-consequence}, Corollary \ref{cor:LTF-SO} and
 the stable version of Weyl's integration formula, we can compute as follows:
 \begin{align*}
  &\sum_{\rho\in\Pi_\phi^J}\Tr\Bigl(\frac{\mathbf{1}_{gK}}{\vol(K)};H_{\mathrm{RZ}}[\rho]\Bigr)
  =\sum_{\rho\in\Pi_\phi^J}\Tr(g;H_{\mathrm{RZ}}[\rho]^K)\stackrel{(1)}{=}\sum_{\rho\in\Pi_\phi^J}\Tr(\eta_{gK};\rho)\\
  &\qquad=\sum_{\rho\in\Pi_\phi^J}\int_{J/p^\Z}\eta_{gK}(h)\theta_\rho(h)dh
  =\int_{J/p^\Z}\eta_{gK}(h)\theta_{\Pi_\phi^J}(h)dh\\
  &\qquad\stackrel{(2)}{=}\sum_{\{\mathbf{T}'\}_{\mathrm{st}}\in\mathcal{T}^{\mathrm{ell}}_{\mathbf{J},\mathrm{st}}}\frac{1}{\#W_{\mathbf{T}'}(\Q_p)}\int_{T'^{\mathrm{reg}}/p^\Z} D_{\mathbf{J}}(t')\theta_{\Pi_\phi^J}(t')\mathit{SO}_{t'}(\eta_{gK})dt'\\
  &\qquad\stackrel{(3)}{=}4\sum_{\{\mathbf{T}'\}_{\mathrm{st}}\in\mathcal{T}^{\mathrm{ell}}_{\mathbf{J},\mathrm{st}}}\frac{1}{\#W_{\mathbf{T}'}(\Q_p)}\int_{T'^{\mathrm{reg}}/p^\Z} D_{\mathbf{J}}(t')\theta_{\Pi_\phi^J}(t')\mathit{SO}_{g_{t'}}\Bigl(\frac{\mathbf{1}_{gKp^\Z}}{\vol(K)}\Bigr)dt',
 \end{align*}
 where $g_{t'}\in G^{\mathrm{ell}}$ is an arbitrary element with $g_{t'}\leftrightarrow t'$,
 and $D_{\mathbf{J}}(t')$ is the Weyl denominator (\cf \cite[p.~185]{MR700135}). 
 Here (1) follows from Proposition \ref{prop:eta-property} and (3) from Corollary \ref{cor:LTF-SO}.
 (2) is a consequence of Theorem \ref{thm:LTF-consequence} ii), stability of $\theta_{\Pi^J_\phi}$
 and the stable version of Weyl's integration formula for $J$.
 Other equalities are obvious.

 Recall that we have a natural bijection 
 $\mathcal{T}^{\mathrm{ell}}_{\mathbf{J},\mathrm{st}}\yrightarrow{\cong} \mathcal{T}^{\mathrm{ell}}_{\mathbf{G},\mathrm{st}}$. For an elliptic maximal torus $\mathbf{T}'$ of $\mathbf{J}$ and 
 an elliptic maximal torus $\mathbf{T}$ of $\mathbf{G}$ corresponding to $\mathbf{T}$, 
 choose $t'_0\in T'^{\mathrm{reg}}$, $t_0\in T^{\mathrm{reg}}$ with $t_0\leftrightarrow t'_0$ and
 consider the isomorphism $\iota=\iota_{t'_0,t_0}\colon \mathbf{T}'\yrightarrow{\cong}\mathbf{T}$
 (\cf Section \ref{sec:notation}).
 Then we can take $g_{t'}$ as $t=\iota(t')$. Obviously we have $D_{\mathbf{J}}(t')=D_{\mathbf{G}}(t)$.
 Since $\Pi^G_\phi$ and $\Pi^J_\phi$ satisfy the character relation, we have
 $\theta_{\Pi^J_\phi}(t')=-\theta_{\Pi^G_\phi}(t)$. Moreover, $\iota\colon T'/p^\Z\yrightarrow{\cong} T/p^\Z$
 is compatible with the fixed measures (\cf the last paragraph of Section \ref{subsec:def-RZ}).
 Together with Lemma \ref{lem:Weyl-group}, we can convert
 the sum with respect to $\mathcal{T}^{\mathrm{ell}}_{\mathbf{J},\mathrm{st}}$ to
 that with respect to $\mathcal{T}^{\mathrm{ell}}_{\mathbf{G},\mathrm{st}}$, and obtain
 \begin{align*}
  &\sum_{\rho\in\Pi_\phi^J}\Tr\Bigl(\frac{\mathbf{1}_{gK}}{\vol(K)};H_{\mathrm{RZ}}[\rho]\Bigr)\\
  &\qquad=-4\sum_{\{\mathbf{T}\}_{\mathrm{st}}\in\mathcal{T}^{\mathrm{ell}}_{\mathbf{G},\mathrm{st}}}\frac{1}{\#W_{\mathbf{T}}(\Q_p)}\int_{T^{\mathrm{reg}}/p^\Z} D_{\mathbf{G}}(t)\theta_{\Pi_\phi^G}(t)\mathit{SO}_t\Bigl(\frac{\mathbf{1}_{gKp^\Z}}{\vol(K)}\Bigr)dt\\
  &\qquad=\frac{-4}{\vol(K)}\int_{G/p^\Z}\mathbf{1}_{gKp^\Z}(x)\theta_{\Pi_\phi^G}(x)dx
  =\frac{-4}{\vol(K)}\int_G\mathbf{1}_{gK}(x)\theta_{\Pi_\phi^G}(x)dx\\
  &\qquad =-4\sum_{\pi\in\Pi^G_\phi}\Tr\Bigl(\frac{\mathbf{1}_{gK}}{\vol(K)};\pi\Bigr)
 \end{align*}
 (in the second equality, we use the stable version of Weyl's integration formula for $\mathbf{G}$;
 again we use the fact that $\mathbf{1}_{gKp^\Z}$ is supported on $G^{\mathrm{ell}}$).
 This completes the proof.
\end{prf}

\begin{cor}\label{cor:RZ-char}
 Let $\phi$ be an $L$-parameter as in Theorem \ref{thm:main}.
 Assume that for every $\rho\in \Pi^J_\phi$ and integers $i,j\ge 0$, 
 $H^{i,j}_{\mathrm{RZ}}[\rho]$ is a $G$-module of finite length.
 Then, for each $\rho\in \Pi^J_\phi$, we can consider the character 
 $\theta_{H_{\mathrm{RZ}}[\rho]}$ of the virtual representation $H_{\mathrm{RZ}}[\rho]$.
 This character satisfies the following for every $g\in G^{\mathrm{ell}}$:
 \[
 \sum_{\rho\in\Pi^J_\phi}\theta_{H_{\mathrm{RZ}}[\rho]}(g)=-4\sum_{\pi\in\Pi^G_\phi}\theta_\pi(g).
 \]
\end{cor}

\begin{prf}
 This is an immediate consequence of Theorem \ref{thm:main}.
\end{prf}

\begin{rem}
 \begin{enumerate}
  \item Since the $G$-representation $H^{i,j}_{\mathrm{RZ}}[\rho]$ is admissible,
	it has finite length if and only if it is finitely generated as a $G$-module.
  \item At least if $\rho\in\Pi^J_\phi$ is supercuspidal, the assumption in Corollary \ref{cor:RZ-char}
	will be proved once we establish an analogue of Faltings' 
	isomorphism for the Rapoport-Zink tower $\{M_K\}_K$
	(\cf \cite[Lemma 5.2]{LT-LTF}). In a recent preprint of Scholze and Weinstein \cite{Scholze-Weinstein},
	Faltings' isomorphism for a Rapoport-Zink tower of EL type is obtained.
	It is plausible that a similar method is applicable to our case.
 \end{enumerate}
\end{rem}

\begin{cor}
 Under the setting in Corollary \ref{cor:RZ-char}, assume moreover that $\phi$ is discrete and trivial
 on $\SL_2(\C)$.
 Then, for each $\pi\in\Pi^G_\phi$, the representation $\pi^\vee$ of $G$ appears in $H^3_{\mathrm{RZ}}$.
\end{cor}

\begin{prf}
 By the similar method as in the proof of Theorem \ref{thm:main}, we may reduce to the case where
 $\pi$ is trivial on $p^\Z$. 
 By Theorem \ref{thm:sc-packet-G} and Theorem \ref{thm:sc-packet-J}, 
 $\Pi_\phi^G$ and $\Pi_\phi^J$ consist of supercuspidal representations.
 Therefore $\Ext^j_J(H^i_{\mathrm{RZ}},\rho)=0$ if $j\ge 1$, and thus
 $H_{\mathrm{RZ}}[\rho]=\sum_i(-1)^i\Hom_J(H^i_{\mathrm{RZ}},\rho)^{\mathrm{sm}}$.
 We denote the supercuspidal part of $H_{\mathrm{RZ}}[\rho]$ by $H_{\mathrm{RZ}}[\rho]_{\mathrm{cusp}}$.
 Write $\sum_{\rho\in\Pi^J_\phi}H_{\mathrm{RZ}}[\rho]_{\mathrm{cusp}}=\sum_{\pi'}a_{\pi'}\pi'$,
 where $\pi'$ runs through irreducible supercuspidal representations of $G$.
 Assume that $\pi^\vee$ does not appear in $H^3_{\mathrm{RZ}}$. 
 Then, by \cite[Theorem 1.1]{RZ-GSp4}, $\pi^\vee$ can appear only in $H^2_{\mathrm{RZ}}$ and $H^4_{\mathrm{RZ}}$.
 Hence $\pi$ can appear in $\Hom_J(H^i_{\mathrm{RZ}},\rho)^{\mathrm{sm}}$ only if $i=2,4$,
 and thus $a_\pi\ge 0$. 
 On the other hand, by Corollary \ref{cor:RZ-char} and the orthogonality relation of characters, we have
 \[
 a_\pi=\Bigl\langle \theta_\pi,\sum_{\rho\in \Pi^J_\phi}\theta_{H_{\mathrm{RZ}}[\rho]}\Bigr\rangle_{\!\!\mathrm{ell}}
 =-4\Bigl\langle \theta_\pi,\sum_{\pi'\in \Pi^G_\phi}\theta_{\pi'}\Bigr\rangle_{\!\!\mathrm{ell}}=-4
 \]
 (for the definition of $\langle\ ,\ \rangle_{\mathrm{ell}}$, see \cite[\S 5]{LJLC}). This is a contradiction.
\end{prf}

\begin{rem}
 In a forthcoming paper with Tetsushi Ito, the author will give more precise description of
 the cuspidal part $H^i_{\mathrm{RZ}}[\rho]_{\mathrm{cusp}}$ of the individual cohomology
 $H^i_{\mathrm{RZ}}[\rho]$ via global method. 
 We can also obtain information on the action of Weil group on $H^i_{\mathrm{RZ}}[\rho]_{\mathrm{cusp}}$;
 we find the local Langlands correspondence for $G$ and $J$ in $H^3_{\mathrm{RZ}}[\rho]_{\mathrm{cusp}}$,
 as expected.
\end{rem}

\def\cftil#1{\ifmmode\setbox7\hbox{$\accent"5E#1$}\else
  \setbox7\hbox{\accent"5E#1}\penalty 10000\relax\fi\raise 1\ht7
  \hbox{\lower1.15ex\hbox to 1\wd7{\hss\accent"7E\hss}}\penalty 10000
  \hskip-1\wd7\penalty 10000\box7} \def\cprime{$'$} \def\cprime{$'$}
\providecommand{\bysame}{\leavevmode\hbox to3em{\hrulefill}\thinspace}
\providecommand{\MR}{\relax\ifhmode\unskip\space\fi MR }
\providecommand{\MRhref}[2]{%
  \href{http://www.ams.org/mathscinet-getitem?mr=#1}{#2}
}
\providecommand{\href}[2]{#2}

\end{document}